\documentclass[a4paper,11pt,twoside]{article}
\usepackage{pst-fill,pst-grad}   
\usepackage{textcomp} 
\usepackage[english]{babel}  
\usepackage[utf8]{inputenc}   
\usepackage[titletoc]{appendix} 
\usepackage{titlesec} 
\usepackage{graphicx}  
\usepackage{amsmath} 
\usepackage{float} 
\usepackage{fancyhdr}  
\usepackage[matrix,arrow,curve]{xy}
\usepackage{pstricks} 
\usepackage{amsmath,amsfonts,verbatim,afterpage,theorem,euscript,mathrsfs,amssymb}
\usepackage{array}
\usepackage{dsfont} 
\usepackage[colorlinks=true,linkcolor=blue,citecolor=red]{hyperref}
\usepackage{authblk}
\usepackage{color} 

\newtheorem{Proposition}{Proposition}[section]
\newtheorem{PropositionP}{Proposition}
\newtheorem{Lemme}{Lemma}[section]
\newtheorem{Theoreme}{Theorem}[section]
\newtheorem{TheoremeP}{Theorem}
\newtheorem{CorollaireP}{Corollary}

\newtheorem{Remarque}{Remark}

\def \vu{\vec{u}}

\def \R{\mathbb{R}}

\def \finpv{\hfill $\blacksquare$}
\def \pv{{\bf{Proof.}}~} 

\def \ds{\displaystyle}

\title{\bf Sharp well-posedness and  spatial decaying  for a generalized dispersive-dissipative Kuramoto-type equation and applications to related models} 

\author[1]{ Manuel Fernando Cortez \footnote{corresponding author: manuel.cortez@epn.edu.ec}}

\author[2]{ Oscar Jarr\'in\footnote{oscar.jarrin@udla.edu.ec}} 

\affil[1]{\scriptsize Departamento de Matem\'aticas, Escuela Politécnica Nacional,  Ladr\'on de Guevera E11-253, Quito, Ecuador}
\affil[2]{\scriptsize Escuela de Ciencias Físicas y Matemáticas, Universidad de Las Américas, Vía a Nayón, C.P.170124, Quito, Ecuador.}  
\date{\today}
\begin{document} 
\maketitle

\begin{abstract}
We introduce a fairly general dispersive-dissipative nonlinear equation, which is characterized by fractional Laplacian operators in both the dispersive and dissipative terms. This equation includes some physically relevant models of fluid dynamics as particular cases. Among them are the \emph{dispersive Kuramoto-Velarde}, the \emph{Kuramoto-Sivashinsky} equation, and some nonlocal perturbations of the \emph{KdV} and the \emph{Benjamin-Ono} equations. We thoroughly study the effects of the fractional Laplacian operators in the qualitative study of solutions: on the one hand, we prove a sharp well-posedness result in the framework of Sobolev spaces of negative order, and on the other hand, we investigate the pointwise decaying properties of solutions in the spatial variable, which are optimal in some cases. These last results are of particular interest for the corresponding physical models. Precisely, they align with previous numerical works on the spatial decay of a particular kind of solutions, commonly referred to as solitary waves. \\[3mm]  
\textbf{Keywords:} Dispersive-dissipative models in fluid dynamics; Sharp well-posedness; Spatially decaying of solutions; Solitary waves. \\[3mm]
\textbf{AMS Classification:}  35A01, 35B30. 
\end{abstract} 
{\footnotesize \tableofcontents}
\section{Introduction and motivation of the model}
In the context of physical phenomena, the dispersive Kuramoto-Velarde equation
\begin{equation*}
(KV) \quad \partial_t u + \partial^{2}_{x} u +  \partial^{3}_{x} u + \partial^{4}_{x} u  + \gamma_2 \, \partial^{2}_{x}(u^2) + \gamma_3\, (\partial_x u)^2  =0, 
\end{equation*} 
describes slow space-time variations of disturbances at interfaces, diffusion–reaction fronts and plasma instability fronts \cite{Chris-Ve, Gar, Gar1}.
This equation is also applicable to the study of Benard-Marangoni cells, which occur when there is significant surface tension on the interface, especially in a microgravity environment (see \cite{Hym, Norm, Oer}). This situation arises in crystal growth experiments aboard an orbiting space station, although the free interface is metastable concerning  small perturbations. In particular, the nonlinearities,  $\gamma_3\, (\partial_x u)^2$ and  $\gamma_2 \, \partial^{2}_{x}(u^2)$, model pressure destabilization effects striving to rupture the interface. 
Likewise, the equation (KV) is a variation of the Kuramoto–Sivashinsky equation,
\begin{equation*}
(KS) \quad \partial_t u + \partial^{2}_{x} u + \partial^{3}_{x} u  + \partial^{4}_{x} u  + \gamma_3\, (\partial_x u)^2  =0, 
\end{equation*} 
 which describes  slow space-time variations of disturbances at
interfaces, flame fronts, diffusion-reaction fronts,  plasma instability fronts and  the long waves on the interface between two viscous fluids \cite{Hoop}. In this equation,   the linear terms describe a balance
between long-wave instability and short-wave stability, while the nonlinear term provides  a mechanism for energy transfer between wave modes. Finally, remark that  the (KS) equation agrees with the  (KV) equation at $\gamma_2= 0$.

\medskip 

Taking the periodic case into account, the equation (KS) is one of the simplest partial differential equations which is capable of exhibiting chaotic behavior. The long time behavior of the (KS) equation is characterized by the negative (therefore destabilizing) second-order diffusion, the positive (therefore stabilizing) fourth-order dissipation, and the nonlinear coupling term \cite{Pap1}. 

\medskip 

On the other hand, thinking about models describing  the behavior of other types of fluids such as stratified fluids, relevant equations with nonlocal terms appear. A case of this type of equations is, on the one hand,   the \emph{Ostrovsky, Stepanyams and Tsimring} (OST) equation: 
\begin{equation*}
\label{OST} 
 (OST) \quad \partial_t u + \partial^{3}_{x} u +u  \partial_x u + \eta \mathcal{H} (\partial_x u + \partial^{3}_{x} u) = 0,   
\end{equation*}
which is a nonlocal perturbation of the celebrated \emph{Korteweg-de Vries} (KdV) equation and, on the other hand,   a nonlocal perturbed version of the well-known  \emph{Benjamin-Ono}  equation \cite{Benjamin}:  
\begin{equation*}
(npBO) \quad \partial_t u + \mathcal{H} \partial^{2}_{x} u +u  \partial_x u + \eta \mathcal{H} (\partial_x u + \partial^{3}_{x} u) = 0.  
\end{equation*}
where $\mathcal{H}$ is the Hilbert transform (see (\ref{Hilbert-transform}) for a precise definition) and $\eta>0$ is a physical parameter.   

\medskip

The (OST) equation   describes the  radiational  instability of long non-linear waves in a stratified flow caused by internal wave radiation from a shear layer. The parameter $\eta>0$ represents the importance of amplification and damping relative to dispersion. For a more complete physical description, we refer to  \cite{OST0,OST,OST1}. The nonlocal perturbed \emph{Benjamin-Ono}  is a good approximate model for long-crested unidirectional waves at the interface of a two-layer system of incompressible inviscid fluids. Moreover, it gives an analogous model of the (OST) equation in deep stratified fluids \cite{Benjamin}. 

 \medskip
 
One of the main objectives of this article is to introduce a new general  \emph{theoretical} equation that encompasses the aforementioned equations as well as some other \emph{physically relevant}  variants. We also want to understand or at least shed some light on the interaction between the dispersive term and the dissipative term that our equation presents  in relation to the mathematical  questions: the local and  global well-posedness  and  persistence properties of spatially decaying. This latter is of  great interest for the particular  physics  models contained in our equation. Precisely, when compared with previous numerical studies on the spatially decaying of their  solitary waves.  

\medskip 

For the parameters  $\alpha > \beta >0$ and $\gamma_1, \gamma_2, \gamma_3 \in \R$,  we shall consider the following \emph{dispersive-dissipative}, nonlocal and nonlinear equation:  
 \begin{equation}\label{Equation} 
\begin{cases}
\partial_t 	u + \,  D( \partial_x  u)   + \Big( D^{\alpha}_{x} - D^{\beta}_{x} \Big)u  + \gamma_1\,  \partial_x ( u^2)  + \gamma_2 \, \partial^{2}_{x}(u^2) + \gamma_3\, (\partial_x u)^2  =0, \quad (t,x)\in (0,+\infty)\times \R, \\
u(0,\cdot)= u_0. 
\end{cases}	
\end{equation} 
Here, the function $u : [0, +\infty) \times \R \to \R$ denotes the solution, and the function $u_0: \R \to \R$ is the initial datum.

\medskip

The \emph{dispersive} effects are characterized  by the term $\ds{D(\partial_x u)}$, where the operator $D$ is given by  $D=\partial^{2}_{x}$ or  by  $D=\mathcal{H} \partial_x$. In this last expression $\mathcal{H}$ denotes the Hilbert transform, which is a nonlocal operator defined in the Fourier variable as 
\begin{equation}\label{Hilbert-transform}
 \widehat{\mathcal{H}}(\varphi)= - i\,  \text{sing}(\xi)\, \widehat{\varphi}(\xi),   
\end{equation} 
where $\text{sing}(\xi)$ is the sing function and $\varphi \in \mathcal{S}(\R)$.  Therefore, in the Fourier variable we have 
\begin{equation}\label{def-D}
\widehat{D \varphi}(\xi)= m(\xi)\widehat{\varphi}(\xi), \quad \text{where} \quad m(\xi)= \begin{cases}  -\vert \xi \vert^2, \,\, \text{when}\,\, D=\partial^{2}_{x}, \\
\vert \xi \vert, \,\, \text{when}\,\, D=\mathcal{H}\partial_x. \end{cases}
\end{equation}
The whole term $\ds{D(\partial_x u)}$ describes  the linearized dispersion relation in the  equation (\ref{Equation}). 

\medskip

The \emph{dissipative} action of the equation is given by the term $\ds{D^{\alpha}_{x} - D^{\beta}_{x}}$. These  two fractional derivative operators are easily defined in the Fourier variable by the expressions 
\begin{equation}\label{def-Der-frac}
\widehat{D^{\alpha}_{x} \varphi} (\xi)= c_\alpha \vert \xi \vert^\alpha \, \widehat{\varphi} (\xi), \qquad \widehat{D^{\beta}_{x} \varphi} (\xi)= c_\beta \vert \xi \vert^\beta \, \widehat{\varphi} (\xi).
\end{equation}
Thus, the total dissipative action of equation \eqref{Equation}, in terms of Fourier  variable, is essentially  given  is given by  the symbol $\vert \xi \vert^\alpha - \vert \xi \vert^\beta$.

\medskip

From a purely physical perspective, this model is not unreal    since physical phenomena that are purely dissipative or purely dispersive are rarely found. This same fact makes interesting  to study the equation \eqref{Equation} from a mathematical point of view.

\medskip

Finally, the nonlinear part of equation \eqref{Equation} is described by the term  
$\gamma_1\,  \partial_x ( u^2)$,   which represents the classical transport term in fluid models,  and by the terms  $\gamma_2 \, \partial^{2}_{x}(u^2)$,  $\gamma_3\, (\partial_x u)^2$,  taken from the (KV) model introduced above. In particular, these last terms allow the model to have a greater mathematical richness, on the one hand, a blow-up criterion in the well-posedness theory (see  Proposition \ref{Prop-blow-up} below) and, on the other hand,  some optimal spatial  decaying rates of solutions (see Corollary \ref{Corollary-1} below).  

\medskip

One of the main interests of equation (\ref{Equation}) is based on the fact that it contains the following \emph{physically relevant models} as a particular case. This is not an exhaustive  list, but we shall mention the most representative ones.   We shall divide them into two main groups according to the \emph{dispersive effects} of the  term $D(\partial_x)$.     

\medskip

{\bf Nonlocal dispersive effects}. We consider here $D=\mathcal{H}\partial_x$ and then $D(\partial_x)=\mathcal{H}\partial^{2}_{x}$. The nonlocal effects of this term, are given by the Hilbert transform $\mathcal{H}$ (defined in (\ref{Hilbert-transform})). In this group we have the following models. 

\begin{itemize}
\item   By setting $\alpha=3$, $\beta=1$, and $\gamma_2=\gamma_3=0$, the equation (\ref{Equation}) agrees with  the nonlocal perturbed \emph{Benjamin-Ono} equation:
\begin{equation}\label{nonloca-BO}
\partial_t u + \mathcal{H}\partial^{2}_{x} u + \mathcal{H}\Big( \partial^{3}_{x}  u +  \partial_x u  \Big)  +  \partial_x ( u^2)=0.
\end{equation}
This equation is a good approximated model for long-crested unidirectional waves at the interface of a two-layer system of deep stratified incompressible inviscid fluids  \cite{Benjamin}.  
\item When  $\alpha=4$, $\beta=1$, and $\gamma_2=\gamma_3=0$, the equation  (\ref{Equation}) writes down as another relevant physical model: 
\begin{equation}\label{Plasma-equation}
\partial_t u + \mathcal{H}\partial^{2}_{x} u + \partial^{4}_{x} u + \mathcal{H}\partial_x u + \partial_x(u^2)=0. 
\end{equation} 
This equation provides a successful  model in plasma theory \cite{Qian}. 
\item More generally, for $\alpha>\beta>0$,  and $\gamma_2=\gamma_3=0$, the  equation (\ref{Equation}) becomes the following modified \emph{Benjamin-Ono} equation: 
\begin{equation}\label{BO-type-eq}
\partial_t +\mathcal{H} \partial^{2}_{x} u +\Big( D^{\alpha}_{x} - D^{\beta}_{x} \Big)u  +  \partial_x ( u^2)=0.
\end{equation}
This equation was introduced in \cite{Pastran} as a theoretical model to sharply study the well-posedness issues, which are  driven by the parameters $\alpha$ and $\beta$. 
\end{itemize}

\medskip

{\bf Local dispersive effects}. In this case, we consider $D=\partial^{2}_{x}$ (the classical Laplacian operator) and we obtain $D(\partial_x)= \partial^{3}_{x}$. Among the models containing this dispersive term, it is worth mentioning  the following ones.
\begin{itemize}
\item For $\alpha>\beta>0$, and  $\gamma_2=\gamma_3=0$,  the equation (\ref{Equation}) writes down as the following modified \emph{KdV} equation:
\begin{equation}\label{KdV-generalized}
\partial_t u + \partial^{3}_{x} u + \Big( D^{\alpha}_{x} - D^{\beta}_{x} \Big)u  +  \partial_x ( u^2)=0.
\end{equation}
To the best of our knowledge, this equation has not been studied before; and it is a \emph{KdV}-counterpart of the equation (\ref{BO-type-eq}). Precisely, its main  interest is the study of the dispersive effects of the term $\partial^{3}_{x} u$, when compared with the effects of the  dispersive term $\mathcal{H}\partial^{2}_{x} u$ in the equation (\ref{BO-type-eq}). 
 
\item  When   $\alpha=3$, $\beta=1$, and $\gamma_2=\gamma_3=0$ we  have the \emph{OST} equation: 
\begin{equation}\label{OST-equation}
\partial_t u + \partial^{3}_{x} u + \mathcal{H}\Big( \partial^{3}_{x}  u +  \partial_x u  \Big)  +  \partial_x ( u^2)=0,
\end{equation}
which describes  the radiational instability of long non-linear waves in a stratified flow caused by internal wave radiation from a shear layer \cite{OST0,OST,OST1}.   
\item Finally,   when  $\alpha=4$, $\beta=2$,  and $\gamma_1=0$,  the equation  (\ref{Equation}) becomes the dispersive \emph{Kuramoto-Velarde} equation:
\begin{equation}\label{dispersive-Kuramoto-Velarde}
\partial_t u + \partial^{2}_{x} u +  \partial^{3}_{x} u + \partial^{4}_{x} u  + \gamma_2 \, \partial^{2}_{x}(u^2) + \gamma_3\, (\partial_x u)^2  =0, 
\end{equation}
and moreover, when we set $\gamma_2=0$ we obtain the 1D- \emph{Kuramoto-Sivashinsky} equation: 
\begin{equation}\label{Kuramoto-Sivashinsky}
\partial_t u + \partial^{2}_{x} u + \partial^{3}_{x} u  + \partial^{4}_{x} u  + \gamma_3\, (\partial_x u)^2  =0.
\end{equation}
The physics interest of both models was explained above. 
\end{itemize}

As mentioned, the main objective of this paper is to focus  on two relevant  issues for equation \eqref{Equation}: a well-posedness theory in the setting of the Sobolev spaces and the persistence problem of the spatially decaying of solutions.  It is worth emphasizing these qualitative properties deeply depend on the parameters $\alpha, \beta$ in the dissipative term $\Big( D^{\alpha}_{x} - D^{\beta}_{x} \Big)u$, on the parameters $\gamma_2, \gamma_3$ in the nonlinear term $\gamma_2 \partial^{2}_{x}(u^2)+ \gamma_3 (\partial_x u)^2$ as well as on the operator $D$ in the dispersive term $D(\partial_x u)$.

\section{The main results}

\subsection{Well-posedness}
We recall that the equation  (\ref{Equation})  is locally well-posed in the space $H^s(\R)$ (with $s\in \R$)  if for any initial datum $u_0 \in H^s(\R)$ there exists a time $0<T=T(\Vert u_0 \Vert_{H^s})$ and there exists  a unique solution $u(t,x)$ to the equation   (\ref{Equation})  in a space $E_{T} \subset  \mathcal{C}([0,T], H^s(\R))$, such that  the flow-map data-solution: 
\begin{equation}\label{Map-flow}
S: H^s(\R) \to E_{T} \subset  \mathcal{C}([0,T], H^s(\R)),  \quad u_0 \mapsto S(t)u_0= u(t,\cdot),
\end{equation}
is a locally continuous  function from $H^s(\R)$ to $E_T$.  

\medskip

As mentioned,  the local well-posedness (LWP)  of  equation (\ref{Equation})  is driven by the parameter $\alpha$ in its  dissipative term. Precisely,  the constraint $\alpha > 7/2$ will  allow us  to handle the strong nonlinear effects of the  terms $\gamma_2 \partial^{2}_{x}(u^2)+\gamma_3 (\partial_x u)^2$ (see Remark \ref{Rmk1} below for more technical  details on this fact) while in the case $\gamma_2=\gamma_3=0$ this constraint is relaxed to $\alpha>2$. Thus,  our first result states as follows: 
\begin{TheoremeP}[LWP]\label{Th1} Let $\alpha > \beta >0$ and let $\gamma_1,\gamma_2, \gamma_3$ be the parameters in the equation (\ref{Equation}).
\begin{enumerate}
    \item Let $\gamma_2, \gamma_3 \neq 0$. We set $\alpha>7/2$ and then the equation (\ref{Equation}) is locally well-posed in the Sobolev space $H^s(\R)$ with $s > 1 - \alpha/ 2$. Moreover, we have $u \in \mathcal{C}^{1}(]0, T], \mathcal{C}^{\infty}(\R))$ and the flow-map function $S$ defined in (\ref{Map-flow}) is smooth. 
 \item   Let $\gamma_2=\gamma_3=0$. We set $\alpha>2$ and then the equation is locally well-posed in $H^s(\R)$ with $s>\max(3/2-\alpha,-\alpha/2)$. As above,  we have $u \in \mathcal{C}^{1}(]0, T], \mathcal{C}^{\infty}(\R))$ and  the flow-map function $S$  is smooth. 
\end{enumerate}
\end{TheoremeP}	

It is important to emphasize that one of the main  interests of this theorem lies in the understanding of  the relationship between the  parameters $\alpha,\gamma_2$, and $\gamma_3$ with the well-posedness theory for the equation (\ref{Equation}).

\medskip

This theorem also recovers  some known results on the local well-posedness for  the particular physical models introduced above, among them, the nonlocal perturbed \emph{Benjamin-Ono} equation (\ref{nonloca-BO}) studied in \cite{Fonseca}, the plasma model (\ref{Plasma-equation}) and the modified \emph{Benjamin-Ono} equation (\ref{BO-type-eq}) studied in \cite{Pastran}, the \emph{OST} equation (\ref{OST-equation}) investigated in \cite{Wang}, the dispersive \emph{Kuramoto-Velarde} equation  (\ref{dispersive-Kuramoto-Velarde}) investigated in \cite{Pilod} and the $1D-$\emph{Kuramoto-Sivashinsky} (\ref{Kuramoto-Sivashinsky}) equation studied in \cite{Biagioni}.

\medskip 

Compared with these results, the novelty of this theorem is the fact that even for negative values of $s$, $H^s$-initial data yield classical solutions to the equation (\ref{Equation}) since they also belong to the space $\mathcal{C}^{1}(]0, T], \mathcal{C}^{\infty}(\R))$. In particular, for the equations (\ref{dispersive-Kuramoto-Velarde}) and (\ref{Kuramoto-Sivashinsky}) this theorem improves the result obtained in \cite{Coclite}, where the existence of classical solutions is proven for $H^2$-initial data verifying some additional smallness conditions. 

\medskip

Finally, the second  point above provides us  with a new locally well-posedness result  for the modified  \emph{KdV} equation  (\ref{KdV-generalized}), which (to our knowledge) has not been studied before. 

\medskip

On the other hand, we observe that the minimal regularity (measured by the parameter $s$) to prove the local well-posedness in the space $H^s(\R)$ also depends on the parameter $\alpha$ through   the conditions $s>1-\alpha/2$ when $\gamma_2, \gamma_3\neq 0$ and $s>\min(3/2-\alpha,-\alpha/2)$ when $\gamma_2=\gamma_3=0$. In our second result, we prove that  the quantities $1-\alpha/2$ and $-\alpha/2$ are sharp  in the local well-posedness theory in the following sense: 
\begin{TheoremeP}[Sharp LWP]\label{Th2}   
\begin{enumerate}  
    \item[] 
    \item Let $\alpha > \beta >0$ with $\alpha > 7/2$ and $\gamma_2, \gamma_3 \neq 0$.  Let  $s < 1- \alpha /2$.   If  the equation (\ref{Equation}) is locally well-posed in $H^s(\R)$ then the flow-map function $S$ is not a $C^2-$ function at $u_0=0$.
    \item  Let $\alpha > \beta >0$ with $\alpha >2$ and $\gamma_2 = \gamma_3 = 0$.  Let  $s <- \alpha /2$.   If  the equation (\ref{Equation}) is locally well-posed in $H^s(\R)$ then the flow-map function $S$ is not a $C^2-$ function at $u_0=0$.
\end{enumerate}  
\end{TheoremeP}

Next, we are interested in studying the global well-posedness (GWP) of the equation (\ref{Equation}). We recall that this equation is globally well-posed in $H^s(\R)$ if the properties mentioned above hold for any time $0<T$. In our third result, we show that the GWP  is driven by the parameters $ \gamma_2$ and $\gamma_3$ in the nonlinear term $\gamma_2 \partial^{2}_{x}(u^2)+\gamma_3 (\partial_x u)^2$.
\begin{TheoremeP}[GWP]\label{Th3} Within the framework of Theorem \ref{Th1}, the equation (\ref{Equation}) is globally well-posed in the space $H^s(\R)$  (with $s>1-\alpha/2$ or $s>\max(3/2-\alpha,-\alpha/2)$) when $-2 \gamma_2 + \gamma_3=0$.
\end{TheoremeP}	
In the particular case when $\gamma_2=\gamma_3=0$, we recover the GWP for the set of models (\ref{nonloca-BO}), (\ref{Plasma-equation}), (\ref{BO-type-eq}), (\ref{OST-equation}) (see  the  references mentioned above). Moreover, we give a new GWP result for the \emph{modified KdV} equation (\ref{KdV-generalized}). 

\medskip

Concerning the \emph{Kuramoto-Sivashinsky} equation (\ref{Kuramoto-Sivashinsky}) (where $\gamma_2=0$) the constraint $-2 \gamma_2 + \gamma_3=0$ implies  that $\gamma_3=0$, and consequently, this result trivially holds  for the  linear version of this equation. It is worth emphasizing  this fact is coherent with \cite{Pokhozhaev}, where it is shown that the nonlinear term $\gamma_3 (\partial_x u)^2$ yields a finite blow-up of solutions to the equation (\ref{Kuramoto-Sivashinsky}) associated with a large class of initial data.  

\medskip 

 Finally, for the \emph{dispersive Kuramoto-Velarde} equation (\ref{dispersive-Kuramoto-Velarde}), we are able to  ensure its GWP as long as $-2 \gamma_2 + \gamma_3=0$, which was pointed out in \cite{Pilod}.  The GWP or blow-up phenom in the case  $-2 \gamma_2 + \gamma_3\neq0$ remains an  open question far from obvious and, in future research, we aim to give a deeper understanding of the effects of the nonlinear term  $\gamma_2 \partial^{2}_{x}(u^2)+\gamma_3 (\partial_x u)^2$ in the GWP theory.  However, by performing some new energy estimates we  can  prove the following: 
\begin{PropositionP}[Blow-up criterion]\label{Prop-blow-up} Within the framework of Theorem \ref{Th1}, assume that  $\gamma_2$ and $\gamma_3$ are such that $ -2\gamma_2 + \gamma_3\neq0$. Then,  for a time $0<T^*<+\infty$ we have:
\[ \lim_{t\to T^{*}}\Vert u(t,\cdot)\Vert_{H^s}=+\infty \,\,  \qquad \text{if and only if} \qquad \int_{0}^{T^{*}} \Vert \partial^{2}_{x} u(t,\cdot) \Vert_{L^{\infty}} dt = +\infty. \]   	 
\end{PropositionP}	
This result gives us a new blow-up criterion for the \emph{dispersive Kuramoto-Velarde} equation (\ref{dispersive-Kuramoto-Velarde})  and its related models containing the nonlinear term $\gamma_2 \partial^{2}_{x}(u^2)+\gamma_3 (\partial_x u)^2$. 

\subsection{Spatially  decaying} 
In this section we study another relevant qualitative property of  equation (\ref{Equation}): the pointwise decaying  of solutions $u(t,x)$ for to the spatial variable $x$.  This question gives us a good comprehension of the terms in this equation governing the spatial behavior of solutions, and it is also of physical interest when particularizing in the  models introduced above. Specifically, when comparing with the spatial behavior of a relevant  kind of particular  solutions, the so-called \emph{solitary waves}.   

\medskip 

From the nonlinear differential equations point of view, the existence of the \emph{solitary
wave} describes a perfect balance between the nonlinearity and the dispersive character of its linear part. We refer to the book \cite{LinaresPonce} for more details.  Concerning the physics models introduced above, there exist  previous \emph{numerical works}  on the spatial decaying of solitary waves. These works give some light on the spatial decaying of solutions to these equations. In this context, the main contribution of this work is to use the general framework of  equation (\ref{Equation}) to \emph{analytically} study the spatially decaying of solutions, which simultaneously holds  for the particular  physics models contained in this equation.   

\medskip

As we shall observe, our main remark is that these decaying properties of solutions to the equation (\ref{Equation}) are driven by both the dispersive term $D(\partial_x u)$ and the dissipative term $(D^{\alpha}_{x}-D^{\beta}_{x})u$. Precisely,  by the operator $D$ defined in (\ref{def-D}) and the parameters $\alpha$ and $\beta$.  We introduce here the parameter $n \geq 2$, which depends on $D$, $\alpha$, and $\beta$, as follows:
\begin{equation}\label{parameter-n}
n=\begin{cases} \vspace{2mm}
\min(3, [\beta]+1) , \quad \text{when $D=\mathcal{H}\partial_x$}, \\ \vspace{2mm}
\text{any natural number}, \quad \text{when $D=-\partial^{2}_{2}$ and $\alpha, \beta$ are both even numbers}, \\ \vspace{2mm}
[\alpha]+1, \quad \text{when $D=-\partial^{2}_{x}$, $\beta$ is an even number and not $\alpha$},\\
[\beta]+1, \quad \text{when $D=-\partial^{2}_{x}$ and $\beta$ is not an even number},
\end{cases}
\end{equation}
where $[\alpha]$ and $[\beta]$ denote the integer part of $\alpha$ and $\beta$ respectively. 
The parameter $n$ gives us a detailed description of the \emph{pointwise  decaying rate} of solutions to the equation (\ref{Equation}), and our next result reads as follows: 
\begin{TheoremeP}[Spatially pointwise decaying]\label{Th-Decaying} Let $s>\frac{5}{2}$ and let $u_0 \in H^s(\R)$ be an initial datum. Let $\alpha>\beta \geq 1$, with $\alpha > 7/2$ when $\gamma_2,\gamma_3\neq 0$ and $\alpha> 2$ when $\gamma_2=\gamma_3=0$.  Moreover, let $u \in \mathcal{C}([0,T],H^s(\R))$ be the solution to the equation (\ref{Equation}) associated with $u_0$, given by Theorem \ref{Th1}.  

\medskip

Let $\kappa>1$ and assume that the initial datum $u_0$ verifies
\begin{equation}\label{Decaying-Data}
| u_0(x)| \leq \frac{c_0}{1+| x |^\kappa}, \quad  x\in \R,  \end{equation}
 with a constant  $c_0>0$. Then  the solution $u(t,x)$ verifies the following pointwise estimate
\begin{equation}\label{Decaying-Solution}
| u(t,x)| \leq \frac{c_1(t,u)}{1+| x |^{\min(\kappa,n)}}, \quad 0<t<T, \ \ x \in \R,   \end{equation}
with  a constant $c_1(t,u)>0$ depending on $t$, $u$, and $c_0>0$; and where the parameter $n\geq 2$ is defined in (\ref{parameter-n}). 
\end{TheoremeP}

Let us make the following comments. The assumption of the initial data $u_0 \in H^s(\R)$ with $s>5/2$ ensures that the arising solution verifies $u(t,\cdot)\in H^s(\R)$ for all $t\geq 0$ (see Theorem \ref{Th1}). In particular, the technical constraint $s>5/2$ allows us to handle the nonlinear terms in  equation (\ref{Equation}). 

\medskip 

In expressions (\ref{Decaying-Data}) and (\ref{Decaying-Solution}), we may observe that the parameter $n$ controls the decaying properties of solutions:  the solution $u(t,x)$ fulfills the decaying given by the initial datum only if $\kappa \leq  n$. But, for initial data decaying fast enough ($\kappa > n$) the corresponding solution does not mimic this decaying rate and it decays at infinity like $ 1/ |x|^n$.   

\medskip

From now on, we shall assume initial data decaying fast enough: $\kappa > n$, and we shall  discuss more in detail the decaying  estimate verified by the solution $| u(t,x)| \lesssim 1/|x|^n$. To do this,  recall that  the parameter $n$ ultimately depends on the operator $D$ and the parameters $\alpha, \beta$ according to the expression (\ref{parameter-n}). 
\begin{itemize}
    \item When $D=\mathcal{H}\partial_x$  the nonlocal effects of this operator have a strong influence on the spatially decaying properties of solutions. Precisely, in this case, we have $n=\max(3,[\beta]+1)$, where the number $3$ is due to the presence of  the Hilbert transform $\mathcal{H}$. See Proposition \ref{Prop-Kernel-1} below for more details. Consequently, the physical models containing  the  dispersive term $\mathcal{H} \partial^{2}_{x} u$  verify the estimate 
   \begin{equation}\label{Decaying-BO}
      |u(t,x)| \lesssim  \frac{1}{|x|^{\max(3,[\beta]+1)}}.
      \end{equation} 
In particular, solutions  to the \emph{plasma model} (\ref{Plasma-equation})  and solutions to  the  nonlocal perturbed \emph{Benjamin-Ono} equation (\ref{nonloca-BO})  (in both cases we have $\beta=1$) have the spatially decaying $\ds{| u(t,x)| \lesssim \frac{1}{|x|^2}}$. For this last equation, this information is coherent with \cite{Bao}, where the authors numerically prove that their \emph{solitary waves}   behave at infinity as $1/|x|^2$. On the other hand, the spatially decaying properties of  solutions to  the modified \emph{Benjamin-Ono} equation (\ref{BO-type-eq}) have not been studied before; and they satisfy  the spatially decaying  (\ref{Decaying-BO}). Here,  we also realize the effects of  parameter $\beta$ in the dissipative term, while parameter $\alpha$ does not intervene. 
 
 \item When $D=-\partial^{2}_{x}$ it is interesting to observe that the local effects of this operator do not influence the decaying properties of solutions, which are now driven by the parameters $\alpha$ and $\beta$. Here, we have the following cases. 
 \begin{itemize}
     \item When $\alpha$ and $\beta$ are both even numbers, one can set any parameter $n\in \mathbb{N}$ (with $n\geq 2$); and for initial data verifying (\ref{Decaying-Data}) with $\kappa>N$, solutions to the equation (\ref{Equation}) verify the estimate
 \[ |u(t,x)| \lesssim \frac{1}{|x|^n}. \]
 This \emph{persistence problem} is verified for the \emph{dispersive Kuramoto-Velarde} equation (\ref{dispersive-Kuramoto-Velarde}) and the \emph{Kuramoto-Sivashinsky} equation (\ref{Kuramoto-Sivashinsky}). Moreover,  this fact is in concordance with some numerical studies  on the well-localized solitary waves to these equations. See for instance \cite{Chris-Ve} and \cite{Chris-Ve-2}. 
 \item  When $\beta$ is an even number but not $\alpha$, solutions to the equation (\ref{Equation}) have a decaying rate 
 \[ |u(t,x)| \lesssim \frac{1}{|x|^{[\alpha]+1}}, \]
  while $\beta$ is not an even number it holds
 \[ |u(t,x)| \lesssim  \frac{1}{|x|^{[\beta]+1}}.\] 
These decaying rates are verified by the modified \emph{KdV} equation (\ref{KdV-generalized}) according to these cases of the parameters $\alpha$ and $\beta$. Moreover, it is interesting to get back to the modified \emph{Benjamin-Ono} equation (\ref{BO-type-eq}), which verifies the decaying rate (\ref{Decaying-BO}), to highlight  the stronger effects of the  dissipative term $\mathcal{H}\partial^{2}_{x}u$ compared with  the dissipative term $-\partial^{3}_{x} u$.

\medskip

\item Finally, for $\beta=1$, solutions to the \emph{OST} equation (\ref{OST-equation}) verify the decaying rate  
\[ |u(t,x)| \lesssim \frac{1}{|x|^2}.\]
This spatially decaying   agrees with numerical studies performed on \emph{solitary waves} in \cite{Bao-2}; and it was also analytically proven in our previous work \cite{CorJar0}.
\end{itemize}
\end{itemize} 

\medskip

Now, we are interested in studying the optimality of the decaying rate (\ref{Decaying-Solution}) (with $\kappa>n$). To do this, first,  we shall  find an asymptotic profile for the solution $u(t,x)$ to the equation (\ref{Equation}). In order to state our next theorem, we need to introduce function $K_\alpha, \beta(t,x)$ which is obtained as the solution of the linear problem (when $\gamma_1=\gamma_2=\gamma_3=0$) of the equation (\ref{Equation}):   
\begin{equation}\label{Equation-Kernel}
\begin{cases}
\partial_t K_{\alpha,\beta} +  D(\partial_x K_{\alpha,\beta})+ \Big( D^{\alpha}_{x} - D^{\beta}_{x} \Big) K_{\alpha,\beta}=0, \qquad (t,x) \in (0,+\infty)\times \R,\\
K_{\alpha,\beta}(0,\cdot)=\delta_0, 
\end{cases}
\end{equation} 
where $\delta_0$ denotes the Dirac mass at the origin.  It is thus interesting to observe that the asymptotic profile of the solution $u(t,x)$ to the equation (\ref{Equation}) is essentially given by the function $K_{\alpha,\beta}(t,x)$. Precisely, we start by studying the pointwise decaying (in the spatial variable) of this  function. 
\begin{PropositionP}\label{Prop-Principal-Kernel} Let $n \geq 2$ be the parameter given in the expression (\ref{parameter-n}). Moreover, let $\alpha>\beta \geq 1$, with $\alpha > 2$.  For all $t>0$ fixed, there exists a quantity $I(t)$, which verifies $|I(t)| \leq C_1 e^{\eta_1\, t}$ with constants $C_1, \eta_1 >0$ depending on $\alpha$ and $\beta$, such that the following identity holds:
\begin{equation}\label{Kernel-pointwise}
\vert K_{\alpha,\beta}(t,x)\vert = \frac{\vert I(t)\vert}{\vert x \vert^n}, \quad t>0, \quad  x\neq 0.
\end{equation}
\end{PropositionP}

\medskip 

Then, our next theorem writes down as follows: 
\begin{TheoremeP}[Asymptotic profile]\label{Th-Asymptotic-Profile} Let $u_0 \in H^s(\R)$ (with $s>\frac{5}{2}$) be an initial datum  verifying (\ref{Decaying-Data}) with $\kappa>n$; and where the parameter $n\geq 2$ is defined in (\ref{parameter-n}). Let $u \in \mathcal{C}([0,T], H^s(\R))$ (with $s>\frac{5}{2}$) be the associated solution to the equation (\ref{Equation}) given by Theorem \ref{Th1}.

\medskip

For $0<t\leq T$ fixed, this solution has  the following asymptotic development in the spatial variable 
\begin{equation}\label{Asymptotic}
u(t,x) =  K_{\alpha,\beta}(t,x) \left( \int_{\R} u_0 (y)dy\right) + \gamma_3 \,  K_{\alpha,\beta}(t,x)\int_{0}^{t} \left( \int_{\R} (\partial_x u)^2(\tau,y) dy \right) d \tau  + R(t,x),  \quad \vert x \vert \to +\infty,   
\end{equation}  
where 
\begin{equation}
| R(t,x)| \leq \frac{c_2(t,u)}{| x |^{n+\varepsilon}}, \quad 0<\varepsilon \leq 1,
\end{equation}
with a constant  $c_2(t,u)>0$ depending on $t$ and $u$. 
\end{TheoremeP} 

From this asymptotic development,  we can deduce some   optimally decaying properties of the solution $u(t,x)$. First, the expression $\ds{\gamma_3 \,  K_{\alpha,\beta}(t,x)\int_{0}^{t} \left( \int_{\R} (\partial_x u)^2(\tau,y) dt \right) d \tau}$ highlights interesting effects of the nonlinear $\gamma_3 (\partial_x u)^2$ in the spatially decaying of solutions. Precisely, when $\gamma_3\neq 0$ this expression yields the following estimate from below:
\begin{CorollaireP}\label{Corollary-1} Within the framework of Theorem \ref{Th-Asymptotic-Profile}, assume that $\gamma_3\neq 0$. Then the solution $u(t,x)$ to the  equation (\ref{Equation}) verifies:
\begin{equation}\label{DecayOp}
\frac{c_3(u_0, \gamma_3, t, u)}{|x|^n} \leq \vert u (t,x)\vert,\quad \vert x \vert \to +\infty,  
\end{equation}
where the quantity $c_3(u_0, \gamma_3, t, u)>0$ (given in (\ref{c3})) depends on $u_0, \gamma_3,t$, and $u$ but it is independent of  the variable $x$. 
\end{CorollaireP}
Consequently, the physical models containing the nonlinear term $\gamma_3 (\partial_x u)^2$ have an \emph{optimal decaying rate}:
\begin{equation}\label{Optimal-Decaying}
 |u(t,x)| \sim  \frac{1}{|x|^n}, \qquad |x|\to +\infty.    
\end{equation}
In particular, this optimal decaying rate is verified by the \emph{dispersive Kuramoto-Velarde} equation (\ref{dispersive-Kuramoto-Velarde}) and the \emph{Kuramoto-Sivashinsky} equation (\ref{Kuramoto-Sivashinsky}).

\medskip 

We study now the case when $\gamma_3 = 0$. Here, solutions to the equation (\ref{Equation}) have the asymptotic profile:
\[ u(t,x) =  K_{\alpha,\beta}(t,x) \left( \int_{\R} u_0 (y)dy\right)  + R(t,x),  \quad \vert x \vert \to +\infty, \]
where the optimality properties are now driven by the term $\ds{\int_{\R} u_0 (y) dy }$. 
\begin{CorollaireP}\label{Corollary-2} Within the framework of Theorem \ref{Th-Asymptotic-Profile},  assume that $\gamma_3= 0$. In this case, we have the following scenarios: 
\begin{enumerate}
    \item If the initial datum $u_0$ verifies $\ds{\int_{\R} u_0 (y) dy \neq  0}$, 
then  the solution $u(t,x)$ to the equation (\ref{Equation}) verifies the estimate from below:
\begin{equation}\label{DecayOp}
 \frac{c_4(u_0,t)}{\vert x \vert^{n}} \leq \vert u (t,x)\vert,\quad \vert x \vert \to +\infty,    
\end{equation}
where the quantity  $c_4(u_0,t)>0$ (given in (\ref{c4})) depends on $u_0$ and $t$ but it is independent of $x$. 
\item Otherwise, if the initial data verifies $\ds{\int_{\R} u_0 (y) dy = 0}$, 
then the solution $u(t,x)$ to the equation (\ref{Equation}) verifies the estimate from above:
\begin{equation}\label{Decay1Imp}
 \vert u (t,x)\vert \leq \frac{c_2(t,u)}{\vert x \vert^{n+\varepsilon}},\quad 0<\varepsilon\leq 1, \quad  \vert x \vert \to +\infty.
\end{equation}  
\end{enumerate}
\end{CorollaireP}

In the first point above, we obtain an \emph{optimal decaying rate} of solutions (\ref{Optimal-Decaying}) as long as $\int_{\R} u_0 (y) dy \neq  0$. In particular,  this property is verified by the physical models from  equation (\ref{nonloca-BO}) to  equation (\ref{OST-equation}), with the respective values of the parameter $n$ detailed above. 

\medskip 

On the other hand, the second point above shows us that this decaying rate can be improved to $|u(t,x)| \lesssim 1/ |x|^{n+\varepsilon}$ (with $0<\varepsilon\leq 1$) in the case of zero-mean initial data. To the best of our knowledge, the value $\varepsilon =1$ seems to be the maximal one to improve the  decaying rate. Since solutions of equation (\ref{Equation}) are written in an explicit \emph{mild formulation} (\ref{Equation-Integral}) involving the function $K_{\alpha,\beta}(t,x)$ defined above. So, the sharp spatially decaying  properties of this function  (given in Proposition \ref{Prop-Principal-Kernel} above) eventually  block an improvement in the decaying of the solution for $\varepsilon >1$.  

\medskip

{\bf Notation.}  To get rid of some unsubstantial constants, for $A,B >0$,  the notation $A \lesssim B$ means that $A \leq c B$ with a constant $c>0$ which does not depend on $A$ nor $B$. Similarly, we shall write $A \sim B$ when $c_1 A \leq B \leq c_2 B$. On the other hand, the Fourier transform (in the spatial variable) of a function $f$ is denoted by $\widehat{f}$ or $\mathcal{F}(f)$, while $\mathcal{F}^{-1}(f)$ stands for the inverse Fourier transform. 

\medskip

{\bf Organization of the paper.} This paper is divided into two big sections: in Section \ref{Sec-Well-posedness} we give a proof of all the results stated in the well-posedness theory, while Section \ref{Sec-Spatially-decaying} is devoted to proving all the results stated in the spatially decaying theory.

\section{The well-posedness theory}\label{Sec-Well-posedness}
\subsection{Kernel estimates I} 
Let $K_{\alpha,\beta}(t,x)$ be the solution to the linear problem (\ref{Equation-Kernel}).  Then, by definition of the operators  $D$, $D^{\alpha}_{x}$ and $D^{\beta}_{x}$, given in the formulas (\ref{def-D}) and (\ref{def-Der-frac}) respectively, for all $t>0$ we have 
\begin{equation}\label{kernel}
\begin{split}
K_{\alpha,\beta}(t,x)= &\,  \mathcal{F}^{-1}\Big( e^{- ( i\, m(\xi)  \xi  +  (\vert \xi \vert^\alpha - \vert \xi \vert^\beta)) t } \Big)(x)\\
=&\, \mathcal{F}^{-1}\left( e^{-f(\xi)t}\right)(x), \end{split}
 \qquad   \text{with}\ \ f(\xi)= i\, m(\xi)\xi + (\vert \xi \vert^\alpha - \vert \xi \vert^\beta).  
\end{equation}
 In what follows we summarize some  properties of the kernel $K_{\alpha,\beta}$, which will be useful in the sequel.   

\begin{Lemme}\label{Key-Lemma} Let $\sigma \geq 0$. For all $t>0$ we have $\ds{ \left\Vert \left\vert  t^{\frac{1}{\alpha}} \,  \xi \right\vert^{2 \sigma} \widehat{K_{\alpha,\beta}}(t, \cdot) \right\Vert_{L^{\infty}} \lesssim    t^{\frac{2 \sigma}{\alpha}}\, e^{ t} +t^{\frac{2 \sigma}{\alpha}} + 1}$.
\end{Lemme}	
\pv 
By (\ref{kernel}) for all $\xi \in \R$ we have $\ds{\left\vert \widehat{K_{\alpha,\beta}}(t,\xi) \right\vert = e^{-(\vert \xi \vert^\alpha -\vert \xi \vert^\beta) t}}$.  Then, we write
\begin{equation}
\begin{split}
\left\Vert \left\vert  t^{\frac{1}{\alpha}} \,  \xi \right\vert^{2 \sigma} \widehat{K_{\alpha,\beta}}(t, \cdot) \right\Vert_{L^{\infty}} \leq &   \left\Vert \left\vert  t^{\frac{1}{\alpha}} \,  \xi \right\vert^{2 \sigma} e^{- (\vert \xi \vert^\alpha -\vert \xi \vert^\beta) t}\right\Vert_{L^{\infty}\left(\vert \xi \vert \leq 2^{\frac{1}{\alpha-\beta}}\right)}+  \left\Vert \left\vert  t^{\frac{1}{\alpha}} \,  \xi \right\vert^{2 \sigma}e^{- (\vert \xi \vert^\alpha -\vert \xi \vert^\beta) t} \right\Vert_{L^{\infty}\left(\vert \xi \vert > 2^{\frac{1}{\alpha-\beta}}\right)}\\
=& I_1 + I_2. 
\end{split}
\end{equation} 
We observe that the term $I_1$ above can be split as: 
\begin{equation*}
I_1 \leq \left\Vert \left\vert  t^{\frac{1}{\alpha}} \,  \xi \right\vert^{2 \sigma} e^{- (\vert \xi \vert^\alpha -\vert \xi \vert^\beta) t} \right\Vert_{L^{\infty}\left(\vert \xi \vert \leq 1\right)}+ \left\Vert \left\vert  t^{\frac{1}{\alpha}} \,  \xi \right\vert^{2 \sigma} e^{- (\vert \xi \vert^\alpha -\vert \xi \vert^\beta) t}\right\Vert_{L^{\infty}\left(1<\vert \xi \vert \leq 2^{\frac{1}{\alpha-\beta}}\right)}=I_{1,a}+I_{1,b}. 
\end{equation*}
Here, to estimate the term $I_{1,a}$,  since $\vert \xi \vert \leq 1$  we write $- (\vert \xi \vert^\alpha -\vert \xi \vert^\beta) t = (\vert \xi \vert^\beta - \vert \xi \vert^\alpha) t \leq \vert \xi \vert^\beta t \leq  t$.  We thus get  $\ds{I_{1,a} \lesssim t^{\frac{2\sigma}{\alpha}} e^{ t}}$. Similarly, to estimate the  $I_{1,b}$, since $1 < \vert \xi \vert \leq 2^{\frac{1}{\alpha-\beta}}$   we have  $-(\vert \xi \vert^\alpha - \vert \xi \vert^\beta) <0$ and  we obtain  $\ds{I_{1,b} \lesssim t^{\frac{2\sigma}{\alpha}}}$. On the other hand, to estimate the term $I_2$, since  $\vert \xi \vert >2^{\frac{1}{\alpha-\beta}}$ we get $-(\vert \xi \vert^\alpha - \vert \xi \vert^\beta) \leq -\frac{1}{2} \vert \xi \vert^\alpha$ and we can write
\[ I_2 \lesssim  \left\Vert \left\vert  t^{\frac{1}{\alpha}} \,  \xi \right\vert^{2 \sigma}e^{- \vert \xi \vert^\alpha \,  t} \right\Vert_{L^{\infty}\left(\vert \xi \vert > 2^{\frac{1}{\alpha-\beta}}\right)}  \lesssim \left\Vert \left\vert  t^{\frac{1}{\alpha}} \,  \xi \right\vert^{2 \sigma}e^{- \vert t^{\frac{1}{\alpha}} \xi \vert^\alpha } \right\Vert_{L^{\infty}\left(\R\right)}\lesssim 1. \]
By gathering the estimates on the terms $I_{1,a}$, $I_{1,b}$ and $I_2$ we obtain the wished result.   \finpv

\medskip

With this estimate, we can prove the following result. 

\begin{Lemme}\label{Lem-Ker-1} Let $s \in \R$ and let  $s_1\geq 0$.  There exists a constant $\eta_0>0$, which depends on $s$ and $\alpha$,  such that the following estimate holds $\ds{ \Vert K_{\alpha, \beta}(t,\cdot)\ast \psi \Vert_{H^{s+s_1}} \lesssim \frac{e^{\eta_0\, t }}{t^{\frac{s_1}{\alpha}}}\,  \Vert  \psi \Vert_{H^s} }$.
\end{Lemme}	 
\pv  By the H\"older inequalities   we write:
\begin{equation*}
\left\Vert K_{\alpha,\beta}(t,\cdot)\ast \psi \right\Vert_{H^{s+s_1}}= \left\Vert (1+\vert \xi \vert^2)^{\frac{s+s_1}{2}} \widehat{K_{\alpha,\beta}}(t,\cdot) \, \widehat{u_0} \right\Vert_{L^2} \leq  \left\Vert  (1+\vert \xi \vert^2)^{\frac{s_1}{2}}\widehat{K_{\alpha,\beta}}(t,\cdot) \right\Vert_{L^\infty} \, \left\Vert (1+\vert \xi \vert^2)^{\frac{s}{2}} \widehat{\psi} \right\Vert_{L^2},
\end{equation*}
where we must estimate the quantity $\ds{\left\Vert  (1+\vert \xi \vert)^{\frac{s_1}{2}}\widehat{K_{\alpha,\beta}}(t,\cdot) \right\Vert_{L^\infty}}$. We thus have  
\begin{equation*}
\left\Vert  (1+\vert \xi \vert^2)^{\frac{s_1}{2}}\widehat{K_{\alpha,\beta}}(t,\cdot) \right\Vert_{L^\infty} \lesssim \left\Vert \widehat{K_{\alpha,\beta}}(t,\cdot) \right\Vert_{L^\infty}+ \left\Vert  \vert \xi \vert^{s_1} \widehat{K_{\alpha,\beta}}(t,\cdot) \right\Vert_{L^\infty}.
\end{equation*}
Now, to estimate the first term on the right-hand side,  by  Lemma \ref{Key-Lemma} (with $\sigma =0$)   we obtain 
\begin{equation*}
\left\Vert \widehat{K_{\alpha,\beta}}(t,\cdot) \right\Vert_{L^\infty} \lesssim e^t + 1.
\end{equation*}
For the second term on the right-hand  side,  we use again the Lemma \ref{Key-Lemma}  (with $\sigma = \frac{s_1}{2}$) to get
\begin{equation*}
\left\Vert  \vert \xi \vert^{s_1} \widehat{K_{\alpha,\beta}}(t,\cdot) \right\Vert_{L^\infty}= \frac{1}{t^{\frac{s_1}{\alpha}}}\, \left\Vert  \left\vert t^{\frac{1}{\alpha}} \, \xi \right\vert^{s_1} \widehat{K_{\alpha,\beta}}(t,\cdot) \right\Vert_{L^\infty} \lesssim e^t + 1 + \frac{1}{t^{\frac{s_1}{\alpha}}}. 
\end{equation*}
By gathering these estimates, and by setting a quantity $\eta_0 = \eta_(s_1, \alpha)>0$ big enough,  we finally obtain
\begin{equation*}
\left\Vert  (1+\vert \xi \vert^2)^{\frac{s_1}{2}}\widehat{K_{\alpha,\beta}}(t,\cdot) \right\Vert_{L^\infty} \lesssim e^t + 1 + \frac{1}{t^{\frac{s_1}{\alpha}}} = \frac{e^t\, t^{\frac{s_1}{\alpha}} + t^{\frac{s_1}{\alpha}} + 1}{t^{\frac{s_1}{\alpha}}}\lesssim  \frac{e^{\eta_0\, t }}{t^{\frac{s_1}{\alpha}}}. 
\end{equation*}
\finpv

\medskip

Finally, we state our last technical lemma. The proof  essentially follows the  same computations performed in  \cite[Lemma $4.1$]{CorJar1}. 
\begin{Lemme}\label{Lem-Ker-2} Let $s \in \R$, $\delta \geq 0$ and let  $\varepsilon>0$. Then, there exists a constant $C=C(s,\varepsilon,\delta)>0$,   such that for all $\varepsilon< t_1, t_2 \leq T$ we have: 
	\[ \Vert K_{\alpha,\beta}(t_1, \cdot)\ast \psi - K_{\alpha,\beta}(t_2, \cdot)\ast \psi \Vert_{H^{s+\delta}}  \leq C\,   \vert t_1 - t_2 \vert^{1/2} \Vert \psi \Vert_{H^s}.\]	  
\end{Lemme}	
\subsection{Sharp local well-posedness}
\subsubsection*{Proof of  Theorem \ref{Th1}}
\subsubsection*{The case $\gamma_2, \gamma_3 \neq 0$ and $\alpha>7/2$.} 
We  divide the proof into three main steps, which we will prove in the technical theorems below. Precisely, in Theorem \ref{Th-tech-1} we prove the local well-posedness of the equation (\ref{Equation}) in a space $E^{s,\alpha}_T \subset \mathcal{C}([0,T], H^s(\R))$ for $1-\alpha/2 < s \leq 0$, while in Theorem \ref{Th-tech-2} we prove the local-well posedness in a space $F^{s,\alpha}_T \subset \mathcal{C}([0,T], H^s(\R))$ for $0<s$. Only for technical reasons, we shall divide our study in the cases $1-\alpha/2 < s \leq 0$ and $0<s$. Finally, in Theorem \ref{Th-tech-3} we study the regularity of solutions. 

\medskip

{\bf Local well-posedness}.  Solutions of the equation (\ref{Equation}) are constructed as the solutions of the following (equivalent) problem: 
\begin{equation}\label{Equation-Integral}
u(t,x)= K_{\alpha,\beta}(t,\cdot)\ast u_0(x)- \int_{0}^{t}K_{\alpha,\beta}(t-\tau,\cdot)\ast\Big( \gamma_1 \partial_x (u^2)  + \gamma_2 \, \partial^{2}_{x}(u^2) + \gamma_3\, (\partial_x u)^2   \Big)(\tau, x) \, d\tau,   
\end{equation}
where the kernel $K_{\alpha,\beta}(t,x)$ is defined in (\ref{kernel}).  The parameters $\gamma_1,\gamma_2,\gamma_3$ do not play any substantial  role in   local well-posedness theory, so for the sake of simplicity we shall set them as $\gamma_1=\gamma_2=\gamma_3=1$.  On the other hand,  we recall the following well-known estimate on the Beta function, which we shall fully  use to study the nonlinear terms above. For $a>-1$ and $b>-1$ we have
\begin{equation}\label{Beta}
\int_{0}^{t}  (t-\tau)^{a}{\tau}^{b}  d\tau \lesssim t^{a+b+1}.
\end{equation} 

{\underline{\bf Case $1-\alpha/2 < s \leq 0$}}.  Let $7/2<\alpha$ and let $0<T$.   We define the Banach space 
\begin{equation}\label{Espace}
E^{s,\alpha}_{T}= \{ u \in \mathcal{C}([0,T], H^s(\R)): \, \Vert u \Vert_{s,\alpha} <+\infty \}, 
\end{equation}
with the norm 
\begin{equation}\label{Norm}
\Vert u \Vert_{s,\alpha}= \sup_{0\leq t \leq T} \Vert u(t,\cdot)\Vert_{H^s}+ \sup_{0\leq t \leq T} t^{\frac{\vert s \vert}{\alpha}} \Vert u(t,\cdot) \Vert_{L^2}+ \sup_{0\leq t \leq T} t^{\frac{1+\vert s \vert}{\alpha}} \Vert \partial_x u (t,\cdot) \Vert_{L^2}.  
\end{equation}
The second and the third term of this norm will be useful to handle the nonlinear terms in the equation (\ref{Equation-Integral}) (see Proposition \ref{Prop1} below).  Now, for a time $0<T\leq 1$ small enough, we shall construct a solution $u(t,x)$ to this  equation  in the space $E^{s,\alpha}_{T}$. 
\begin{Theoreme}\label{Th-tech-1}  Let $\alpha>\beta >0$ with $\alpha > 7/2$, let $1-\alpha/2 <s \leq 0$, and moreover, define  the quantity 
\begin{equation}\label{eta}
\eta=\frac{s}{\alpha}-\frac{5}{2\alpha}+1>0.
\end{equation}
For any  $u_0 \in H^s(\R)$  there exists a time  
\begin{equation}\label{Control-Tiempo-LWP}
T=T(\Vert u_0\Vert_{H^s}) < \min\left(  \frac{1}{4^{1/\eta} \Vert u_0 \Vert^{1 /  \eta}_{H^s}}, 1  \right),
\end{equation}	and a function $u \in E^{s,\alpha}_{T}$, which is the unique solution  of the equation (\ref{Equation-Integral}). Moreover, the flow-map function $S: H^s(\R) \to E^{s,\alpha}_{T}\subset \mathcal{C}([0,T], H^s(\R))$ defined in (\ref{Map-flow}) is smooth. \end{Theoreme}	
\begin{Remarque}\label{Rmk1} The quantity $\eta$ appears in the estimates to handle the whole nonlinear term $\partial_x(u^2)+\partial^{2}_{x}(u^2)+(\partial_x u)^2$. See the Proposition \ref{Prop-Non-Lin} below. Thus, the constraints $7/2<\alpha $ and $1-\alpha/2 < s$ ensure that $0<\eta$. Indeed, by (\ref{eta}) the inequality $0<\eta$ is equivalent to the inequality $\frac{5}{2}-\alpha<s$, but since $1-\alpha/2 < s$ and $7/2<\alpha$ we can write $\frac{5}{2}-\alpha< 1-\alpha/2<s$.  
\end{Remarque}
\pv  We start by studying the linear term in  the equation (\ref{Equation-Integral}).  
\begin{Proposition}\label{Prop-Lin-1} We have $K_{\alpha,\beta}(t,\cdot)\ast u_0  \in E^{s,\alpha}_{T}$ and $\left\Vert K_{\alpha,\beta}(t,\cdot)\ast u_0 \right\Vert_{s,\alpha} \lesssim \Vert u_0 \Vert_{H^s}$.
\end{Proposition}	
\pv  We shall study separately each term in the norm  $\left\Vert K_{\alpha,\beta}(t,\cdot)\ast u_0 \right\Vert_{s,\alpha}$ defined in (\ref{Norm}).  For the first term, by  Lemma \ref{Lem-Ker-1} (with $s_1=0$) and since  $0<T\leq 1$  we get 
\begin{equation}\label{estim-01}
\sup_{0\leq t \leq T} \left\Vert K_{\alpha,\beta}(t,\cdot)\ast u_0 \right\Vert_{H^s} \lesssim  \Vert u_0 \Vert_{H^s}. 
\end{equation}
We  also  have $K_{\alpha,\beta}(t,\cdot) \ast u_0 \in \mathcal{C}([0,T], H^s(\R))$. Indeed,  on the one hand, for  $t=0$  by  a standard convergence-dominated argument  we get $\ds{\lim_{t \to 0^{+}} \Vert K_{\alpha,\beta}(t,\cdot)\ast u_0 - u_0 \Vert_{H^s}=0}$.  On the other hand, by Lemma \ref{Lem-Ker-2}  (with $s_1=0$) we obtain $K_{\alpha,\beta}(t,\cdot) \ast u_0 \in \mathcal{C}((0,T], H^s(\R))$.

\medskip

To estimate the  second term, first we need  to verify the  following pointwise estimate:
\begin{equation}\label{estim-02}
t^{\frac{\vert s \vert}{\alpha}} \leq \frac{(1+ \left\vert t^{\frac{1}{\alpha}} \xi \right\vert^2)^{\frac{\vert s \vert}{2}}}{(1+\vert \xi \vert^2)^{\frac{\vert s \vert}{2}}}. 
\end{equation} 
Indeed, again by the fact that $0\leq t \leq T \leq 1$,  we just write
\begin{equation*}
t^{\frac{\vert s \vert}{\alpha}} ( 1+\vert \xi \vert^2)^{\frac{\vert s \vert}{2}}=\left( t^{\frac{2}{\alpha}} \right)^{\frac{\vert s \vert}{2}}( 1+\vert \xi \vert^2)^{\frac{\vert s \vert}{2}}= \left(t^{\frac{2}{\alpha}}+ \left\vert t^{\frac{1}{\alpha}} \xi \right\vert^2\right)^{\frac{\vert s \vert}{2}} \leq \left(1+ \left\vert t^{\frac{1}{\alpha}} \xi \right\vert^2\right)^{\frac{\vert s \vert}{2}}. 
\end{equation*}
Once we have the estimate (\ref{estim-02}), we  obtain
\begin{equation*}
\begin{split}
t^{\frac{\vert s \vert}{\alpha}} \left\Vert K_{\alpha,\beta}(t,\cdot) \ast u_0 \right\Vert_{L^2} \leq & \left\Vert  t^{\frac{\vert s \vert}{\alpha}}  \widehat{K_{\alpha,\beta}}(t,\cdot) \, \widehat{u_0} \right\Vert_{L^2} \leq \left\Vert \frac{(1+ \left\vert t^{\frac{1}{\alpha}} \xi \right\vert^2)^{\frac{\vert s \vert}{2}}}{(1+\vert \xi \vert^2)^{\frac{\vert s \vert}{2}}}\,  \widehat{K_{\alpha,\beta}}(t,\cdot) \, \widehat{u_0} \right\Vert_{L^2} \\
\leq & \left\Vert \frac{(1+ \left\vert t^{\frac{1}{\alpha}} \xi \right\vert^2)^{\frac{\vert s \vert}{2}}}{(1+\vert \xi \vert^2)^{\frac{\vert s \vert}{2}} (1+\vert \xi \vert^2)^{s/2}} \, \widehat{K_{\alpha,\beta}}(t,\cdot) \,  (1+\vert \xi \vert^2)^{s/2}\,\widehat{u_0} \right\Vert_{L^2}\\
= &\, \left\Vert  (1+ \left\vert t^{\frac{1}{\alpha}} \xi \right\vert^2)^{\frac{\vert s \vert}{2}} \, \widehat{K_{\alpha,\beta}}(t,\cdot) \,  (1+\vert \xi \vert^2)^{s/2}\,\widehat{u_0} \right\Vert_{L^2}=(A). 
\end{split}
\end{equation*}
Then, by the H\"older inequalities and by Lemma \ref{Key-Lemma} (by setting first $\sigma=0$ and then  $\sigma =\frac{\vert s \vert}{2}$)  for all $0\leq t \leq T \leq 1$  we get
\begin{equation*}
\begin{split}
(A) \lesssim &\,\,  \left\Vert \widehat{K_{\alpha,\beta}}(t,\cdot) \,  (1+\vert \xi \vert^2)^{s/2}\,\widehat{u_0} \right\Vert_{L^2}  + \left\Vert ( \left\vert t^{\frac{1}{\alpha}} \xi \right\vert^2)^{\frac{\vert s \vert}{2}}\, \widehat{K_{\alpha,\beta}}(t,\cdot) \,  (1+\vert \xi \vert^2)^{s/2}\,\widehat{u_0} \right\Vert_{L^2}\\
\lesssim & \,\,\left\Vert \widehat{K_{\alpha,\beta}}(t,\cdot) \right\Vert_{L^\infty} \Vert u_0 \Vert_{H^s}+ \left\Vert ( \left\vert t^{\frac{1}{\alpha}} \xi \right\vert^2)^{\frac{\vert s \vert}{2}}\, \widehat{K_{\alpha,\beta}}(t,\cdot) \right\Vert_{L^\infty} \Vert u_0 \Vert_{H^s} \\
\lesssim & \,\, (e^t + 1) \Vert u_0 \Vert_{H^s}+ (t^{\frac{\vert s \vert}{\alpha}}e^{t}+ t^{\frac{\vert s \vert}{\alpha}} +1 ) \Vert u_0 \Vert_{H^s}\\
\lesssim &  \,\,\Vert u_0 \Vert_{H^s}.
\end{split}
\end{equation*}
We thus obtain 
\begin{equation}\label{estim-03}
\sup_{0\leq t \leq T} t^{\frac{\vert s \vert}{\alpha}} \left\Vert K_{\alpha,\beta}(t,\cdot)\ast u_0 \right\Vert_{L^2} \lesssim \Vert u_0 \Vert_{H^s}. 
\end{equation}

We  study  the third term. First,  we remark   that by  (\ref{estim-02}) we have 
\[ t^{\frac{1+\vert s \vert}{\alpha}} = t^{\frac{1}{\alpha}} \, t^{\frac{\vert s \vert}{\alpha}} \leq t^{\frac{1}{\alpha}}\, \frac{(1+ \left\vert t^{\frac{1}{\alpha}} \xi \right\vert^2)^{\frac{\vert s \vert}{2}}}{(1+\vert \xi \vert^2)^{\frac{\vert s \vert}{2}}}.\]
Then, by following the same computations to prove the estimate (\ref{estim-03}),  we  write
\begin{equation*}
\begin{split}
t^{\frac{1+\vert s \vert}{\alpha}} \left\Vert \partial_x K_{\alpha,\beta}(t,\cdot) \ast u_0 \right\Vert_{L^2} \leq & \,\, \left\Vert t^{\frac{1+\vert s \vert}{\alpha}}  \vert \xi \vert \,  \widehat{K_{\alpha,\beta}}(t,\cdot) \, \widehat{u_0} \right\Vert_{L^2} \\
\leq & \,\,  \left\Vert \left\vert t^{\frac{1}{\alpha}} \xi \right\vert   \frac{(1+ \left\vert t^{\frac{1}{\alpha}} \xi \right\vert^2)^{\frac{\vert s \vert}{2}}}{(1+\vert \xi \vert^2)^{\frac{\vert s \vert}{2}} (1+\vert \xi \vert^2)^{\frac{s}{2}}} \, \widehat{K_{\alpha,\beta}}(t,\cdot) \,  (1+\vert \xi \vert^2)^{\frac{s}{2}}\,\widehat{u_0}  \right\Vert_{L^2}\\
\lesssim & \,\, \left\Vert  \frac{(1+ \left\vert t^{\frac{1}{\alpha}} \xi \right\vert^2)^{\frac{\vert s \vert+1}{2}}}{(1+\vert \xi \vert^2)^{\frac{\vert s \vert}{2}} (1+\vert \xi \vert^2)^{\frac{s}{2}}} \, \widehat{K_{\alpha,\beta}}(t,\cdot) \,  (1+\vert \xi \vert^2)^{\frac{s}{2}}\,\widehat{u_0}  \right\Vert_{L^2} \lesssim \Vert u_0 \Vert_{H^s}. 
\end{split}
\end{equation*}
We thus obtain 
\begin{equation}\label{estim-04}
\sup_{0\leq t \leq T} t^{\frac{1+\vert s \vert}{\alpha}} \left\Vert \partial_x K_{\alpha, \beta}(t,\cdot) \ast u_0 \right\Vert_{L^2} \lesssim \Vert u_0 \Vert_{H^s}. 
\end{equation}
The desired estimate follows from (\ref{estim-01}), (\ref{estim-03}), and (\ref{estim-04}). Proposition \ref{Prop-Lin-1} is proven. \finpv 

\medskip

We study now the nonlinear terms in the equation (\ref{Equation-Integral}). For this, we shall need the following useful technical  estimates. Particularly, we shall observe the use of the second and the third terms in the norm $\Vert \cdot \Vert_{s,\alpha}$ given in (\ref{Norm}). 
\begin{Proposition}\label{Prop1} Let $1-\alpha/2<s\leq 0$, and let $0\leq \sigma \leq 3$. For all $0<t\leq T \leq 1$ the following estimates hold:
\begin{enumerate}
\item $\ds{\left\Vert  \int_{0}^{t} K_{\alpha,\beta}(t-\tau,\cdot)\ast (-\partial^{2}_{x})^{\frac{\sigma}{2}} (u^2)(\tau,\cdot) \, d \tau \right\Vert_{H^s} \lesssim \ t^{ \frac{s}{\alpha}  - \frac{(\sigma+1/2)}{\alpha}+1}\,   \left( \sup_{0\leq t\leq T} t^{\frac{\vert s \vert}{2}} \Vert u(t, \cdot) \Vert_{L^2}\right)^2}$. 
\item $\ds{\left\Vert  \int_{0}^{t} K_{\alpha,\beta}(t-\tau,\cdot)\ast (-\partial^{2}_{x})^{\frac{\sigma}{2}} (u^2)(\tau,\cdot) \, d \tau \right\Vert_{L^2} \lesssim \  t^{ \frac{2s}{\alpha}  - \frac{(\sigma+1/2)}{\alpha}+1}\, \left( \sup_{0\leq t\leq T} t^{\frac{\vert s \vert}{2}} \Vert u(t, \cdot) \Vert_{L^2}\right)^2}$.	
\item  $\ds{\left\Vert  \int_{0}^{t} K_{\alpha,\beta}(t-\tau,\cdot)\ast (-\partial^{2}_{x})^{\frac{\sigma}{2}} ((\partial_x u)^2)(\tau,\cdot) \, d \tau \right\Vert_{H^s} \lesssim \ t^{\frac{s}{\alpha}  - \frac{(\sigma+5/2)}{\alpha} +1}\, \left( \sup_{0\leq t \leq T} t^{\frac{1+\vert s \vert}{\alpha}} \Vert \partial_x u(t, \cdot) \Vert_{L^2}\right)^2}$.
\item  $\ds{\left\Vert  \int_{0}^{t} K_{\alpha,\beta}(t-\tau,\cdot)\ast (-\partial^{2}_{x})^{\frac{\sigma}{2}} ((\partial_x u)^2)(\tau,\cdot) \, d \tau \right\Vert_{L^2} \lesssim \ t^{\frac{2s}{\alpha}  - \frac{(\sigma+5/2)}{\alpha} +1}\, \left( \sup_{0\leq t \leq T} t^{\frac{1+\vert s \vert}{\alpha}} \Vert \partial_x u(t, \cdot) \Vert_{L^2}\right)^2}$.
\end{enumerate}		
\end{Proposition}	
\pv Let us prove the first point. Since $s \leq 0$   we can write 
\begin{equation}\label{estim-05}
\begin{split}
&\left\Vert \int_{0}^{t} K_{\alpha,\beta}(t-\tau,\cdot)\ast (-\partial^{2}_{x})^{\frac{\sigma}{2}} (u^2)(\tau,\cdot)d \tau \right\Vert_{H^s}  \leq \int_{0}^{t} \left\Vert (1+\vert \xi \vert^2)^{s/2} \vert \xi \vert^{\sigma} \, \widehat{K_{\alpha,\beta}}(t-\tau,\cdot) (\widehat{u} \ast \widehat{u})(\tau,\cdot) \right\Vert_{L^2} d \tau \\
\leq & \int_{0}^{t} \left\Vert \vert \xi \vert^{s+\sigma}  \widehat{K_{\alpha,\beta}}(t-\tau,\cdot) (\widehat{u} \ast \widehat{u})(\tau,\cdot) \right\Vert_{L^2} d \tau \leq  \int_{0}^{t} \left\Vert \vert \xi \vert^{s+\sigma} \widehat{K_{\alpha,\beta}}(t-\tau,\cdot) \right\Vert_{L^2}\, \left\Vert  (\widehat{u} \ast \widehat{u})(\tau,\cdot) \right\Vert_{L^\infty} d \tau.
\end{split}
\end{equation} 
The first expression on the right-hand side can be estimated as follows:
\begin{equation}\label{estim-06}
\left\Vert \vert \xi \vert^{s+\sigma} \widehat{K_{\alpha,\beta}}(t-\tau,\cdot) \right\Vert_{L^2} \lesssim (t-\tau)^{-\frac{s+\sigma}{\alpha}-\frac{1}{2\alpha}}. 
\end{equation}
Indeed, we split  $\ds{\left\vert  \widehat{K_{\alpha,\beta}}(t-\tau, \xi) \right\vert= e^{(\vert \xi \vert^\alpha - \vert \xi \vert^\beta) (t-\tau)}=e^{-\frac{\vert \xi \vert^\alpha}{2} (t-\tau)}\, e^{-\left( \frac{\vert \xi \vert^\alpha}{2} - \vert \xi \vert^\beta\right) (t-\tau)}}$; and  for $\kappa= (t-\tau)^{\frac{1}{\alpha}} \xi$ we  write
\begin{equation*}
\begin{split}
\left\Vert \vert \xi \vert^{s+\sigma} \widehat{K_{\alpha,\beta}}(t-\tau,\cdot) \right\Vert_{L^2} \leq & \,\, \left\Vert \vert \xi \vert^{s+\sigma}  e^{-\frac{\vert \xi \vert^\alpha}{2}(t-\tau)} \right\Vert_{L^2}\, \left\Vert e^{-\left( \frac{\vert \xi \vert^\alpha}{2} - \vert \xi \vert^\beta\right)(t-\tau)} \right\Vert_{L^\infty} \\
\leq  & \,\, (t-\tau)^{-\frac{s+\sigma}{\alpha}}  \left\Vert \left\vert (t-\tau)^{\frac{1}{\alpha}} \xi \right\vert^{s+\sigma} e^{-\frac{\left\vert (t-\tau)^{\frac{1}{\alpha}} \xi \right\vert^\alpha}{2}}  \right\Vert_{L^2}\,   \left\Vert e^{-\left( \frac{\vert \xi \vert^\alpha}{2} - \vert \xi \vert^\beta\right)(t-\tau)} \right\Vert_{L^\infty} \\
\leq  & \,\, (t-\tau)^{-\frac{s+\sigma}{\alpha}-\frac{1}{2\alpha}} \,  \left\Vert \left\vert \kappa \right\vert^{s+\sigma} e^{-\frac{\left\vert \kappa \right\vert^\alpha}{2}}  \right\Vert_{L^2}\,   \left\Vert e^{-\left( \frac{\vert \xi \vert^\alpha}{2} - \vert \xi \vert^\beta\right)(t-\tau)} \right\Vert_{L^\infty} \\
\lesssim &\,\,  (t-\tau)^{-\frac{s+\sigma}{\alpha}-\frac{1}{2\alpha}}. 
\end{split}
\end{equation*}
To estimate the second expression on the right-hand side, by the Young inequalities (with $1+1/ \infty= 1/2+1/2$), by the Plancherel's identity and by the second term in the  norm $\Vert \cdot \Vert_{s,\alpha}$,  we have:
\begin{equation}\label{estim-07}
\left\Vert ( \widehat{u} \ast \widehat{u})(\tau,\cdot) \right\Vert_{L^\infty} \lesssim \tau^{-\frac{2 \vert s \vert}{\alpha}}\, \left( \sup_{0\leq \tau \leq T} \tau^{\frac{2\vert s \vert}{\alpha}} \Vert u(\tau, \cdot) \Vert^{2}_{L^2} \right).
\end{equation} 
With the estimates (\ref{estim-06}) and (\ref{estim-07}) at hand, we get back to (\ref{estim-05}) and  write 
\begin{equation}\label{Estim-01}
\begin{split}
&\, \left\Vert \int_{0}^{t} K_{\alpha,\beta}(t-\tau,\cdot)\ast (-\partial^{2}_{x})^{\frac{\sigma}{2}} (u^2)(\tau,\cdot)d \tau \right\Vert_{H^s} \\
\lesssim &\,   \left(  \int_{0}^{t}  (t-\tau)^{-\frac{s+\sigma}{\alpha}-\frac{1}{2\alpha}}\, \tau^{-\frac{2 \vert s \vert }{\alpha}}\, d \tau\right)\,  \left( \sup_{0\leq t\leq T} t^{\frac{2\vert s \vert}{\alpha}} \Vert u(t, \cdot) \Vert^{2}_{L^2}\right).
\end{split}
\end{equation} 
To estimate the integral above, in the formula  (\ref{Beta}) we  set $a= -\frac{s+\sigma}{\alpha}-\frac{1}{2\alpha}$ and $b=-\frac{2\vert s \vert}{\alpha}$.
\begin{Remarque} Since $0\leq \sigma \leq 3$  and $s \leq 0$ we have $a= -\frac{s+\sigma}{\alpha}-\frac{1}{2\alpha}>-1$ as long as $\alpha > 7/2$. Moreover, since  $s>1-\alpha/2 > - \alpha/2$ we have $-\frac{2\vert s \vert}{\alpha}=\frac{2s}{\alpha}>-1$.
\end{Remarque}
Thus, a direct application of   (\ref{Beta}) yields  $\ds{\int_{0}^{t}  (t-\tau)^{-\frac{s+\sigma}{\alpha}-\frac{1}{2\alpha}}\, \tau^{\frac{2 s }{\alpha}}\, d \tau \lesssim t^{s/\alpha -(\sigma+1/2)/\alpha +1}}$. With this estimate  the first point of this proposition follows.  

\medskip

The other points of this proposition follow similar estimates. Indeed, for the second point, we just remark that  in  the estimate (\ref{estim-06}) we  set now  $s=0$, which yields the integral $\ds{\int_{0}^{t}  (t-\tau)^{-\frac{\sigma}{\alpha}-\frac{1}{2\alpha}}\, \tau^{\frac{2 s }{\alpha}}\, d \tau \lesssim t^{2s/\alpha -(\sigma+1/2)/\alpha +1}}$.  For the third point,  we follow the same computations in the estimates (\ref{estim-05}) and  (\ref{estim-06}) to write 
\begin{equation}
\begin{split}
&\left\Vert  \int_{0}^{t} K_{\alpha,\beta}(t-\tau,\cdot)\ast (-\partial^{2}_{x})^{\frac{\sigma}{2}} ((\partial_x u)^2)(\tau,\cdot) \, d \tau \right\Vert_{H^s} \lesssim \,\, \int_{0}^{t}  (t-\tau)^{-\frac{s+\sigma}{\alpha}-\frac{1}{2\alpha}} \, \left\Vert (\widehat{\partial_x u} \ast \widehat{\partial_x u}) (\tau,\cdot) \right\Vert_{L^\infty} d \tau \\
\lesssim &\,\,   \left( \int_{0}^{t} (t-\tau)^{-\frac{s+\sigma}{\alpha}-\frac{1}{2\alpha}} \, \tau^{-\frac{2(1+\vert s \vert)}{\alpha}} d \tau \right)\, \left( \sup_{0\leq t \leq T} t^{\frac{2(1+\vert s \vert)}{\alpha}} \Vert \partial_x u(t, \cdot) \Vert^{2}_{L^2}\right).
\end{split}
\end{equation}
In the formula  (\ref{Beta}) we set now $b=-\frac{2(1+\vert s \vert)}{\alpha}$.
\begin{Remarque}
We have $-\frac{2(1+\vert s \vert)}{\alpha}>-1$ as long as $1-\alpha / 2 <s\leq 0$.
\end{Remarque}
Then we  obtain $\ds{\int_{0}^{t} \tau^{-\frac{s+\sigma}{\alpha}-\frac{1}{2\alpha}} \, (t-\tau)^{-\frac{2(1+\vert s \vert)}{\alpha}} d \tau  \lesssim t^{s/\alpha - (\sigma+5/2)/\alpha +1}}$. Finally, for the fourth point, we also  follow the same computations above, where we have the integral $\ds{\int_{0}^{t} (t-\tau)^{-\frac{\sigma}{\alpha}-\frac{1}{2\alpha}} \, \tau^{-\frac{2(1+\vert s \vert)}{\alpha}} d \tau  \lesssim t^{2s/\alpha - (\sigma+5/2)/\alpha +1}}$. Proposition \ref{Prop1} is proven. \finpv 

\medskip

With these estimates at our disposal, we directly obtain the following proposition. 

\begin{Proposition}\label{Prop-Non-Lin} Let $\eta>0$ be the quantity defined in (\ref{eta}). The following estimates hold:
\begin{enumerate}
\item  	$\ds{\left\Vert \int_{0}^{t} K_{\alpha,\beta}(\tau, \cdot)\ast \partial_x (u^2)(t-\tau,\cdot)\, d \tau \right\Vert_{s,\alpha} \lesssim T^{\eta}\, \Vert u \Vert^{2}_{s,\alpha}}$.
\item  $\ds{\left\Vert \int_{0}^{t} K_{\alpha,\beta}(\tau,\cdot ) \ast\partial^{2}_{x} (u^2)(t-\tau,\cdot)\, d \tau \right\Vert_{s,\alpha} \lesssim T^{\,\eta}\, \Vert u \Vert^{2}_{s,\alpha}}$.
\item  $\ds{\left\Vert \int_{0}^{t} K_{\alpha,\beta}(\tau, \cdot) \ast (\partial_x u)^2(t-\tau,\cdot)\, d \tau \right\Vert_{s,\alpha} \lesssim T^{\eta}\, \Vert u \Vert^{2}_{s,\alpha}}$.
\end{enumerate}		
\end{Proposition} 
\pv  The first estimate follows from the first point  and the second  point  of Proposition \ref{Prop1} with $\sigma=1$ and $\sigma =2$. The second estimate also follows from these same points whit $\sigma =2$ and $\sigma=3$. Finally, the third estimate follows from the third point  and the fourth point of Proposition \ref{Prop1} with $\sigma=0$ and $\sigma=1$.  \finpv

\medskip 

Consequently, the existence and uniqueness of a local in-time solution  $u\in E^{s,\alpha}_{T} \subset  \mathcal{C}([0,T], H^s(\R))$ follow from standard arguments, provided that the condition (\ref{Control-Tiempo-LWP}) holds. Moreover, the smoothness of the flow-map function $S: H^s(\R) \to E^{s,\alpha}_{T}$ also follows from well-known  arguments, see for instance \cite{Kenig-Ponce-Vega}. Theorem \ref{Th-tech-1} is proven. \finpv 

\medskip

{ \bf \underline{Case $s>0$}.}  The key idea to prove the local-well posedness, in this case, is to use the estimates performed above.  We thus start by proving the following useful lemma, which is a product law-type in the Sobolev spaces. Let us mention  that for a parameter $z\in \R$ we shall denote  the Bessel potential $(I_d - \partial^{2}_{x})^{\frac{z}{2}}= \mathcal{J}^{z}$, which is defined in the Fourier level by the symbol $(1+|\xi|^2)^{\frac{z}{2}}$. 
\begin{Lemme}\label{Lema-tech1} Let $s_1 \leq 0 < s$.  The following estimate holds: $\Vert f\,g \Vert_{H^s} \lesssim \left\Vert  (\mathcal{J}^{s-s_1} f) \,  g \right\Vert_{H^{s_1}} + \left\Vert f \, (\mathcal{J}^{s-s_1} g) \right\Vert_{H^{s_1}} $. 
\end{Lemme}	
\pv The proof follows from the  pointwise estimate:
\[  (1+\vert \xi \vert^2)^{\frac{s}{2}} \vert (\widehat{f} \ast \widehat{g}) (\xi) \vert \lesssim (1+\vert \xi \vert^2)^{\frac{s_1}{2}} \ \left( \left( (1+\vert \xi \vert^2)^{\frac{s-s_1}{2}} \vert  \widehat{f} \vert \right) \ast  \vert \widehat{g} \vert \right) (\xi)   +  \left( \vert \widehat{f} \vert \ast  \left( (1+\vert \xi \vert^2)^{\frac{s-s_1}{2}}  \vert \widehat{g} \vert\right)\right) (\xi). \qquad \blacksquare  \]
With this lemma at our disposal,  we are able  to estimate the product $f\, g $ in the norm of the space  $H^{s}(\R)$ (with $s>0$) in terms of the products $(\mathcal{J}^{s-s_1} f) g$ and $f  ( \mathcal{J}^{s-s_1} g)$ in the norm of the space $H^{s_1}(\R)$ (with $s_1 \leq 0$). Consequently, we can use the estimates above as follows:    for $\alpha > \frac{7}{2}$ we set  $1-\frac{\alpha}{2}< s_1  \leq 0$. Then, for $s>0$ and for a time $0\leq T\leq 1$ small enough, we define the Banach space
\[ F^{s,\alpha}_{T}=\{  u \in \mathcal{C}([0,T], H^s(\R)) : \Vert u \Vert_{s,\alpha,s_1}<+\infty \}, \]
with the norm
\begin{equation}\label{Norm-F}
\begin{split}
\Vert u \Vert_{s,\alpha,s_1}= & \,  \sup_{0\leq t \leq T} \Vert u(t,\cdot)\Vert_{H^s}+
\sup_{0\leq t \leq T} t^{\frac{\vert s_1 \vert}{\alpha}} \Vert u(t,\cdot) \Vert_{L^2} + \sup_{0\leq t \leq T} t^{\frac{1+\vert s_1 \vert}{\alpha}} \Vert \partial_x u(t,\cdot) \Vert_{L^2} \\
&+  \sup_{0\leq t \leq T} t^{\frac{\vert s_1 \vert}{\alpha}} \Vert  \mathcal{J}^{s-s_1}\, u(t,\cdot) \Vert_{L^2}+ \sup_{0\leq t \leq T} t^{\frac{1+\vert s_1 \vert}{\alpha}} \Vert \partial_x  \, \mathcal{J}^{s-s_1}\, u(t,\cdot)   \Vert_{L^2}. 
\end{split} 
\end{equation}
Let us briefly explain this norm. The second and the third terms are the same as  used in the  norm $\Vert \cdot \Vert_{\sigma,\alpha}$ defined in (\ref{Norm}). Moreover, as we consider here $s>0$,   the fourth and the fifth term  will allow us to easily estimate the nonlinear terms of the equation (\ref{Equation-Integral}) in the  norm of the space $H^s(\R)$. Thus, we can state our second technical theorem.
\begin{Theoreme}\label{Th-tech-2}  Let $\alpha>\beta >0$ with $\alpha > 7/2$, let $s>0$,  $1-\frac{\alpha}{2}<s_1\leq 0$, and  let $\eta>0$  be the quantity given in (\ref{eta}).

\medskip

For any  $u_0 \in H^s(\R)$, there exists a time  $T=T(\Vert u_0\Vert_{H^s})$ given in (\ref{Control-Tiempo-LWP}) and  there exists $u \in F^{s,\alpha}_{T}$ a unique solution  to the equation (\ref{Equation-Integral}). Moreover, the flow-map function $S: H^s(\R) \to F^{s,\alpha,s_1}_{T} \subset \mathcal{C}([0,T], H^s(\R))$ defined in (\ref{Map-flow}) is smooth.  
\end{Theoreme}	 
\pv As mentioned, the proof  uses the estimates already proven in the previous case when $1-\alpha/2 < s \leq 0$; and it follows very similar ideas. So, it is enough to give a brief proof. The linear term in the equation (\ref{Equation-Integral})) is easy to estimate and for $u_0 \in H^s(\R)$ we have  
\begin{equation}\label{Lin-Norm-F}
 \Vert K_{\alpha,\beta}(t,\cdot) \ast u_0  \Vert_{s,\alpha,s_1} \lesssim \Vert u_0 \Vert_{H^s}.   
\end{equation}

\medskip

We study now the nonlinear terms.  For the first term in the norm $\Vert \cdot \Vert_{s,\alpha,\sigma}$ (given in (\ref{Norm-F})),  by Lemma \ref{Lema-tech1}, by the first point of Proposition \ref{Prop1} (with $\sigma=1$), and moreover, by recalling that $\eta= \frac{s}{\alpha}-\frac{5}{2\alpha}+1 <  \frac{s}{\alpha}-\frac{3}{2\alpha}+1$, for $0<t\leq T \leq 1$  we write 
\begin{equation*}
\begin{split}
& \left\Vert \int_{0}^{t} K_{\alpha,\beta}(t-\tau, \cdot)\ast \partial_{x}(u^2)(\tau,\cdot) d\tau  \right\Vert_{H^s} \\
\lesssim  &  \left\Vert \int_{0}^{t} K_{\alpha,\beta}(t-\tau, \cdot)\ast \partial_{x}( (\mathcal{J}^{s-s_1} u )   u)(\tau,\cdot) d\tau  \right\Vert_{H^{s_1}}  + \left\Vert \int_{0}^{t} K_{\alpha,\beta}(t-\tau, \cdot)\ast \partial_{x}( u \, (\mathcal{J}^{s-s_1} u))(\tau,\cdot) d\tau  \right\Vert_{H^{s_1}}\\
\lesssim & \,\,   t^{\frac{s}{\alpha}-\frac{3}{2\alpha}+1}     \left( \sup_{0\leq t \leq T} t^{\frac{\vert \sigma \vert}{\alpha}} \, \Vert \mathcal{J}^{s-\sigma} u (t,\cdot) \Vert_{L^2}\right) \left( \sup_{0\leq t \leq T} t^{\frac{\vert s_1 \vert}{\alpha}} \, \Vert u (t,\cdot) \Vert_{L^2}\right) \lesssim  \,\,  t^{\frac{s}{\alpha}-\frac{5}{2\alpha}+1}  \,  \Vert u \Vert^{2}_{s,\alpha,s_1}\leq T^{\eta}\, \Vert u \Vert^{2}_{s,\alpha,s_1}. 
\end{split} 
\end{equation*} 
 The other  terms  $\ds{\left\Vert \int_{0}^{t} K_{\alpha,\beta}(t-\tau, \cdot)\ast \partial^{2}_{x}(u^2)(\tau,\cdot) d\tau  \right\Vert_{H^s}}$  and $\ds{\left\Vert \int_{0}^{t} K_{\alpha,\beta}(t-\tau, \cdot)\ast (\partial_x u)^2(\tau,\cdot) d\tau  \right\Vert_{H^s}}$  are  treated similarly, where we use again   Lemma \ref{Lema-tech1} as well as  Proposition \ref{Prop1}.   Moreover, remark that  the second to the fifth expressions  in the norm $\Vert \cdot \Vert_{s,\alpha,\sigma}$ (see (\ref{Norm-F})) were already estimated in Proposition \ref{Prop1}.  Thus,  the following estimate holds: 
\begin{equation}\label{NonLin-Norm-F}
\left\Vert \int_{0}^{t}K_{\alpha,\beta}(\tau,\cdot)\ast\Big(  \partial_x (u^2)  +  \partial^{2}_{x}(u^2) + (\partial_x u)^2   \Big)(t-\tau, \cdot) \, d\tau \right\Vert_{s,\alpha,s_1} \lesssim T^\eta \,  \Vert u \Vert^{2}_{s,\alpha,s_1}.    
\end{equation}
 
Consequently,  Theorem \ref{Th-tech-2} follows from arguments already studied in the previous case when $1-\frac{\alpha}{2}<s\leq 0$.   \finpv 

\medskip

{\bf Regularity of solutions}. In  our last technical theorem,  we study the regularity (in the spatial variable) of solutions constructed above. We recall the standard  notation $\ds{H^{\infty}(\R)=\bigcap_{r\geq s}H^r(\R)}$.
\begin{Theoreme}\label{Th-tech-3} Let $\alpha>\beta>0$, with $\alpha >7/2$.  Let $u \in E^{s,\alpha}_{T}$ (when $1-\alpha/ 2 <s \leq 0$) or let $u \in F^{s,\alpha,\sigma}_{T}$ (when $0<s$) be the solution of the integral equation (\ref{Equation-Integral}) given by Theorems \ref{Th-tech-1} and \ref{Th-tech-2} respectively. Then we have $u \in \mathcal{C}((0,T], H^{\infty}(\R))$.  
\end{Theoreme}	
\pv  We shall prove that each term on the right-hand side of the equation (\ref{Equation-Integral}) belongs to the space $\mathcal{C}((0,T], H^{\infty}(\R))$.  For the linear term $K_{\alpha,\beta}(t,\cdot)\ast u_0$ (with $u_0 \in H^s(\R)$)  by Lemmas \ref{Lem-Ker-1} and \ref{Lem-Ker-2} we directly have  $K_{\alpha,\beta}(t,\cdot)\ast u_0 \in \mathcal{C}((0,T], H^{\infty}(\R))$.

\medskip 

We study now the nonlinear term in  (\ref{Equation-Integral}), where (for the sake of clearness) we shall consider the cases $1-\alpha/2 s \leq 0$ and $0<s$ separately. 

\medskip

\underline{{\bf Case $1-\alpha/2 < s \leq 0$}}. For the sake of simplicity, we shall write 
\begin{equation}\label{def-B}
B(u,u)=  \partial_x(u^2) + \partial^{2}_{x}(u^2) +(\partial_x u)^2,
\end{equation}
and  for all $0<t\leq T$ fixed,  we will prove that there exists $0<\delta$ small enough such that  we have
\begin{equation}\label{estim-reg} 
\begin{split}
\left\Vert  \int_{0}^{t} K_{\alpha,\beta}(t-\tau, \cdot) \ast  B(u,u)(\tau,\cdot) d \tau  \right\Vert_{H^{s+\delta}} \lesssim  t^{\frac{2s}{\alpha} - \frac{s+\delta+5/2}{\alpha}+1}\,  \Vert u \Vert^{2}_{s,\alpha}\, ,  \quad \text{with}\quad \frac{2s}{\alpha} - \frac{s+\delta+5/2}{\alpha}+1>0.
\end{split}
\end{equation}
We consider  here the following subcases: first  when $s+\delta \leq 0$ and thereafter when $0<s+\delta$. 

\medskip

In the case  $s+\delta \leq 0$, by the first point of Proposition \ref{Prop1} (with $\sigma=1$ and $\sigma=2$ respectively), and moreover, by the third point of Proposition \ref{Prop1}  (with $\sigma=0$) we directly have
\[ \left\Vert  \int_{0}^{t} K_{\alpha,\beta}(t-\tau,\cdot) \ast B(u,u) (\tau,\cdot) d \tau \right\Vert_{H^{s+\delta}} \lesssim  \,\, \left( t^{\frac{s+\delta}{\alpha}-\frac{3}{2 \alpha} +1} + t^{\frac{s+\delta}{\alpha}-\frac{5}{2 \alpha} +1}\right) \Vert u \Vert^{2}_{s,\alpha}.  \]
Moreover, since  $1-\alpha/2 < s < s+\delta \leq 0$; and as we have $0<t\leq T \leq 1$, the term on the right-hand side is estimated from above by $t^{\frac{s+\delta}{\alpha}-\frac{5}{2 \alpha} +1}$. We thus get: 
\begin{equation}\label{case-1}
\left\Vert  \int_{0}^{t} K_{\alpha,\beta}(t-\tau,\cdot) \ast B(u,u) (\tau,\cdot) d \tau \right\Vert_{H^{s+\delta}} \lesssim  \,\,  t^{\frac{s+\delta}{\alpha}-\frac{5}{2 \alpha} +1}\, \Vert u \Vert^{2}_{s,\alpha}, \qquad \text{with}\quad  \frac{s+\delta}{\alpha}-\frac{5}{2 \alpha} +1>0. 
\end{equation}
We consider now the case when $0<s+\delta$. Here  we write:
\begin{equation}\label{estim-reg-0}
\begin{split}
\left\Vert  \int_{0}^{t} K_{\alpha,\beta}(t-\tau, \cdot) \ast  B(u,u)(\tau,\cdot) d \tau  \right\Vert_{H^{s+\delta}}
=&\,\,  \left\Vert  \int_{0}^{t} K_{\alpha,\beta}(t-\tau, \cdot) \ast  B(u,u)(\tau,\cdot) d \tau  \right\Vert_{L^{2}} \\
&\,\, + \left\Vert   \int_{0}^{t} K_{\alpha,\beta}(t-\tau, \cdot) \ast (-\partial^{2}_{x})^{\frac{s+\delta}{2}}  B(u,u)(\tau,\cdot) d \tau \right\Vert_{L^{2}},
\end{split}
\end{equation}
where we must estimate each term on the right-hand side. For the first term, by the second point of Proposition \ref{Prop1} (with $\sigma=1$ and $\sigma=2$) and by the fourth term of Proposition \ref{Prop1} (with $\sigma=0$) we have 
\begin{equation}\label{estim-reg-1}
\left\Vert  \int_{0}^{t} K_{\alpha,\beta}(t-\tau, \cdot) \ast  B(u,u)(\tau,\cdot) d \tau  \right\Vert_{L^{2}}   \lesssim  \,\,  t^{\frac{2s}{\alpha}-\frac{5}{2 \alpha} +1}\, \Vert u \Vert^{2}_{s,\alpha}\, , \qquad \text{with}\quad  \frac{s}{\alpha}-\frac{5}{2 \alpha} +1>0.
\end{equation}
For the second term, we use again the second point of Proposition \ref{Prop1}  (with $\sigma=s+\delta+1$ and  $\sigma=s+\delta+2$) and  we use again the fourth point of Proposition \ref{Prop1} (with $\sigma=\frac{s+\delta}{2}$).  Moreover, we set $0<\delta < s-\frac{5}{2}+\alpha$ (since  $\alpha >7/2$ and $1-\alpha/2<s$ we have  $0<s-\frac{5}{2}+\alpha$) to obtain that $\frac{2s}{\alpha} - \frac{s+\delta+5/2}{\alpha}+1>0$. Then we have
\begin{equation}\label{estim-reg-2}
\left\Vert   \int_{0}^{t} K_{\alpha,\beta}(t-\tau, \cdot) \ast (-\partial^{2}_{x})^{\frac{s+\delta}{2}}  B(u,u)(\tau,\cdot) d \tau \right\Vert_{L^{2}} \lesssim  t^{\frac{2s}{\alpha} - \frac{s+\delta+5/2}{\alpha}+1}\, \Vert u \Vert^{2}_{s,\alpha}.
\end{equation}
Once we have the estimates (\ref{estim-reg-1}) and (\ref{estim-reg-2}), we remark that $\frac{2s}{\alpha}-\frac{5}{2 \alpha} +1> \frac{2s}{\alpha} - \frac{s+\delta+5/2}{\alpha}+1$, hence, since  $0<t\leq T \leq 1$  we get $ t^{\frac{2s}{\alpha}-\frac{5}{2 \alpha} +1} \leq t^{\frac{2s}{\alpha} - \frac{s+\delta+5/2}{\alpha}+1}$.  We  get back  to  (\ref{estim-reg-0}) and we obtain
\begin{equation}\label{case-2}
\left\Vert  \int_{0}^{t} K_{\alpha,\beta}(t-\tau, \cdot) \ast  B(u,u)(\tau,\cdot) d \tau  \right\Vert_{H^{s+\delta}} \lesssim  t^{\frac{2s}{\alpha} - \frac{s+\delta+5/2}{\alpha}+1}\, \Vert u \Vert^{2}_{s,\alpha}, \quad \text{with}\quad \frac{2s}{\alpha} - \frac{s+\delta+5/2}{\alpha}+1>0.
\end{equation}

Thus,  for both cases when $s+\delta \leq 0$ and $0<s+\delta$, the wished  estimate (\ref{estim-reg}) follows from  (\ref{case-1}) and (\ref{case-2}) respectively. Moreover, remark that we have $t^{\frac{s+\delta}{\alpha}-\frac{5}{2 \alpha} +1} \leq t^{\frac{2s}{\alpha} - \frac{s+\delta+5/2}{\alpha}+1}$. 

\medskip

We study now the continuity in the time variable. Let  $\varepsilon>0$ and let $\varepsilon<t_1,t_2\leq T\leq 1$, without loss of generality we assume that $t_1<t_2$. Then  we write 
\begin{equation}
\begin{split}
&\left\Vert  \int_{0}^{t_2} K_{\alpha,\beta}(t_2-\tau, \cdot)\ast B(u,u)(\tau, \cdot) d \tau -  \int_{0}^{t_1} K_{\alpha,\beta}(t_1-\tau, \cdot)\ast B(u,u)(\tau, \cdot) d \tau \right\Vert_{H^{s+\delta}}  \\
\leq  & \left\Vert  \int_{t_1}^{t_2} K_{\alpha,\beta}(t_2-\tau, \cdot)\ast B(u,u)(\tau, \cdot) d \tau \right\Vert_{H^{s+\delta}} \\
&+ \left\Vert  \int_{0}^{t_1} ( K_{\alpha,\beta}(t_2-\tau, \cdot)-K_{\alpha,\beta}(t_1-\tau, \cdot))\ast B(u,u)(\tau, \cdot) d \tau \right\Vert_{H^{s+\delta}}.  
\end{split} 
\end{equation}
For the first term above, by following  the same estimates performed in (\ref{estim-reg}) we have 
\begin{equation}\label{estim-cont-1}
\left\Vert  \int_{t_1}^{t_2} K_{\alpha,\beta}(t_2-\tau, \cdot)\ast B(u,u)(\tau, \cdot) d \tau \right\Vert_{H^{s+\delta}} \lesssim (t_2-t_1)^{\frac{2s}{\alpha} -\frac{s+\delta +5/2}{\alpha}+1}\, \Vert u \Vert^{2}_{s,\alpha}.
\end{equation} 
We estimate  the second term above. Getting back to the expression (\ref{kernel}), we have 
\begin{equation*}
\begin{split}
&\,\left\Vert \int_{0}^{t_1} ( K_{\alpha,\beta}(t_2-\tau, \cdot)-K_{\alpha,\beta}(t_1-\tau, \cdot))\ast B(u,u)(\tau, \cdot) d \tau  \right\Vert_{H^{s+\delta}} \\
\leq &\,\,  \int_{0}^{t_1} \left\Vert  (1+\vert \xi \vert^2)^{\frac{s+\delta}{2}}\, \left\vert  e^{-f(\xi)(t_2-\tau)} - e^{-f(\xi)(t_1-\tau)} \right\vert \widehat{B(u,u)}(\tau,\cdot)\right\Vert_{L^2}\, d \tau,
\end{split}
\end{equation*}
where we study the expression $ \left\vert e^{-f(\xi)(t_2-\tau)} - e^{-f(\xi)(t_1-\tau)} \right\vert$. First, we remark that the function $\vert f(\xi) \vert$ (given in (\ref{kernel})) is of polynomial growth and, for $i=1,2$ we have $\vert e^{-f(\xi) (t_i -\tau)} \vert = e^{-(\vert \xi \vert^\alpha - \vert \xi \vert^\beta)(t_i -\tau)} \lesssim 1$.  Then,  by the mean value theorem in the temporal variable, there exists $t_0 \in (0, t_2-t_1)$ such that 
\begin{equation*}
\begin{split}
&\left\vert   e^{-f(\xi)(t_2-\tau)} - e^{-f(\xi)(t_1-\tau)} \right\vert  = \,\, \vert  e^{-f(\xi)(t_1-\tau)}\vert \left\vert  \left( e^{-f(\xi)(t_2-t_1)} - 1 \right) \right\vert \\
\lesssim &\,\,\vert  e^{-f(\xi)(t_1-\tau)}\vert\,  \vert f(\xi) \vert e^{-(\vert \xi \vert^\alpha - \vert \xi \vert^\beta)t_0 } (t_2-t_1)  \lesssim    \vert f(\xi) \vert  \, e^{-(\vert \xi \vert^\alpha - \vert \xi \vert^\beta) t_0} (t_2-t_1)\\
\lesssim  & \,\, \left( \vert f(\xi) \vert  \, e^{-(\vert \xi \vert^\alpha - \vert \xi \vert^\beta) \frac{t_0}{2}} \right) e^{-(\vert \xi \vert^\alpha - \vert \xi \vert^\beta) \frac{t_0}{2}}  (t_2-t_1) \lesssim e^{-(\vert \xi \vert^\alpha - \vert \xi \vert^\beta) \frac{t_0}{2}}  (t_2-t_1). 
\end{split} 
\end{equation*}
We get back to the previous estimate   and we use the definition of $B(u,u)$ given in (\ref{def-B})  to write
\begin{equation*}
\begin{split}
&\int_{0}^{t_1} \left\Vert  (1+\vert \xi \vert^2)^{\frac{s+\delta}{2}}\, \left\vert  e^{-f(\xi)(t_2-\tau)} - e^{-f(\xi)(t_1-\tau)} \right\vert \widehat{B(u,u)}(\tau,\cdot)\right\Vert_{L^2}\, d \tau\\
\lesssim & \,\,  (t_2-t_1) \int_{0}^{t_1} \left\Vert  (1+\vert \xi \vert^2)^{\frac{s+\delta}{2}}\,  e^{-(\vert \xi \vert^\alpha - \vert \xi \vert^\beta) \frac{t_0}{2}}  (\vert \xi \vert + \vert \xi \vert^2) (\widehat{u}\ast \widehat{u})(\tau,\cdot)  \right\Vert_{L^2}\, d \tau \\
&\,\, + (t_2-t_1) \int_{0}^{t_1} \left\Vert  (1+\vert \xi \vert^2)^{\frac{s+\delta}{2}}\,  e^{-(\vert \xi \vert^\alpha - \vert \xi \vert^\beta) \frac{t_0}{2}}   (\widehat{\partial_x u}\ast \widehat{\partial_x u})(\tau,\cdot)  \right\Vert_{L^2}\, d \tau  \\
\lesssim &\,\, (t_2-t_1)\Vert u \Vert^{2}_{s,\alpha}. 
\end{split}
\end{equation*}
We thus have 
\begin{equation}\label{estim-cont-2}
\left\Vert \int_{0}^{t_1} ( K_{\alpha,\beta}(t_2-\tau, \cdot)-K_{\alpha,\beta}(t_1-\tau, \cdot))\ast B(u,u)(\tau, \cdot) d \tau  \right\Vert_{H^{s+\delta}}  \lesssim \,\, (t_2-t_1)\Vert u \Vert^{2}_{s,\alpha}. 
\end{equation}
Finally, by the estimates (\ref{estim-cont-1}) and (\ref{estim-cont-2}) we obtain $\ds{\int_{0}^{t} K_{\alpha,\beta}(t-\tau, \cdot) \ast  B(u,u)(\tau,\cdot) d \tau  \in \mathcal{C}((0,T], H^{s+\delta}(\R))}$, and consequently we have $u \in \mathcal{C}((0,T], H^{s+\delta}(\R))$ for $0<\delta < s-\frac{5}{2}+\alpha$. By  bootstrapping this procedure  (to obtain
a gain of regularity for the nonlinear term) we conclude that $u \in \mathcal{C}((0,T], H^{\infty}(\R))$.

\medskip 

\underline{{\bf The case $0<s$}}. In this case, we have very similar estimates to the previous ones: we  essentially follow  the ideas of the proof of Theorem \ref{Th-tech-2}, we use  the Lemma \ref{Lema-tech1} and the norm $\Vert u \Vert^{2}_{s,\alpha,s_1}$ instead of the norm $\Vert u \Vert^{2}_{s,\alpha}$. Theorem \ref{Th-tech-3} is now proven.  \finpv

\medskip

Once we have proven Theorems \ref{Th-tech-1}, \ref{Th-tech-2}, and \ref{Th-tech-3}, to conclude with the proof of the whole Theorem \ref{Th1} we must verify that  $u \in \mathcal{C}^{1}((0,T], \mathcal{C}^{\infty}(\R))$. Indeed, by Theorem \ref{Th-tech-3} we have $u \in \mathcal{C}((0,T], H^{\infty}(\R))$ and then for  $0<t\leq T$  the solution $u$ of the integral equation (\ref{Equation-Integral}) (constructed in Theorems \ref{Th-tech-1} and \ref{Th-tech-2}) also solves the differential  equation (\ref{Equation}) in the classical sense. We thus write  $\partial_t 	u  = - \,  D( \partial_x  u)  -  \Big( D^{\alpha}_{x} - D^{\beta}_{x} \Big)u  - \partial_x ( u^2)  - \partial^{2}_{x}(u^2) - (\partial_x u)^2$ to get that $\partial_t u \in \mathcal{C}(]0,T], H^{\infty}(\R))$. Thereafter, we can follow the same ideas at the end of the proof of \cite[Proposition $4.2$]{CorJar1}  to obtain $\partial_t u \in  \mathcal{C}(]0,T], \mathcal{C}^{\infty}(\R))$  and therefore  $u \in \mathcal{C}^{1}(]0,T], \mathcal{C}^{\infty}(\R))$.  
\subsubsection*{The case $\gamma_2=\gamma_3=0$ and $\alpha>2$.}
In this case, recall that  mild solutions of the equation (\ref{Equation}) write down as 
\begin{equation}
u(t,x)=K_{\alpha,\beta}(t,\cdot)\ast u_0(x) -\int_{0}^{t}K_{\alpha,\beta}(t-\tau,\cdot)\ast \partial_x (u^2)(\tau, x) d \tau.  
\end{equation}
As is the proof of Theorem \ref{Th-tech-1}, for $\alpha>\beta>0$ with $\alpha > 2$  and for  $-\frac{\alpha}{2}<s\leq0$ this equation is locally well-posed in the space (with $0<T<1$):
\begin{equation*}
 \mathcal{E}^{s,\alpha}_{T}=\{ u \in \mathcal{C}([0,T], H^s(\R)): \ \ \| u \|_{s,\alpha,2}<+\infty \},   
\end{equation*}
with the norm
\begin{equation*}
 \| u \|_{s,\alpha,2} = \sup_{0\leq t \leq T} \| u(t,\cdot)\|_{H^s}+ \sup_{0<t\leq T} t^{\frac{|s|}{\alpha}}\| u(t,\cdot)\|_{L^2}.\end{equation*}
 Indeed, we shall detail the bilinear estimates. We get back to the estimate (\ref{Estim-01}) (with $\sigma=1$) to obtain
\begin{equation}
\left\Vert \int_{0}^{t} K_{\alpha,\beta}(t-\tau,\cdot)\ast \partial_x (u^2)(\tau,\cdot)d \tau \right\Vert_{H^s} 
\lesssim \,   \left(  \int_{0}^{t}  (t-\tau)^{-\frac{s+1}{\alpha}-\frac{1}{2\alpha}}\, \tau^{-\frac{2 \vert s \vert }{\alpha}}\, d \tau\right)\,  \left( \sup_{0\leq t\leq T} t^{\frac{2\vert s \vert}{\alpha}} \Vert u(t, \cdot) \Vert^{2}_{L^2}\right).
\end{equation} 
To study the integral above, observe that since $s \leq 0$ we have  $-\frac{s+1}{\alpha}-\frac{1}{2\alpha}>-1$ as long as $\alpha>2$, and we have $-\frac{2 \vert s \vert }{\alpha}>-1$ as long as $-\frac{\alpha}{2}<s$. Then, we apply the estimate (\ref{Beta}) to obtain this integral computes down as $t^{1+\frac{s}{\alpha}-\frac{3}{2\alpha}}$. Moreover, remark that $1+\frac{s}{\alpha}-\frac{3}{2\alpha}>0$ as long as $s> \frac{3}{2}-\alpha$. Consequently, we have 
\begin{equation*}
\sup_{0\leq t \leq T} \left\Vert \int_{0}^{t} K_{\alpha,\beta}(t-\tau,\cdot)\ast \partial_x (u^2)(\tau,\cdot)d \tau \right\Vert_{H^s} \lesssim T^{1+\frac{s}{\alpha}-\frac{3}{2\alpha}}\, \| u \|^{2}_{s,\alpha,2}.    
\end{equation*}
Then, by the second point of Proposition \ref{Prop1} (with $\sigma=1$) we have 
\begin{equation*}
\sup_{0<t\leq T}t^{\frac{|s|}{\alpha}} \left\| \int_{0}^{t} K_{\alpha,\beta}(t-\tau,\cdot)\ast \partial_x (u^2)(\tau,\cdot)d \tau \right\|_{L^2} \lesssim T^{1+\frac{s}{\alpha}-\frac{3}{2\alpha}}\, \| u \|^{2}_{s,\alpha,2}.
\end{equation*}
By these estimates we have the local well-posed in the space $\mathcal{E}^{s,\alpha}_{T}$, with $\max\left( 3/2 - \alpha, -\alpha /2 \right)<s\leq 0$.

\medskip

Thereafter, in the case $s>0$, by following the same arguments in the proof of Theorem \ref{Th-tech-2} we also have the locally well-posedness in the space 
\begin{equation*}
  \mathcal{F}^{s,\alpha}_{T}= \{  u \in \mathcal{C}([0,T], H^s(\R)): \ \| u \|_{s,\alpha, s_1, 2} < +\infty \}, 
\end{equation*}
where for $\max\left( 3/2 - \alpha, -\alpha /2 \right)<s_1\leq 0$  we define 
\begin{equation*}
\Vert u \Vert_{s,\alpha,s_1,2}= \sup_{0\leq t \leq T} \Vert u(t,\cdot)\Vert_{H^s}+
\sup_{0\leq t \leq T} t^{\frac{\vert s_1 \vert}{\alpha}} \Vert u(t,\cdot) \Vert_{L^2}  + \sup_{0\leq t \leq T} t^{\frac{\vert s_1 \vert}{\alpha}} \Vert  \mathcal{J}^{s-s_1}\, u(t,\cdot) \Vert_{L^2}.
\end{equation*}
Finally, with minor modifications, the statement of Theorem \ref{Th-tech-3} (regularity of solutions) also holds in this case.  Theorem \ref{Th1} is now  proven. \finpv     
 
\subsubsection*{Proof of Theorem \ref{Th2}}  
Let us briefly explain the strategy of the proof.  We shall  assume that equation (\ref{Equation}) is locally well-posed in the space $H^s(\R)$ when $s<1-\alpha/2$ (when $\gamma_2,\gamma_3 \neq 0$) and $s<-\alpha/2$ (when $\gamma_2=\gamma_3=0$) respectively. Moreover,   we shall  assume that the flow-map function $S: H^s(\R) \to \mathcal{C}([0,T], H^s(\R))$ (defined in (\ref{Map-flow}))  is  a $\mathcal{C}^2-$ function at $u_0=0$. In particular,  this implies that the  second Fréchet derivative of $S(t)$ at $u_0=0$, defined as  $D^{2}_{0}S(t): H^s(\R)\times H^s(\R) \to  H^s(\R), \,\, (v_0, w_0)\mapsto D^{2}_{0}(v_0, w_0)$,  is a linear and bounded operator. Our general strategy is to construct  well-prepared initial data $(v_0, w_0) \in H^s(\R)\times H^s(\R) $ to contradict the boundness of the operator $D^{2}_{0}S(t)$.  The proof is divided into three steps: first, we shall explicitly compute the operator $D^{2}_{0}S(t)$. Then, we shall construct the well-prepared initial data and, in the last step, we shall prove the unboundedness of the operator $D^{2}_{0}S(t)$. 

\medskip

{\bf The operator $D^{2}_{0}S(t)$}. Our starting point is to explicitly compute this  operator. In all the computations   below, the limit is understood in the strong topology of the space $H^s(\R)$. Let us start by computing the first Fréchet derivative of $S(t)$ at $u_0 \in H^s(\R)$  in the direction $v_0\in H^s(\R)$.  Recall that by  (\ref{Map-flow}) and (\ref{Equation-Integral}), and moreover, for   the bilinear form $B(\cdot , \cdot)$  defined in (\ref{def-B}),  we have 
\[ S(t)u_0= u(t,\cdot)= K_{\alpha,\beta}(t,\cdot)\ast u_0  \underbrace{- \int_{0}^{t} K_{\alpha,\beta}(t-\tau,\cdot)\ast B(u,u)(\tau,\cdot) \, d\tau}_{\mathcal{B}(u,u)},\]
with 
\[ B(u,u)(\tau,\cdot)  = \Big(\gamma_1 \,\partial_x(u^2)+ \gamma_2 \, \partial^{2}_{x}(u^2)+\gamma_3 (\partial_x u)^2\Big)(\tau,\cdot).\]
For the sake of simplicity, we shall write  $\ds{ S(t)u_0= K_{\alpha,\beta}(t,\cdot)\ast u_0 + \mathcal{B}\Big(S(t)u_0, S(t)u_0 \Big)}$, where this bilinear form is symmetric.  Then, we have
\begin{equation}\label{Der-1}
\begin{split}
D_{u_0} S(t)v_0 =&\,  \lim_{h\to 0} \frac{S(t)(u_0+h\, v_0) - S(t)u_0}{h} \\
=& \, \lim_{h\to 0} \frac{K_{\alpha,\beta}(t,\cdot)\ast (u_0+hv_0) -K_{\alpha,\beta}(t,\cdot)\ast u_0}{h}\\
&\, + \lim_{h\to 0}\frac{\mathcal{B}\Big( S(t)(u_0 + h v_0), S(t)(u_0 + h v_0)  \Big) - \mathcal{B}\Big( S(t)(u_0), S(t)(u_0)  \Big)}{h}\\
=&\, K_{\alpha,\beta}(t,\cdot)\ast v_0+2 \mathcal{B}\Big( S(t)u_0, S(t)v_0 \Big).
\end{split}
\end{equation}
We compute now the second derivative  $D^{2}_{u_0}S(t)$ at $u_0=0$. First,  for $u_0, v_0 \in H^s(\R)$ fixed,  we define the function $ x  \in \R \mapsto  D_{x u_0}S(t) v_0 \in H^s(\R)$; and by following similar computations performed in (\ref{Der-1})  we have 
\begin{equation*}
\partial_{x} D_{x u_0}S(t) v_0=2 \mathcal{B}\Big( D_{xu_0}S(t)u_0, D_{xu_0}S(t)v_0 \Big)+ 2\mathcal{B}\Big( S(t)(xu_0), D^{2}_{xu_0}S(t)(u_0,v_0) \Big).
\end{equation*}
We thus set  $x=0$, and moreover,  by the identity   $S(t)0=0$ and since by (\ref{Der-1}) we have $D_0 S(t)v_0=K_{\alpha,\beta}(t,\cdot)\ast v_0$,  we obtain 
\begin{equation}\label{Der-2}
\begin{split}
D^{2}_{0}S(t)(u_0,v_0) = &\,  2 \mathcal{B}\Big( K_{\alpha,\beta}(t,\cdot)\ast u_0, \, K_{\alpha,\beta}(t,\cdot)\ast v_0 \Big)\\
= &\, 2\int_{0}^{t}K_{\alpha,\beta}(t-\tau,\cdot)\ast B\Big( K_{\alpha,\beta}(\tau,\cdot)\ast u_0, \, K_{\alpha,\beta}(\tau,\cdot)\ast v_0  \Big)\, d \tau. 
\end{split}
\end{equation}

\medskip

{\bf Well-prepared initial data}. Let $N \in \mathbb{N}^{*}$ be fixed such that $N \gg 1$. Moreover, let $r\in \R$ fixed such that $r \sim  1$. We consider  the disjoint intervals $[-N, -N+r]$ and $[N+r, N+2r]$. Then,  we define the functions $v_0$ and $w_0$ as 
\begin{equation}\label{initial-data}
 v_0 = r^{-1/2}\, N^{-s}\,  \mathcal{F}^{-1} \left( \mathds{1}_{[-N, -N+r]}(\xi)\right), \quad  v_0 = r^{-1/2}\, N^{-s}\,  \mathcal{F}^{-1} \left( \mathds{1}_{[N+r, N+2r]}(\xi)\right).
\end{equation}
We will verify that $\Vert v_0 \Vert_{H^s} \sim 1$ and  $\Vert w_0 \Vert_{H^s} \sim 1$. Indeed,  for the function $v_0$ defined above  we write 
\[ \Vert v_0 \Vert^{2}_{H^s}= \int_{\R} (1+\vert \xi \vert^2)^{s} r^{-1}\, N^{-2s} \mathds{1}_{I}(\xi) \, d \xi =  r^{-1} N^{-2s}\,\int_{-N}^{-N+r}   (1+\vert \xi \vert^2)^{s} d\xi.  \]
Here, as $\xi \in [-N, -N+r]$, and moreover, as $N\gg 1$ and  $r\sim 1$, we have $\vert \xi \vert \sim N$ and  $1+\vert \xi \vert^2 \sim  N^2$. Consequently,   $(1+\vert \xi \vert^2)^{s} \sim N^{2s}$.  We thus obtain
\[    r^{-1} N^{-2s}\,\int_{-N}^{-N+r}   (1+\vert \xi \vert^2)^{s} d\xi \sim r^{-1} N^{-2s} \, N^{2s}\, \int_{-N}^{-N+r}  d \xi =1. \]
The function $w_0$ follows the same estimates and we also have $\Vert w_0 \Vert_{H^s}\approx 1$. 

\medskip

{\bf The  unboundedness of the operator $D^{2}_{0}S(t)$}.   With the particular initial data constructed above, we shall prove the following estimates from below: for $N \in \mathbb{N}$ such that  $N \gg 1$ 
\begin{equation}\label{Estim}
\begin{cases}\vspace{2mm}
N^{2(1-s-\alpha/2)} \lesssim  \Vert D^{2}_{0} S(t) (u_0, v_0) \Vert_{H^s}, \quad  \mbox{when} \ \  s<1-\alpha/2, \\
N^{2(-s-\alpha/2)} \lesssim \Vert D^{2}_{0} S(t) (u_0, v_0) \Vert_{H^s}, \quad \mbox{when} \ \  s<-\alpha/2. 
\end{cases}
\end{equation}
Indeed, by the identity (\ref{Der-2}) we write 
\begin{equation}\label{Iden-aux}
\Vert D^{2}_{0} S(t) (u_0, v_0) \Vert_{H^s}=  2 \left\Vert (1+\vert \xi \vert^2)^{\frac{s}{2}}\, \mathcal{F}\left( \int_{0}^{t}K_{\alpha,\beta}(t-\tau,\cdot)\ast B\Big( K_{\alpha,\beta}(\tau,\cdot)\ast u_0, \, K_{\alpha,\beta}(\tau,\cdot)\ast v_0  \Big)\, d \tau \right) \right\Vert_{L^2}.
\end{equation}
Then, for $t>0$ and $\xi \in \R$ we define the function 
\[ g(t,\xi)= \mathcal{F}\left( \int_{0}^{t}K_{\alpha,\beta}(t-\tau,\cdot)\ast B\Big( K_{\alpha,\beta}(\tau,\cdot)\ast u_0, \, K_{\alpha,\beta}(\tau,\cdot)\ast v_0  \Big)\, d \tau \right)(\xi), \]
and recalling that $\widehat{K_{\alpha,\beta}}(t,\xi)= e^{-f(\xi) t}$  (with $f(\xi)=i\, m(\xi)\, \xi +(\vert \xi \vert^\alpha - \vert \xi \vert^\beta)$) we can prove the following identity which, for the reader's convenience, will be  postponed to Appendix \ref{AppendixA}:
\begin{equation}\label{Iden}
g(t,\xi)= \int_{\R} \left( \gamma_1\,  i \xi -\gamma_2\, \xi^2-\gamma_3(\xi - \eta)\eta \right) \widehat{u_0}(\xi-\eta)\widehat{v_0}(\eta)\, \frac{e^{-f(\eta)t-f(\xi-\eta)t}-e^{-f(\xi)t}}{f(\xi)-f(\eta)-f(\xi-\eta)}\, d \eta.
\end{equation}
Now, we prove the following estimate from below.
\begin{Lemme} The following estimates hold:
\begin{enumerate}
    \item When $\gamma_2,\gamma_3\neq 0$, we have $\ds{N^{2-2s-\alpha} \lesssim |g(t,\xi)|}$. 
    \item When $\gamma_2=\gamma_3=0$, we have $\ds{N^{-2s-\alpha} \lesssim |g(t,\xi)|}$. 
\end{enumerate}
\end{Lemme}
\pv  We must study each term inside the integral (\ref{Iden}). 
\begin{enumerate}
    \item Assume that $\gamma_2, \gamma_3 \neq 0$.  For the first and the second  term  we  have 
\begin{equation}\label{Sim1}
 \left( \gamma_1\, i \xi -\gamma_2\, \xi^2  -\gamma_3(\xi - \eta)\eta \right) \widehat{u_0}(\xi-\eta)\widehat{v_0}(\eta)   \sim   N^{2-2s}.    
\end{equation}
Indeed, recall that $u_0$ and $v_0$ are defined in (\ref{initial-data}) and we thus have $\widehat{u_0}(\xi-\eta)=r^{-1/2}N^{-s}\, \mathds{1}_{[-N, -N+r  ]}(\xi-\eta)$. Here $\xi-\eta \in [-N, -N+r]$ is equivalent to  $N-r+\xi \leq \eta \leq N+\xi$. Moreover, we also have  $\widehat{v_0}(\eta)= r^{-1/2}\, N^{-s}\, \mathds{1}_{[N+r, N+2r]}(\eta)$. Then, since  $N\in \mathbb{N}$ such that $N\gg 1$ and $r \sim 1$, we remark  that the intervals $N-r+\xi \leq \eta \leq N+\xi$ and $N+r \leq \eta \leq N+2r$ are not disjoint, provided that $r < \xi < 3 r$. Hence, we obtain $\xi \sim r$.  On the other hand, since  $N+r \leq \eta \leq N+2r$ we are able to write $\eta \sim  N$.  Consequently, we  obtain the estimates $ \left(\gamma_1 i \xi -\gamma_2\xi^2 -\gamma_3(\xi - \eta)\eta   \right) \widehat{u_0}(\xi-\eta)\widehat{v_0}(\eta)  \sim   \left( r + N^2 + r^2 \right) r^{-1}N^{-2s} \sim  N^{2}r^{-1} N^{-2s} \sim  N^{2-2s}$.  
\item Assume that $\gamma_2=\gamma_3=0$. By following the same arguments above we have 
\begin{equation}\label{Sim2}
 \left( \gamma_1 i \xi   \right) \widehat{u_0}(\xi-\eta)\widehat{v_0}(\eta)   \sim   N^{-2s}.    
\end{equation}
\end{enumerate}
On the other hand, for the third term, we  have the estimate 
\begin{equation}\label{Sim3}
f(\xi)-f(\eta)-f(\xi-\eta) \sim N^\alpha.     
\end{equation}
Indeed, recall that $f(\xi)=i\, m(\xi) \, \xi +(\vert \xi \vert^\alpha - \vert \xi \vert^\beta)$, where the symbol $m(\xi)$ is defined in  (\ref{def-D}) and  we have $i\, m(\xi) \, \xi  \sim  \xi^{3}$ when $D=\partial^{2}_{x}$; or $i\, m(\xi) \, \xi  \sim  \xi^{2}$ when $D=\mathcal{H}\partial_x$.   Moreover,  recall that $\xi \sim  1$ and $\eta \sim  N$. Then, since  $\alpha>\beta$ and $\alpha > 7/2$, for both cases  $D=\partial^{2}_{x}$ and  $D=\mathcal{H}\partial_x$ we have
$ f(\xi)-f(\eta)-f(\xi-\eta) \sim   f(\eta) \sim  N^\alpha.$

\medskip

With these estimates (\ref{Sim1}), (\ref{Sim2}), and (\ref{Sim3}) at hand, we get back to identity  (\ref{Iden}) to obtain the wished  estimates from below.  \finpv 

\medskip

Finally, we get back to the identity (\ref{Iden-aux}), hence we get the desired estimate (\ref{Estim}). In this   estimate, we consider first the case  $s<1-\alpha/2$, hence  we have $1-s-\alpha/2>0$. Moreover,  since   $\Vert  u_0 \Vert \sim  \Vert v_0 \Vert \sim  1$  we have $N^{2(1-s-\alpha/2)} \lesssim  \Vert D^{2}_{0} S(t) (u_0, v_0) \Vert_{H^s} \lesssim \Vert u_0 \Vert_{H^s}\, \Vert v_0 \Vert_{H^s} \lesssim 1$, which is a contradiction as long as  $N \gg 1$. Consequently, the flow-map function $S: H^s(\R) \to \mathcal{C}([0,T], H^s(\R))$ (given  in (\ref{Map-flow}))  is  not  a $\mathcal{C}^2-$ function at $u_0=0$.  The 
case $s<-\alpha/2$ follows the same ideas. Theorem \ref{Th2} is now proven. \finpv

\subsection{Global well-posedness}
\subsubsection*{Proof of Theorem \ref{Th3}}

As  the proof of Theorem \ref{Th1},  we shall consider  the two cases of the parameter $s$:  

\medskip 

{\underline{ \bf The case $1-\alpha/2 < s \leq 0$ (when $\gamma_2,\gamma_3\neq 0$) or $\max(3/2-\alpha,-\alpha/2)<s\leq 0$ (when $\gamma_2=\gamma_3=0$)}. 

\medskip

By Theorem \ref{Th-tech-3} the solution $u\in E^{s,\alpha}_{T}$ (constructed in Theorem \ref{Th-tech-1}) is regular enough and then, by multiplying the equation (\ref{Equation}) by $u(t,x)$, and after some integration by parts (in the spatial variable), we obtain
\begin{equation}\label{Global}
\frac{1}{2} \frac{d}{dt} \Vert u(t,\cdot) \Vert^{2}_{L^2}= - \int_{\R}  \Big( D^{\alpha}_{x} - D^{\beta}_{x} \Big)u \, u \,  dx - (2 \gamma_2 - \gamma_3) \int_{\R} (\partial_x u)^2\, u\, dx.
\end{equation}
Here we assume  that $-2 \gamma_2+\gamma_3=0$ to get 
\[ \frac{1}{2} \frac{d}{dt} \Vert u(t,\cdot) \Vert^{2}_{L^2}= - \int_{\R}  \Big( D^{\alpha}_{x} - D^{\beta}_{x} \Big)u \, u \,  dx. \]
We estimate now the term on the right-hand side. By the Parseval's identity and for $M=2^{\frac{1}{\alpha-\beta}}$ (note that for $ \vert \xi \vert \geq M$ we have $\vert \xi \vert^\beta - \vert \xi \vert^\alpha \leq - \vert \xi \vert^\beta$) we write
\begin{equation}\label{Estim-Dissipative}
\begin{split}
&- \int_{\R}  \Big( D^{\alpha}_{x} - D^{\beta}_{x} \Big)u \, u \,  dx =\, - \int_{\R} \Big( \vert \xi \vert^\alpha - \vert \xi \vert^\beta \Big) \vert \widehat{u}\vert^2 \, d \xi  \\
=&\, \int_{\vert \xi \vert \leq M} \Big( \vert \xi \vert^\beta - \vert \xi \vert^\alpha \Big) \vert \widehat{u}\vert^2 \, d \xi + \int_{\vert \xi \vert \geq M} \Big( \vert \xi \vert^\beta - \vert \xi \vert^\alpha \Big) \vert \widehat{u}\vert^2 \, d \xi \leq  \, \int_{\vert \xi \vert \leq M} \vert \xi \vert^\beta \vert \widehat{u} \vert^2 \, d \xi   - \int_{\vert \xi \vert \geq M} \vert \xi \vert^\beta \vert \widehat{u}\vert^2\, d \xi \\
\leq &\,  \int_{\vert \xi \vert \leq M} \vert \xi \vert^\beta \vert \widehat{u} \vert^2 \, d \xi  \leq M^\beta  \Vert u \Vert^{2}_{L^2}. 
\end{split}
\end{equation}
With this estimate and by the Gr\"onwall inequality, for all $ \frac{T}{2}<t<T$ we obtain
\[ \Vert u(t,\cdot) \Vert^{2}_{H^s} \leq  \Vert u(t,\cdot) \Vert^{2}_{L^2} \leq \Vert u( T / 2, \cdot) \Vert^{2}_{L^2} e^{2 M^\beta\, t}. \]
hence, the solution can be extended to the whole interval $[0, +\infty[$. 

\medskip

{\underline{ \bf The case $0< s$}}. We shall follow similar ideas of \cite[Proposition $4.3$]{CorJar1}. Let   $u_0 \in H^s(\R)$ (with $s>1-\alpha/2$ or $s>\max(3/2-\alpha,-\alpha/2)$ respectively) be an initial datum.  We define the time $T^{*}$ as follows:
\[ T^*= \sup\left\{  T>0: \,\, \text{there exists a unique solution}\,\, u \in \mathcal{C}([0,T], H^s(\R)) \,\, \text{of}\,\, (\ref{Equation})\,\, \text{arising from}\,\, u_0\right\}.\]
We assume  the relationship  $-2 \gamma_2+\gamma_3=0$  and we will prove that $T^*=+\infty$. Our strategy is to assume that  $T^*<+\infty$ to obtain a contradiction. Always by Theorem \ref{Th-tech-3} and by following the same estimates above, we have the estimate 
\begin{equation}\label{Gronwall}
\Vert u(t,\cdot) \Vert^{2}_{L^2} \leq \Vert u_0 \Vert^{2}_{L^2} e^{2 M^\beta\, T^{*}}. 
\end{equation}
where we set the constant $M_0= \Vert u_0 \Vert^{2}_{L^2} e^{2 M^\beta\, T^{*}}>0$. 

\medskip

On the other hand, recall that by Theorem \ref{Th-tech-2} for any initial datum $v_0 \in H^s(\R)$  there exists  $v \in F^{s,\alpha}_{T} \subset  \mathcal{C}([0,T], H^s(\R))$ an arising solution of the equation (\ref{Equation}), where the time $T=T(v_0)$ is given by the expression (\ref{Control-Tiempo-LWP}). Precisely, we have the bound from above   $T(v_0)< \frac{1}{4^{1/\eta} \Vert v_0 \Vert^{1/\eta}_{H^s}}$ and since   $\Vert v_0 \Vert_{L^2}\leq \Vert v_0 \Vert_{H^s}$, we obtain $T(v_0)< \frac{1}{4^{1/\eta} \Vert v_0 \Vert^{1/\eta}_{L^2}}$.  Consequently, the time $T(v_0)$ is a decreasing function of $\Vert v_0 \Vert_{L^2}$.  This decreasing property yields that we can find  
 a time $0<T_1 < T^{*}$ such that for all initial datum $v_0 \in H^{s}(\R)$  verifying $\Vert v_0 \Vert_{L^2} \leq M_0$ the associated solution $v \in  \mathcal{C}([0, T[, H^s(\R))$ exists at least on the interval  $[0,T_1]$; and it verifies  $v \in \mathcal{C}([0,T_1], L^2(\R))$.  
 
 \medskip

 Now, for $0< \varepsilon < T_1$ and for the solution $u(t,x)$ (arising from $u_0$)   we  the initial datum $v_0= u(T^{*}-\varepsilon, \cdot) \in H^s(\R)$, which by (\ref{Gronwall}) verifies $\Vert v_0 \Vert_{L^2}\leq M_0$. So, there exists a solution $v$ arising from  $v_0= u(T^{*}-\varepsilon, \cdot)$ which is defined at least on $[0,T_1]$. Thus, by  gathering  the functions  $u(t,x)$ and $v(t,x)$  we get a solution   
\begin{equation*}
\tilde{u}(t,\cdot) = \begin{cases}\vspace{1mm} u(t,\cdot),  \quad  \text{when} \quad 0\leq t \leq T^{*}-\varepsilon, \\
v(t,\cdot),  \quad   \text{when}\quad T^{*}-\varepsilon \leq t \leq T^{*} - \varepsilon +T_1, \end{cases}
\end{equation*}
which arises from the datum $u_0$  and which  is defined on the interval $[0, T^{*} - \varepsilon +T_1]$. But, since $0<\varepsilon <T_1$ we have $T^{*} - \varepsilon +T_1 > T^{*}$, which contradicts the definition of the time $T^{*}$.  We thus have $T^{*}=+\infty$. Theorem \ref{Th3} is now proven.  \finpv 

\subsubsection*{Proof of Proposition \ref{Prop-blow-up}} 
Our starting point is the identity (\ref{Global}), where we must estimate the second  term  on the right-hand side:  $\ds{\int_{\R} (\partial_x u)^2\, u\, dx}$. We write 
\begin{equation*}
\begin{split}
&\, \int_{\R} (\partial_x u)^2\, u\, dx =  \int_{\R} \partial_x u \, u\, \partial_x u \, dx = \frac{1}{2}\, \int_{\R} \partial_x u \partial_x (u^2)\, dx = - \frac{1}{2} \int_{\R} \partial^{2}_{x} u \, u^2\, dx \\
\leq & \,  \Vert \partial^{2}_{x} u \Vert_{L^\infty}\, \Vert u^2 \Vert_{L^1} \leq \Vert \partial^{2}_{x} u \Vert_{L^\infty}\, \Vert u \Vert^{2}_{L^2}. 
\end{split}
\end{equation*}
With this estimate at hand, we get back to (\ref{Global}), hence, together with the  estimate (\ref{Estim-Dissipative}) we get
\begin{equation}
\begin{split}
\frac{1}{2}\frac{d}{dt}\Vert u(t,\cdot)\Vert^{2}_{L^2} \leq &\, M^\beta \, \Vert u(t,\cdot)\Vert^{2}_{L^2}+ \frac{2\gamma_2-\gamma_3}{2}\, \Vert \partial^{2}_{x} u(t,\cdot)\Vert_{L^\infty}\, \Vert u(t,\cdot) \Vert^{2}_{L^2}\\
\lesssim &\, (1+ \Vert \partial^{2}_{x} u(t,\cdot) \Vert_{L^\infty})\, \Vert u(t,\cdot)\Vert^{2}_{L^2}. 
\end{split}
\end{equation}
Then, by the Gr\"onwall inequality  for all  $t>T/2$ (where the time $T$ is given by (\ref{Control-Tiempo-LWP})) we have
\begin{equation*}
\Vert u(t,\cdot) \Vert^{2}_{L^2} \lesssim \Vert u(T/ 2,\cdot) \Vert^{2}_{L^2}\, e^{t-T/2 + \int_{T/2}^{t} \Vert \partial_x u(s,\cdot) \Vert_{L^\infty} ds},
\end{equation*}
hence we obtain 
\begin{equation}\label{Gronwall-2}
\Vert u(t,\cdot) \Vert^{2}_{L^2} \lesssim \Vert u(T/ 2,\cdot) \Vert^{2}_{L^2}\, e^{t + \int_{0}^{t} \Vert \partial_x u(s,\cdot) \Vert_{L^\infty} ds}.
\end{equation}
From this estimate the blow-up criterion stated in Proposition \ref{Prop-blow-up}  is obtained as follows: first, let us assume that $\ds{ \lim_{t\to T^*} \Vert u(t,\cdot)\Vert_{H^s}=+\infty}$. This fact yields   $\int_{0}^{T^*}\Vert \partial^{2}_{x} u(t,\cdot) \Vert_{L^\infty} dt = +\infty$. Indeed, if we assume that  $\int_{0}^{T^*}\Vert \partial^{2}_{x} u(t,\cdot) \Vert_{L^\infty} dt < +\infty$, then by (\ref{Gronwall-2}) we get that the quantity $\Vert u(t,\cdot)\Vert^{2}_{L^2}$ can be extended beyond the time $T^*$; and by following the same arguments in the proof of Theorem \ref{Th3}, we get that the quantity $\Vert u(t,\cdot) \Vert_{H^s}$ extends beyond $T^*$, which contradicts  the definition of $T^*$.  

\medskip

Now, let us assume that $\int_{0}^{T^*}\Vert \partial^{2}_{x} u(t,\cdot) \Vert_{L^\infty} dt = +\infty$, which yields $\lim_{t\to T^*} \Vert u(t,\cdot)\Vert_{H^s}=+\infty$. Indeed, if we  assume that  $\ds{\lim_{t\to T^*} \Vert u(t,\cdot)\Vert_{H^s}<+\infty}$, then  by Theorem \ref{Th-tech-3} we have $u \in  \mathcal{C}(0, T^*+\varepsilon], H^{\infty} (\R))$ with $\varepsilon>0$.  Consequently, for $\sigma >  1/2 $ we have $\int_{0}^{T^*} \Vert u(t,\cdot) \Vert_{H^{2+\sigma}} dt <+\infty$. Then, by the continuous Sobolev embedding $L^\infty(\R) \subset H^\sigma (\R)$, we write
\[ \int_{0}^{T^*} \Vert \partial^{2}_{x} u(t,\cdot)\Vert_{L^{\infty}} dt \leq \int_{0}^{T^*} \Vert \partial^{2}_{x} u(t,\cdot)\Vert_{H^\sigma} dt  \leq \int_{0}^{T^*} \Vert  u(t,\cdot)\Vert_{H^{2+\sigma }} dt <+\infty,  \]
which is a contradiction.  We thus have  $\ds{\lim_{t\to T^*} \Vert u(t,\cdot)\Vert_{H^s}=+\infty}$. Proposition \ref{Prop-blow-up} is proven. \finpv

\section{Spatially decaying properties}\label{Sec-Spatially-decaying}
\subsection{Kernel estimates II}
\subsubsection*{Proof of Proposition \ref{Prop-Principal-Kernel}} 
Since the definition of this kernel involves the operator $D$ defined in (\ref{def-D}), we shall consider   the following  cases: when  $D=\mathcal{H}\partial_x$, we shall refer to the nonlocal dispersive effects due the presence of  the Hilbert transform. On the other hand, when $D=\partial^{2}_{x}$ we shall refer to the local dispersive effects. Moreover, recall that the action of the operator $D$ is given in the Fourier level by the symbol $m(\xi)$, which is also given in the expression (\ref{def-D}).  Then, for the sake of clearness, we shall prove the identity (\ref{Kernel-pointwise}) in the following technical propositions.

\begin{Proposition}[The nonlocal dispersive effects]\label{Prop-Kernel-1} Let $D = \mathcal{H}\partial_x$, where we have $m(\xi)=\vert \xi \vert$. Let $\alpha > \beta \geq 1$  with $\alpha>2$. For $t>0$ there exists a quantity $I(t)$, which verifies $\vert I(t)\vert \leq C e^{\eta_1\, t}$ with $C>0$ and $\eta_1>0$ depending on $\alpha$ and $\beta$, such that for all $x\neq 0$ the following estimate holds: 
\[\ds{\vert K_{\alpha, \beta}(t,x) \vert = \frac{\vert I(t) \vert}{\vert x \vert^{\min(3, [\beta]+1)}}},\] where $[\beta]$ denotes the integer part of $\beta$.		
\end{Proposition}
\pv  We start by explaining the general idea of the proof. This idea was inspired by the previous works \cite{CorJar0,CorJar1}.  By the expression (\ref{kernel}),  for $x \neq 0$ and $t>0$ we write 
\begin{equation*}
K_{\alpha,\beta}(t,x) =  \int_{-\infty}^{0} e^{2 \pi i\, x \xi } e^{- f(\xi) t}  d \xi + \int_{0}^{+\infty} e^{2 \pi i\, x \xi } e^{- f(\xi) t}  d \xi. 
\end{equation*}
In each term on the right-hand side,  we multiply and divide by $2 \pi  i\,x$ to get:
\begin{equation*}
\begin{split}
K_{\alpha,\beta}(t,x) =  \frac{1}{2 \pi i x} \int_{-\infty}^{0}  2 \pi i x\,  e^{2 \pi i\, x \xi } e^{- f(\xi) t}  d \xi +  \frac{1}{2 \pi i x} \int_{0}^{+\infty}   2 \pi i x\, e^{2 \pi i\, x \xi } e^{- f(\xi) t}  d \xi.
\end{split}
\end{equation*}
Then, since  $2 \pi i  x \, e^{2 \pi i \, x \xi}= \partial_\xi \left( e^{2 \pi i \, x\xi}\right)$ we  integrate by parts respect to the variable $\xi$ to obtain
\begin{equation*}
\begin{split}
K_{\alpha, \beta}(t,x)= &\,\frac{1}{2\pi i x} \left(-1\,  \int_{-\infty}^{0} e^{2 \pi i \, x\xi} \, \partial_{\xi} \left(e^{-f(\xi) t} \right) d \xi + \left. \left( e^{2 \pi i \, x\xi}\,  \left( e^{-f(\xi) t}\right)\right) \right\vert^{0}_{-\infty}   \right) \\
&\, + \frac{1}{2\pi i x} \left(-1\,  \int_{0}^{+\infty} e^{2 \pi i \, x\xi} \, \partial_{\xi} \left(e^{-f(\xi) t} \right) d \xi + \left. \left( e^{2 \pi i \, x\xi}\,  \left( e^{-f(\xi) t}\right)\right) \right\vert^{+\infty}_{0}  \right).
\end{split}
\end{equation*}
By iterating this process $n$ times, we formally obtain the following expression:  
\begin{equation}\label{Iden-K-n}
\begin{split}
K_{\alpha, \beta}(t,x)= &\,\frac{1}{(2\pi i x)^n} \left((-1)^n\,  \int_{-\infty}^{0} e^{2 \pi i \, x\xi} \, \partial^{n}_{\xi} \left(e^{-f(\xi) t} \right) d \xi + \left. \left( e^{2 \pi i \, x\xi}\, \partial^{n-1}_{\xi} \left( e^{-f(\xi) t}\right)\right) \right\vert^{0}_{-\infty}   \right) \\
&\, + \frac{1}{(2\pi i x)^n} \left((-1)^n\,  \int_{0}^{+\infty} e^{2 \pi i \, x\xi} \, \partial^{n}_{\xi} \left(e^{-f(\xi) t} \right) d \xi + \left. \left( e^{2 \pi i \, x\xi}\, \partial^{n-1}_{\xi} \left( e^{-f(\xi) t}\right)\right) \right\vert^{+\infty}_{0}  \right),
\end{split}
\end{equation}
and this iterative process continues until we have one of the following scenarios:  
\begin{enumerate}
\item[$\bullet$] On the one hand, this process  stops at the step  $n$ when  for the next step $n+1$ we have 
\begin{equation}\label{Int-diverge}
\int_{-\infty}^{0} e^{2 \pi i \, x\xi} \, \partial^{n+1}_{\xi} \left(e^{-f(\xi) t} \right) d \xi = +\infty, \quad  \int_{0}^{+\infty} e^{2 \pi i \, x\xi} \, \partial^{n+1}_{\xi} \left(e^{-f(\xi) t} \right) d \xi=+\infty.
\end{equation}
Precisely, when  both integrals diverge at $\xi =0$ depending on  the behavior  of the function $f^{(n+1)}(\xi)$ when $\xi \to  0^{-}$ and $\xi \to 0^{+}$. 

\medskip

\item[$\bullet$] On the other hand, this process stops at the step $n$ when we have 
\begin{equation}\label{Limite}
L_n:=\left. \left( e^{2 \pi i \, x\xi}\, \partial^{n-1}_{\xi} \left( e^{-f(\xi) t}\right)\right) \right\vert^{0}_{-\infty} + \left. \left( e^{2 \pi i \, x\xi}\, \partial^{n-1}_{\xi} \left( e^{-f(\xi) t}\right)\right) \right\vert^{+\infty}_{0}  \neq 0. 
\end{equation}
We thus obtain 
\begin{equation}\label{Identity-Kernel}
K_{\alpha,\beta}(t,x) = \frac{ L_n }{(2\pi \, i x )^n} + \frac{(-1)^n}{2\pi \, i x )^n}  \int_{-\infty}^{+\infty} e^{2\pi i\, x \xi} \partial^{n}_{\xi}\left( e^{-f(\xi)t } \right) d \xi. \end{equation}
\end{enumerate} 
In both scenarios,   we conclude the wished identity stated in this proposition: $\vert K_{\alpha,\beta}(t,x)\vert = \frac{\vert I(t) \vert}{\vert x \vert^n}$. The generic quantity $I(t)$ may change in the different cases that we shall consider below, but we always have the control $\vert I(t)\vert \leq C e^{\eta_1\, t}$. 

\medskip 

Now, we can prove this proposition. For the sake of clearness, we shall consider separately the following cases of the parameter $\beta \geq 1$. 
\begin{enumerate}

\item[$\bullet$] The case $1\leq \beta <2$. By the expression (\ref{Identity-Kernel})   (with $n=2$) we write
\begin{equation*}
 K_{\alpha,\beta}(t,x)  = \frac{L_2}{(2\pi\, i x)^2} +  \frac{1}{(2\pi i x)^2} \left(  \int_{-\infty}^{+\infty} e^{2 \pi i x \xi } \partial^{2}_{\xi} \left( e^{-f(\xi) t } \right)  d \xi  \right).
\end{equation*}
Recall that the term $L_2$ (given in (\ref{Limite})) involves the expression $f'(\xi)$; and by a simple computation we have: 
\begin{equation}\label{1-der}
f'(\xi)= \begin{cases}\vspace{2mm} - 2  i \xi  - \alpha (-\xi)^{\alpha-1} + \beta (-\xi)^{\beta -1}, \quad \xi <0, \\
2 i \xi   + \alpha \xi^{\alpha-1} - \beta \xi^{\beta -1}, \quad \xi >0.
\end{cases}
\end{equation}
Thus, when $\beta=1$ by this expression we obtain $L_2 = - 2 t \neq 0$, and thus we can write  
\begin{equation*}
\vert K_{\alpha,\beta}(t,x)\vert= \frac{1}{\vert 2\pi \, i x \vert^2} \left\vert -2 t + \int_{-\infty}^{+\infty} e^{2 \pi i x \xi } \partial^{2}_{\xi} \left( e^{-f(\xi) t } \right)  d \xi \right\vert = \frac{\vert I(t)\vert}{\vert x \vert^2}=  \frac{\vert I(t)\vert}{\vert x \vert^{\min(3, [\beta]+1)}},  
\end{equation*}
where $[\beta]=1$.  Moreover, by the good decaying properties of the function $e^{-f(\xi)t}$  and  by following the same computations performed in \cite[Lemma $3.1$]{CorJar0}, we have $\vert I(t) \vert \leq C^{\eta_1\, t}$.

\medskip

On the other hand, when $1<\beta < 2$ by the identity (\ref{1-der}) we have $L_2=0$ and we obtain
\begin{equation*}
  K_{\alpha,\beta}(t,x)  =   \frac{1}{(2\pi i x)^2} \left(  \int_{-\infty}^{+\infty} e^{2 \pi i x \xi } \partial^{2}_{\xi} \left( e^{-f(\xi) t } \right)  d \xi  \right).   
\end{equation*}
To study the integral above we need to compute $f''(\xi)$ and we have: 
\begin{equation}\label{2-der}
f''(\xi)= \begin{cases}\vspace{2mm} - 2  i   +  \alpha(\alpha-1) (-\xi)^{\alpha-2} - \beta(\beta-1) (-\xi)^{\beta -2}, \quad \xi <0, \\
2 i   + \alpha(\alpha-1) \xi^{\alpha-2} - \beta(\beta-1) \xi^{\beta -2}, \quad \xi >0.
\end{cases}
\end{equation}
In particular, we have $f''(\xi) \sim \xi^{\beta-2}$ when $\xi \to 0$ and since $1<\beta < 2$ this integral converges. 

\medskip

Finally, we remark that  for the next value $n=3$, by the expression  (\ref{Identity-Kernel}) we formally  have 
\begin{equation*} 
 K_{\alpha,\beta}(t,x)  = \frac{L_3}{(2\pi\, i x)^3} +  \frac{1}{(2\pi i x)^2} \left(  \int_{-\infty}^{+\infty} e^{2 \pi i x \xi } \partial^{3}_{\xi} \left( e^{-f(\xi) t } \right)  d \xi  \right),  
\end{equation*}
but the last integral diverges. Indeed, to study this integral we need to compute $f'''(\xi)$ and we have
\begin{equation}\label{n-der}
f^{'''}(\xi)= \begin{cases}\vspace{2mm} - c_\alpha (-\xi)^{\alpha -3} + c_\beta (-\xi)^{\beta -3}, \quad \xi <0, \\
c_\alpha \, \xi^{\alpha-3} - c_{\beta} \, \xi^{\beta-3}, \quad \xi >0.
\end{cases}
\end{equation}
By the expression (\ref{n-der}) we observe that  $f^{(3)}(\xi)\sim \xi^{\beta-3}$ when $\xi\to 0$, and since $1<\beta<2$ this fact yields (\ref{Int-diverge}). We thus obtain $\vert K_{\alpha,\beta}(t,x)\vert = \frac{\vert I(t)\vert}{\vert x \vert^{\min(3, [\beta]+1)}}$. 

\item[$\bullet$] The case $2\leq \beta$.  By  (\ref{Identity-Kernel}) (with $n=3$) we   can write:  
\begin{equation*}
 K_{\alpha,\beta}(t,x)  = \frac{L_3}{(2\pi\, i x)^3} +  \frac{1}{(2\pi i x)^2} \left(  \int_{-\infty}^{+\infty} e^{2 \pi i x \xi } \partial^{3}_{\xi} \left( e^{-f(\xi) t } \right)  d \xi  \right). 
\end{equation*} 
where by the expressions (\ref{Limite}) and (\ref{2-der})   we always  have $L_3= -4 i\, t \neq 0$. Then, by the identity (\ref{Identity-Kernel}) (with $n=3$)  we  obtain 
\begin{equation*}
\vert K_{\alpha,\beta}(t,x)\vert = \left\vert  \frac{ -4 i \, t  }{(2\pi \, i x )^3} + \frac{-1}{2\pi \, i x )^3}  \int_{-\infty}^{+\infty} e^{2\pi i\, x \xi} \partial^{3}_{\xi}\left( e^{-f(\xi)t } \right)  \right\vert = \frac{\vert I(t) \vert}{\vert x \vert^3}=\frac{\vert I(t) \vert}{\vert x \vert^{\min(3, [\beta]+1)}}.    \end{equation*}
\end{enumerate}  
Proposition  \ref{Prop-Kernel-1} is proven.  \finpv

  

\begin{Proposition}[The local dispersive effects I ]\label{Prop-Kernel-2} Let $D = \partial^{2}_{x}$, where we have $m(\xi)=-\vert \xi \vert^2$. Let $\alpha > \beta \geq 1$ with $\alpha>2$. Moreover, we assume that  $\alpha$ and $\beta$ are both even numbers. Then we have $ K_{\alpha,\beta}(t,\cdot) \in \mathcal{S}(\R)$. 
\end{Proposition}
\pv  Since  $\alpha$ and $\beta$ are both even numbers, and moreover, since  $m(\xi)=-\vert \xi \vert^2$, the function  $f(\xi)$ given in (\ref{kernel}) verifies $f \in \mathcal{C}^{\infty}(\R)$. Consequently, by the good decaying properties of the function $e^{-f(\xi)t}$ when $\vert \xi \vert \to +\infty$, we have $e^{-f(\xi) t} \in \mathcal{S}(\R)$. Then, always by (\ref{kernel}) we conclude that $K_{\alpha, \beta}(t, \cdot) \in \mathcal{S}(\R)$. \finpv 

\medskip

\begin{Proposition}[The local dispersive effects II ]\label{Prop-Kernel-3} Let $D = \partial^{2}_{x}$, where we have $m(\xi)=-\vert \xi \vert^2$. Let $\alpha > \beta \geq 1$ with $\alpha>2$. Moreover, we assume that  $\alpha$ and $\beta$ are not \emph{both} even numbers. 

\medskip

For $t>0$ there exists a  quantity $I(t)$  which verifies $\vert I(t)\vert \leq C_1 e^{\eta_1\, t}$  with the constants  $C, \eta_1> 0$  depending on $\alpha$ and $\beta$, such that   for $x \neq 0$ the following estimates hold:
\begin{enumerate}
\item[1] If $\beta>0$ is not an even number, we have  $\ds{\vert K_{\alpha,\beta}(t,x)\vert = \frac{\vert I(t) \vert}{\vert x \vert^{[\beta]+1}}}$, 

\medskip

 \item[2] If $\beta>0$ is an even number, we have  $   \ds{\vert K_{\alpha,\beta}(t,x)\vert = \frac{\vert I(t) \vert}{\vert x \vert^{[\alpha]+1}}}$,
\end{enumerate} 
 where $[\beta]$ and $[\alpha]$ denote the integer part of $\alpha$ and $\beta$ respectively.   
   \end{Proposition}
\pv  The proof essentially follows the same ideas as the proof of Proposition \ref{Prop-Kernel-1}. Our starting point is the expression (\ref{Iden-K-n}) with the value $n=[\beta]+1$: 
\begin{equation}\label{Iden-K-beta}
\begin{split}
K_{\alpha, \beta}(t,x)= &\,\frac{1}{(2\pi i x)^{[\beta]+1}} \left((-1)^{[\beta]+1}\,  \int_{-\infty}^{0} e^{2 \pi i \, x\xi} \, \partial^{[\beta]+1}_{\xi} \left(e^{-f(\xi) t} \right) d \xi + \left. \left( e^{2 \pi i \, x\xi}\, \partial^{[\beta]}_{\xi} \left( e^{-f(\xi) t}\right)\right) \right\vert^{0}_{-\infty}   \right) \\
&\, + \frac{1}{(2\pi i x)^{[\beta]+1}} \left((-1)^{[\beta]+1}\,  \int_{0}^{+\infty} e^{2 \pi i \, x\xi} \, \partial^{[\beta]+1}_{\xi} \left(e^{-f(\xi) t} \right) d \xi + \left. \left( e^{2 \pi i \, x\xi}\, \partial^{[\beta]}_{\xi} \left( e^{-f(\xi) t}\right)\right) \right\vert^{+\infty}_{0}  \right).
\end{split}
\end{equation}
To study  the terms 
\[ \left. \left( e^{2 \pi i \, x\xi}\, \partial^{[\beta]}_{\xi} \left( e^{-f(\xi) t}\right)\right) \right\vert^{0}_{-\infty}, \qquad  \left. \left( e^{2 \pi i \, x\xi}\, \partial^{[\beta]}_{\xi} \left( e^{-f(\xi) t}\right)\right) \right\vert^{+\infty}_{0},\]
we need to compute the expression $f^{([\beta])}$ and we have 
\begin{equation}\label{n-der-2}
f^{([\beta])}(\xi)= \begin{cases}\vspace{2mm} i\,c_\beta  \, \xi^{3-[\beta]}+ (-1)^{[\beta]} c_\alpha (-\xi)^{\alpha -[\beta]} - (-1)^{[\beta]} c_\beta (-\xi)^{\beta - [\beta]}, \quad \xi <0, \\
i\,c_\beta \, \xi^{3-[\beta]}+ c_\alpha \, \xi^{\alpha-[\beta]} - c_{\beta} \, \xi^{\beta-[\beta]}, \quad \xi >0, 
\end{cases} \qquad  1\leq [\beta] \leq 3,
\end{equation}
\begin{equation}\label{n-der-3}
f^{[\beta]}(\xi)= \begin{cases}\vspace{2mm}  (-1)^{[\beta]} c_\alpha (-\xi)^{\alpha -[\beta]} - (-1)^{[\beta]} c_\beta (-\xi)^{\beta -[\beta]}, \quad \xi <0, \\
 c_\alpha \, \xi^{\alpha-[\beta]} - c_{\beta} \, \xi^{\beta-[\beta]}, \quad \xi >0, 
\end{cases} \qquad 4 \leq [\beta]. 
\end{equation} 
At this point, we shall consider the following cases of  parameter $\beta$.
\begin{enumerate}
\item When $\beta$ is not an even number. Within this setting, we still need to consider the next sub-cases.
\begin{itemize}
    \item[1.1] When $\beta$ is an integer number. We thus have $\beta=[\beta]$ where $[\beta]$ is not an even number. Then, by the expressions (\ref{n-der-2}) and (\ref{n-der-3}) we have 
    \begin{equation*}
  \left. \left( e^{2 \pi i \, x\xi}\, \partial^{[\beta]}_{\xi} \left( e^{-f(\xi) t}\right)\right) \right\vert^{0}_{-\infty}+ \left. \left( e^{2 \pi i \, x\xi}\, \partial^{[\beta]}_{\xi} \left( e^{-f(\xi) t}\right)\right) \right\vert^{+\infty}_{0} =-2 c_\beta\, t \neq 0.      
    \end{equation*}
Consequently, we obtain 
\begin{equation*}
\begin{split}
\vert K_{\alpha,\beta}(t,x)\vert= &\, \left\vert \frac{-2c_\beta \, t}{(2\pi \, i x)^{[\beta]+1}}+ \frac{(-1)^{[\beta]+1}}{(2\pi \, i x)^{[\beta]+1}}\int_{-\infty}^{+\infty} e^{2 \pi i \, x\xi} \,\partial^{[\beta]+1}_{\xi} \left( e^{-f(\xi)t}\right) d \xi \right\vert\\
=&\, \frac{\vert I(t)\vert}{\vert x \vert^{[\beta]+1}}, \qquad \vert I(t)\vert \leq C e^{\eta\, t}. 
\end{split}
\end{equation*}
\item[1.2] When $\beta$ is not an integer number. In this case, we have $[\beta]<\beta$ and by the expressions (\ref{n-der-2}) and (\ref{n-der-3}), we obtain
    \begin{equation*}
  \left. \left( e^{2 \pi i \, x\xi}\, \partial^{[\beta]}_{\xi} \left( e^{-f(\xi) t}\right)\right) \right\vert^{0}_{-\infty}+ \left. \left( e^{2 \pi i \, x\xi}\, \partial^{[\beta]}_{\xi} \left( e^{-f(\xi) t}\right)\right) \right\vert^{+\infty}_{0} =0.      \end{equation*}
Therefore, we can write 
\begin{equation*}
\begin{split}
\vert K_{\alpha,\beta}(t,x)\vert= &\, \left\vert \frac{(-1)^{[\beta]+1}}{(2\pi \, i x)^{[\beta]+1}}\int_{-\infty}^{+\infty} e^{2 \pi i \, x\xi} \,\partial^{[\beta]+1}_{\xi} \left( e^{-f(\xi)t}\right) d \xi \right\vert = \frac{\vert I(t)\vert}{\vert x \vert^{[\beta]+1}}, \quad \vert I(t)\vert \leq C e^{\eta_1\, t}. 
\end{split}
\end{equation*}
Moreover, we remark that we cannot continue the iterative process described in (\ref{Iden-K-n}) with the next step $n=[\beta]+2$: the resulting integrals involve the expression $\partial^{[\beta]+2}_{\xi}\left(e^{-f(\xi)t}\right)$, which contains the term $f^{([\beta]+2)}(\xi)$. But,  by a simple calculation, from  the expressions   (\ref{n-der-2}) and (\ref{n-der-3}) for $[\beta]+2=3$ we have
\begin{equation*}
f^{([\beta]+2)}(\xi)= \begin{cases}\vspace{2mm} i\,c_\beta  \, \xi^{3-([\beta]+2)}+ (-1)^{[\beta]} c_\alpha (-\xi)^{\alpha -[\beta]-2} - (-1)^{[\beta]+2} c_\beta (-\xi)^{\beta - [\beta]-2}, \quad \xi <0, \\
i\,c_\beta \, \xi^{3-([\beta]+2)}+ c_\alpha \, \xi^{\alpha-[\beta]-2} - c_{\beta} \, \xi^{\beta-[\beta]-2}, \quad \xi >0, 
\end{cases}
\end{equation*}
and for $4 \leq [\beta]+2$ we have 
\begin{equation*}
f^{([\beta]+2)}(\xi)= \begin{cases}\vspace{2mm}  (-1)^{[\beta]+2} c_\alpha (-\xi)^{\alpha -[\beta]-2} - (-1)^{[\beta]+2} c_\beta (-\xi)^{\beta -[\beta]-2}, \quad \xi <0, \\
 c_\alpha \, \xi^{\alpha-[\beta]-2} - c_{\beta} \, \xi^{\beta-[\beta]-2}, \quad \xi >0, 
\end{cases} 
\end{equation*} 
In both cases, we get that $f^{([\beta]+2)}(\xi)\sim \xi^{\beta-[\beta]-2}$ when $\xi \to 0$. Consequently, these integrals are not convergent. 
\end{itemize}
\item When $\beta$ is an even number. We get back to the expression (\ref{Iden-K-beta}) and since $\beta$ is an even number  the expressions (\ref{n-der-2}) and (\ref{n-der-3}) write down as: 
\begin{equation*}
f^{(\beta)}(\xi)= \begin{cases}\vspace{2mm} i\,c_\beta  \, \xi^{3-\beta}+  c_\alpha (-\xi)^{\alpha -\beta} -  c_\beta, \quad \xi <0, \\
i\,c_\beta \, \xi^{3-\beta}+ c_\alpha \, \xi^{\alpha-\beta} - c_{\beta}, \quad \xi >0, 
\end{cases} \quad 1 \leq \beta \leq 3,
\end{equation*}
\begin{equation*}
f^{(\beta)}(\xi)= \begin{cases}\vspace{2mm}  
c_\alpha (-\xi)^{\alpha -\beta} -  c_\beta , \quad \xi <0, \\
 c_\alpha \, \xi^{\alpha-\beta} - c_{\beta}, \quad \xi >0, 
\end{cases} \quad 4 \leq   \beta,
\end{equation*}
Therefore, we obtain 
\[ \left. \left( e^{2 \pi i \, x\xi}\, \partial^{[\beta]}_{\xi} \left( e^{-f(\xi) t}\right)\right) \right\vert^{0}_{-\infty} +  \left. \left( e^{2 \pi i \, x\xi}\, \partial^{[\beta]}_{\xi} \left( e^{-f(\xi) t}\right)\right) \right\vert^{+\infty}_{0} = 0,\]
and we can continue with the iterative process described in (\ref{Iden-K-n}) until the step  $n = [\alpha]+1$ to write
\begin{equation*}    
\begin{split}
K_{\alpha, \beta}(t,x)= &\,\frac{1}{(2\pi i x)^{[\alpha]+1}} \left((-1)^{[\alpha]+1}\,  \int_{-\infty}^{0} e^{2 \pi i \, x\xi} \, \partial^{[\alpha]+1}_{\xi} \left(e^{-f(\xi) t} \right) d \xi + \left. \left( e^{2 \pi i \, x\xi}\, \partial^{[\alpha]}_{\xi} \left( e^{-f(\xi) t}\right)\right) \right\vert^{0}_{-\infty}   \right) \\
&\, + \frac{1}{(2\pi i x)^{[\alpha]+1}} \left((-1)^{[\alpha]+1}\,  \int_{0}^{+\infty} e^{2 \pi i \, x\xi} \, \partial^{[\alpha]+1}_{\xi} \left(e^{-f(\xi) t} \right) d \xi + \left. \left( e^{2 \pi i \, x\xi}\, \partial^{[\alpha]}_{\xi} \left( e^{-f(\xi) t}\right)\right) \right\vert^{+\infty}_{0}  \right).
\end{split}  
\end{equation*} 
From this identity and by following the same arguments detailed  at  points $1.1$ and $1.2$ above (with $\alpha$ instead of $\beta$) we finally obtain  
$\ds{\vert K_{\alpha,\beta}(t,x)\vert=\frac{\vert I(t)\vert}{\vert x \vert^{[\alpha]+1}}}$, with $\vert I(t) \vert \leq C e^{\eta_1\, t}$. 
\end{enumerate}

Proposition \ref{Prop-Kernel-3} is proven. \finpv

\medskip

Summarizing, by the pointwise identities proven in Propositions \ref{Prop-Kernel-1}, \ref{Prop-Kernel-2}, and \ref{Prop-Kernel-3},  and by the parameter $n \geq 2$  defined in the expression (\ref{parameter-n}) we obtain the   unified identity (\ref{Kernel-pointwise}). Proposition \ref{Prop-Principal-Kernel} is now proven. \finpv  

\medskip 

As a corollary of this identity, we can easily estimate the kernel $K_{\alpha, \beta}$ in the $L^p-$ norms, for both the nonlocal dispersive (when $D=\mathcal{H}\partial_x$) and the local dispersive  (when $D=-\partial^{2}_x$) cases. 
\begin{Proposition}\label{Prop-Kernel-4} Let $\alpha > \beta \geq 1$ with $\alpha>2$. For all $t>0$ fixed, and for $1\leq p \leq +\infty$ the following estimate hold true:
\begin{equation*}
\Vert K_{\alpha, \beta}(t,\cdot) \Vert_{L^p} \leq C\, \frac{e^{\eta_1\, t}}{t^{\frac{1}{\alpha}}},   
\end{equation*}
where the constants $C>0$ and $\eta_1>0$ depend on $\alpha,\beta$ and $p$.
\end{Proposition}	
\pv  For $t>0$ fixed, we start by estimating the quantity $\Vert K_{\alpha,\beta}(t,\cdot)\Vert_{L^\infty}$. We recall  that $f(\xi) = i\, m(\xi)\xi + (\vert \xi \vert^\alpha - \vert \xi \vert^\beta)$ (with $m(\xi)=\vert \xi \vert$ or $m(\xi)=\vert \xi \vert^2$)  and since $\alpha>\beta$ for $M>0$ big enough  we can write
\begin{equation*}
\begin{split}
\Vert K_{\alpha,\beta}(t,\cdot) \Vert_{L^\infty}  \leq &\,  C\,\Vert e^{-f(\cdot) t} \Vert_{L^1} \leq \,  C\, \int_{\vert \xi \vert \leq M} e^{-(\vert \xi \vert^\alpha - \vert \xi \vert^\beta) t} \, d\xi +  C\, \int_{\vert \xi \vert > M} e^{-(\vert \xi \vert^\alpha - \vert \xi \vert^\beta) t} \, d\xi\\
\leq &\, C\, \int_{\vert \xi \vert \leq M} e^{\vert \xi \vert^\beta t} d\xi + C\, \int_{\vert \xi \vert >M} e^{-\vert \xi \vert^{\alpha} t} d\xi  \leq C\,  e^{\eta \, t} + \frac{C}{t^{\frac{1}{\alpha}}} \leq \frac{C( t^{\frac{1}{\alpha}} e^{\eta \, t} + 1 )}{t^{\frac{1}{\alpha}}} \leq C \,\frac{e^{\eta_1 \, t}}{t^{\frac{1}{\alpha}}}. 
\end{split}
\end{equation*} 
We thus have
\begin{equation}\label{Kernel-inf}
\Vert K_{\alpha,\beta}(t,\cdot)\Vert_{L^\infty} \leq C \,\frac{e^{\eta_1\, t}}{t^{\frac{1}{\alpha}}}. 
\end{equation}

On the other hand, we estimate now the quantity $\Vert K_{\alpha, \beta}(t,\cdot)\Vert_{L^1}$: 
\begin{equation*}
\Vert K_{\alpha,\beta}(t,\cdot) \Vert_{L^1}= \int_{\vert x \vert \leq 2} \vert K_{\alpha,\beta}(t,x)\vert \, dx + \int_{\vert x \vert > 2} \vert K_{\alpha,\beta}(t,x)\vert \, dx \leq C\, \Vert K_{\alpha,\beta}(t,\cdot) \Vert_{L^\infty} + \int_{\vert x \vert >2} \vert K_{\alpha,\beta}(t,x)\vert dx. 
\end{equation*}
The first term on the right-hand side was already estimated  in  (\ref{Kernel-inf}).  On the other hand, by  (\ref{Kernel-pointwise}) we have   
\begin{equation}\label{Estim-kernel}
\int_{\vert x \vert >2} \vert K_{\alpha,\beta}(t,x)\vert dx \leq C\, e^{\eta_0\, t} \, \int_{\vert x \vert >2} \frac{dx}{\vert x \vert^n} \leq C\, e^{\eta_1 \, t}.
\end{equation} 
We thus obtain $\Vert K_{\alpha,\beta}(t,\cdot)\Vert_{L^1} \leq \frac{C}{t^{\frac{1}{\alpha}}} +  C\,  e^{\eta_1 \, t}   \leq  C \,\frac{e^{\eta_1 \, t}}{t^{\frac{1}{\alpha}}}$.  

\medskip

Finally, the quantity $\Vert K_{\alpha,\beta}(t,\cdot)\Vert_{L^p}$, with $1<p<+\infty$, follows from  the standard interpolation inequalities. Proposition \ref{Prop-Kernel-4} is proven. \finpv

\medskip

To close this section, remark that by the identity (\ref{Kernel-pointwise}) and by the estimate $\Vert K_{\alpha,\beta}(t,\cdot)\Vert_{L^\infty} \leq C\, \frac{e^{\eta_1\, t}}{t^{\frac{1}{\alpha}}}$ proven above, for $t>0$ and $x\in \R$ we have the following pointwise estimate: 
\begin{equation}\label{Kernel-pointwise-estimates}
\vert K_{\alpha,\beta}(t,x) \vert \leq C\, \frac{e^{\eta_1\, t}}{t^{\frac{1}{\alpha}}}\, \frac{1}{1+\vert x \vert^n}, \quad \text{with} \quad \alpha>2 \ \ \mbox{and} \ \  n\geq 2.
\end{equation}
This estimate will be very useful in the following section. 
\subsection{Spatial pointwise decaying} 
\subsection*{Proof of Theorem \ref{Th-Decaying}}
Given an initial datum  $u_0 \in H^s(\R)$, with $s>\frac{5}{2}$,  we assume now that  it  verifies  $u_0 \in L^{\infty}((1+\vert \cdot \vert^{\kappa}) dx)$, with $\kappa>1$. Then, for a time $0<T<1$ we will construct a solution $u(t,x)$ of the equation (\ref{Equation-Integral}) in the following  Banach  space
\begin{equation}\label{def-E}
E_{T}= \left\{ u \in \mathcal{C}([0,T], H^s(\R)): \,\,  \Vert u \Vert_{T} <+\infty \right\},
\end{equation}
where the norm $\Vert \cdot \Vert_{T}$ depends on the parameter $\kappa$, the parameter $\alpha$,  the parameter   $n$ given in (\ref{parameter-n}), and moreover, it   also depends on the previous norm $\| \cdot \|_{s,\alpha,0}$ (defined in (\ref{Norm-F}))   as follows:
\begin{equation}\label{Def-Norm}
\Vert u \Vert_{T}= \| u \|_{s,\alpha,0}+ \sup_{0\leq t \leq T}t^{\frac{1}{\alpha}}\|(1+\vert \cdot \vert^{\min(\kappa,n)})u(t,\cdot)\|_{L^\infty}
 + \sup_{0\leq t \leq T} t^{\frac{2}{\alpha}}\Vert (1+\vert \cdot \vert^{\min(\kappa,n+1)})\partial_x u(t,\cdot)\Vert_{L^\infty}. 
\end{equation}
In this expression, the first term  norm $\| \cdot \|_{\alpha,s,0}$ will allow us to  control in the space $H^s(\R)$ each term of the nonlinear part of the equation (\ref{Equation-Integral}).   The second term  characterizes the spatially decaying properties of solutions, while the third and the fourth terms are meant to treat the (more delicate) nonlinear term $(\partial_x u)^2$. Finally,   the weights in the temporal variable  $t^{\frac{1}{\alpha}}$ and $t^{\frac{2}{\alpha}}$ are  essentially technical (due to the kernel estimates (\ref{Kernel-pointwise-estimates}))  and they will be useful to carry up all our estimates.  

\medskip

Let us start by studying the linear term in the mild formulation (\ref{Equation-Integral}). 

\begin{Proposition}\label{Prop-Decay-Lin} We have $K_{\alpha,\beta}\ast u_0 \in E_T$ and $\Vert K_{\alpha,\beta}\ast u_0 \Vert_{T} \lesssim (\Vert \vu_0 \Vert_{H^s}+\Vert (1+\vert \cdot\vert^{\kappa}) \vu_0 \Vert_{L^\infty})$. 
\end{Proposition} 
\pv We must estimate each term in the norm given in (\ref{Norm}), but recall that the first term was already considered in  (\ref{Lin-Norm-F}), and consequently, we shall focus on the second and the third term. 

\medskip 

For the second term, since  $u_0 \in  L^{\infty}((1+\vert \cdot \vert^{\kappa})dx)$ and  by the kernel estimate  (\ref{Kernel-pointwise-estimates}), for $0 \leq t\leq T$ and $x\in \R$ fixed we write  
\begin{equation*}
\begin{split}
&\, \vert K_{\alpha,\beta}(t,\cdot)\ast u_0(x) \vert  \leq   \int_{\R} \vert K_{\alpha,\beta}(t, x-y)\vert \vert u_0(y)\vert dy \leq \int_{\R} \vert K_{\alpha,\beta}(t, x-y)\vert \frac{1+\vert y \vert^{\kappa}}{1+\vert y \vert^{\kappa}} \vert u_0(y)\vert dy\\
\leq & \,  \Vert (1+ \vert \cdot \vert^{\kappa})  u_0 \Vert_{L^{\infty}} \int_{\R}  \frac{\vert K_{\alpha,\beta}(t,x-y)\vert }{1+\vert y \vert^{\kappa}} dy \lesssim \Vert (1+ \vert \cdot \vert^{\kappa})u_0 \Vert_{L^{\infty}} \, \frac{e^{c_\eta\, t}}{t^{\frac{1}{\alpha}}} \int_{\R} \frac{dy}{(1+\vert x-y\vert^{n})(1+\vert y \vert^{\kappa})} \\
\lesssim & \, \Vert (1+ \vert \cdot \vert^{\kappa})  u_0 \Vert_{L^{\infty}} \,\frac{e^{\eta_0\, t}}{ t^{\frac{1}{\alpha}}} \, \frac{1}{1+\vert x \vert^{\min(\kappa,n)}},
\end{split}  
\end{equation*}	hence we have 
\begin{equation}\label{Estim-2}
\sup_{0\leq t \leq T}t^{\frac{1}{\alpha}}\left\Vert (1+\vert \cdot \vert^{\min(\kappa,n)})K_{\alpha,\beta}(t,\cdot)\ast u_0 \right\Vert_{L^{\infty}} \lesssim  \Vert (1+\vert \cdot\vert^{\kappa}) u_0 \Vert_{L^{\infty}}.  
\end{equation}
 
For the third term,  we shall need the following technical lemma, which was essentially proven in  \cite[Lemma $4.2$]{CorJar0}:
\begin{Lemme}\label{Der-Ker-estimates}
Let $\alpha>\beta \geq 1$ with $\alpha>2$.  Let $K_{\alpha,\beta}$ the kernel given in (\ref{kernel}) with $m(\xi)=\vert \xi \vert$ or $ m(\xi)=-\vert \xi \vert^2$. Moreover, let $n \geq 2$ be the parameter defined in (\ref{parameter-n}). Then for $t>0$ we have $K_{\alpha,\beta}(t,\cdot) \in \mathcal{C}^{1}(\R)$  and the following estimates hold:
\begin{enumerate}
\item  For all $x\neq 0$, $\ds{\vert \partial_x K_{\alpha,\beta}(t,x) \vert  \leq C\, \frac{e^{\eta_1\, t}}{\vert x \vert^{n+1}}}$.
\item For all $x\in \R$, $\ds{\vert \partial_x K_{\alpha,\beta}(t,x) \vert  \leq C\, \frac{e^{\eta_1\, t}}{t^{\frac{2}{\alpha}}} \frac{1}{1+\vert x \vert^{n+1}}}$,
\end{enumerate}
for two constants $C, \eta_1>0$ which depend on $\alpha$ and $\beta$. 
\end{Lemme}

By the second point above  and by following very similar estimates done to prove (\ref{Estim-2}) we obtain 
\begin{equation}\label{Estim-3}
\sup_{0\leq t \leq T}t^{\frac{2}{\alpha}}\left\Vert (1+\vert \cdot \vert^{\min(\kappa,n+1)}) \partial_x \left(K_{\alpha,\beta}(t,\cdot)\ast u_0 \right) \right\Vert_{L^{\infty}} \lesssim  \Vert (1+\vert \cdot\vert^{\kappa}) u_0 \Vert_{L^{\infty}}. 
\end{equation}
Thus, the wished estimate follows from (\ref{Lin-Norm-F}), (\ref{Estim-2}),  and (\ref{Estim-3}).  Proposition \ref{Prop-Decay-Lin} is now proven. \finpv

\medskip

We study now the nonlinear term in the mild formulation (\ref{Equation-Integral}). For the sake of simplicity, we shall only consider the case $\gamma_1=\gamma_2=\gamma_3=1$ with $\alpha>7/2$. The other case: $\gamma_1=1$, $\gamma_2=\gamma_3=0$  with $\alpha>2$, essentially follows the same estimates with the obvious minor modifications.  

\begin{Proposition}\label{Prop-Decay-Nonlin} Lat   $\alpha>\frac{7}{2}$,   $s>\frac{5}{2}$  and let $\eta>0$ be quantity given in (\ref{eta}). Define   $0<\eta_2 < \min(\eta, 1-\frac{3}{\alpha})$. Then the following estimate holds:
\begin{equation*}
\left\Vert \int_{0}^{t} K_{\alpha,\beta}(t-\tau,\cdot)\ast \left( \partial_x(u^2)+\partial^{2}_{x}(u^2)+(\partial_x u )^2 \right)(\tau, \cdot) d \tau  \right\Vert_{T} \leq C\, T^{\eta_2}\, \Vert u \Vert^{2}_{T}. \end{equation*}
\end{Proposition}
\pv We get back to the definition of the norm $\Vert \cdot \Vert_T$ given in (\ref{Def-Norm}), where we must estimate each term in this expression.  We recall that the first term $\| \cdot \|_{s,\alpha,0}$ was already estimated in (\ref{NonLin-Norm-F}); and for the quantity $\eta>0$ given in (\ref{eta}) we have 
\begin{equation}\label{Estim-NonLin-1}
\left\Vert \int_{0}^{t}K_{\alpha,\beta}(t-\tau,\cdot)\ast\Big(  \partial_x (u^2)  +  \partial^{2}_{x}(u^2) + (\partial_x u)^2   \Big)(\tau, \cdot) \, d\tau \right\Vert_{s,\alpha,0} \lesssim T^\eta \,  \Vert u \Vert^{2}_{T}.    
\end{equation}

For the second term, the following estimate holds:
\begin{equation}\label{Estim-NonLin-2}
\sup_{0\leq t \leq T} t^{\frac{1}{\alpha}} \left\| (1+|\cdot|^{\min(\kappa,n)}) \int_{0}^{t}K_{\alpha,\beta}(t-\tau,\cdot)\ast (\partial_x(u^2)+\partial^{2}_{x}(u^2)+(\partial_x u)^2)(\tau,\cdot) d \tau  \right\|_{L^\infty}  \lesssim T^{1-\frac{2}{\alpha}}\, \| u \|^{2}_{T}.    
\end{equation}
Indeed, to estimate the expression $\partial_x(u^2)$, 
for $t>0$ and $x \in \R$ fixed,  by the kernel estimate (\ref{Kernel-pointwise-estimates}) and by the first and the  second  expressions in (\ref{Def-Norm}), and by recalling that $s>\frac{5}{2}$ and we have the continuous embedding $\| \partial_x u \|_{L^\infty} \lesssim \| u \|_{H^s}$, we write
\begin{equation*}
\begin{split}
&\, \left|\int_{0}^{t} K_{\alpha,\beta}(t-\tau,\cdot)\ast \partial_{x}(u^2)(\tau,x) \right| \\
\leq &\, \int_{0}^{t} \int_{\R} |K_{\alpha,\beta}(t-\tau, x-y)| | u(t,y)| | \partial_y u(\tau,y)| dy\, d\tau  \\
\lesssim & \| u \|_{T}\,e^{\eta\, t} \int_{0}^{t} \frac{1}{(t-\tau)^{\frac{1}{\alpha}}}\, \frac{1}{\tau^{\frac{1}{\alpha}}} \left( \int_{\R} \frac{1}{1+| x-y |^{n}}\frac{1}{1+|y|^{\min(\kappa,n)}}\,  dy  \right) \| \partial_y u (\tau,\cdot) \|_{L^\infty}  d \tau \\
\lesssim & \| u \|^{2}_{T}\,e^{\eta\, t} \int_{0}^{t} \frac{1}{(t-\tau)^{\frac{1}{\alpha}}}\, \frac{1}{\tau^{\frac{1}{\alpha}}} \left( \int_{\R} \frac{1}{1+| x-y |^{n}}\frac{1}{1+|y|^{\min(\kappa,n)}}\, dy  \right) d \tau \\
\lesssim &\| u \|^{2}_{T}\, e^{\eta\, t} t^{-\frac{2}{\alpha}+1}\,\frac{1}{1+|x|^{\min(\kappa,x)}},
\end{split}
\end{equation*}
which yields the estimate $\ds{ \sup_{0\leq t \leq T} t^{\frac{1}{\alpha}} \left\Vert (1+\vert \cdot\vert^{\min(k,n)}) \int_{0}^{t}K_{\alpha,\beta}(t-\tau,\cdot)\ast \partial_x(u^2)(\tau,\cdot) d \tau  \right\Vert_{L^\infty} \lesssim T^{1-\frac{1}{\alpha}} \, \| u \|^{2}_{T}}$.

To estimate the expression $\partial^{2}_{x}(u^2)$ we remark that we can write $K_{\alpha,\beta}\ast (\partial^{2}_{x}(u^2))=\partial_x K_{\alpha,\beta}\ast  2\left( u \partial_x u \right)$. By the second point of Lemma \ref{Der-Ker-estimates},  we obtain 
\begin{equation*}
\begin{split}
&\, \left|\int_{0}^{t} K_{\alpha,\beta}(t-\tau,\cdot)\ast \partial^{2}_{x}(u^2)(\tau,x) \right| \\
\lesssim &\, \int_{0}^{t} \int_{\R} |\partial_{x} K_{\alpha,\beta}(t-\tau, x-y)| | u(t,y)| | \partial_y u(\tau,y)| dy\, d\tau  \\
\lesssim & \| u \|_{T}\,e^{\eta\, t} \int_{0}^{t} \frac{1}{(t-\tau)^{\frac{2}{\alpha}}}\, \frac{1}{\tau^{\frac{1}{\alpha}}} \left( \int_{\R} \frac{1}{1+| x-y |^{n+1}}\frac{1}{1+|y|^{\min(\kappa,n)}}\,  dy  \right) \| \partial_x u(\tau, \cdot) \|_{L^\infty} d \tau \\
\lesssim &\| u \|^{2}_{T}\, e^{\eta_1\, t} t^{-\frac{3}{\alpha}+1}\,\frac{1}{1+|x|^{\min(\kappa,x)}}.
\end{split}
\end{equation*}
We thus have  $\ds{ \sup_{0\leq t \leq T} t^{\frac{1}{\alpha}} \left\Vert (1+\vert \cdot\vert^{\min(k,n)}) \int_{0}^{t}K_{\alpha,\beta}(t-\tau,\cdot)\ast \partial^{2}_{x}(u^2)(\tau,\cdot) d \tau  \right\Vert_{L^\infty} \lesssim T^{1-\frac{2}{\alpha}} \, \| u \|^{2}_{T}}$.

\medskip

Similarly, for the expression $K_{\alpha,\beta}\ast (\partial_x u)^2$ we just write $K_{\alpha,\beta}\ast (\partial_x u)^2= K_{\alpha,\beta}\ast ((\partial_x u) (\partial_x u))$ and we have the estimate $\ds{ \sup_{0\leq t \leq T} t^{\frac{1}{\alpha}} \left\Vert (1+\vert \cdot\vert^{\min(k,n)}) \int_{0}^{t}K_{\alpha,\beta}(t-\tau,\cdot)\ast (\partial_x u)^2(\tau,\cdot) d \tau  \right\Vert_{L^\infty} \lesssim T^{1-\frac{2}{\alpha}} \, \| u \|^{2}_{T}}$.

\medskip

Finally, we recall that since  $0<T<1$ we have $T^{1-\frac{1}{\alpha}}<T^{1-\frac{2}{\alpha}}$; and we thus obtain  the wished estimate stated in  (\ref{Estim-NonLin-2}). 

\medskip

For the third term in the norm $\| \cdot \|_T$ (given in (\ref{Def-Norm})) we essentially follow the same arguments exposed above to obtain the estimate 
\begin{equation}\label{Estim-NonLin-3}
\begin{split}
&\,\sup_{0\leq t \leq T} t^{\frac{2}{\alpha}} \left\| (1+|\cdot|^{\min(\kappa,n+1)})  \int_{0}^{t} \partial_x  K_{\alpha,\beta}(t-\tau,\cdot)\ast (\partial_x(u^2)+\partial^{2}_{x}(u^2)+(\partial_x u)^2)(\tau,\cdot) d \tau \right\|_{L^\infty} \\
\lesssim &\,  T^{1-\frac{3}{\alpha}}\, \| u \|^{2}_{T}.    
\end{split}
\end{equation}
Indeed, for the reader's convenience, we shall only mention that to treat the expression $\partial^{2}_{x}(u^2)$ we write 
$\partial_x K_{\alpha,\beta}\ast \partial^{2}_{x}(u^2)= \partial_x K_{\alpha,\beta}\ast (\partial_x u )^2 +  \partial_x K_{\alpha,\beta}\ast (u \, \partial^{2}_{x}u) $. Here, to control the last term $\partial^{2}_{x}u$ we use the continuous embedding $\| \partial^{2}_{x}u \|_{L^\infty} \lesssim \| u \|_{H^s}$, which is valid for  $s>\frac{5}{2}$.

\medskip

To finish the proof, we set $\eta_2=\min(\eta, 1-\frac{3}{\alpha})$ from which we get the  desired estimate. Proposition \ref{Prop-Decay-Nonlin} is proven. \finpv  

\medskip

With Propositions \ref{Prop-Decay-Lin} and \ref{Prop-Decay-Nonlin} at our disposal, there exists a solution $u \in E_{T_0}$ to the equation (\ref{Equation-Integral}), for a time $0<T_0< T\leq 1$ small enough. But, by the embedding  $E_{T_0} \subset \mathcal{C}([0,T_0],H^s(\R))$ and since the equation (\ref{Equation}) is locally well-posed in this space (in particular we have the uniqueness of the solution) this solution is the same to the one constructed in Theorem \ref{Th1}.  

\medskip

Until now we have proven  the estimate (\ref{Decaying-Solution}) for all time $0<t\leq T_0$. Thereafter,  by following the same arguments of  \cite[Theorem $4.2$]{CorJar0}  this estimate is extended to the time $T$. Theorem \ref{Th-Decaying} is now proven. \finpv  

\subsection{Asymptotic profiles and optimality}
Once the problem of  the spatial pointwise decaying of our equation \eqref{Equation} is finished, our next objective is to show in which cases we could speak of optimally  of this decaying. To do this, we will start by giving an asymptotic profile of the solution in the spatial variable. 

\subsection*{Proof of Theorem \ref{Th-Asymptotic-Profile}}
Since the solution $u(t,x)$ writes down as in the integral formulation (\ref{Equation-Integral}), we start by proving that the first term on the right-hand side in has  the following asymptotic development: 
\begin{equation}\label{eq18}
K_{\alpha,\beta}(t,\cdot) \ast u_0 (x)= K_{\alpha,\beta}(t,x)\left(  \int_{\R} u_0 (y)dy\right)+ R_{1}(t,x),  \quad \vert x \vert \to +\infty,   
\end{equation} with $\ds{| R_1(t,x)|=o(t)\left( 1/|x|^n \right)}$. Indeed,  for $t>0$ and $x \in \R$ fix  this term can be decomposed as follows:
\begin{equation}\label{Split-Linear}
\begin{split}
\int_{\R} K_{\alpha,\beta}(t,x-y)u_0(y) dy =&\, K_{\alpha,\beta}(t,x)\left(  \int_{\R} u_0 (y)dy\right)  + \int_{\vert y \vert < \frac{\vert x \vert}{2}} (K_{\alpha,\beta} (t, x-y)- K_{\alpha,\beta} (t, x)) u_0(y) dy  \\
&\, + \int_{\vert y \vert > \frac{\vert x \vert}{2}} K_{\alpha,\beta}(t, x-y) u_0(y) dy - K_{\alpha,\beta}(t,x) \left( \int_{\vert y \vert > \frac{\vert x \vert}{2}} u_0 (y) dy\right)\\
=&\,   K_{\alpha,\beta}(t,x)\left(  \int_{\R} u_0 (y)dy\right)  + I_1 + I_2 + I_3,
\end{split}
\end{equation}
hence, we define $\ds{R_1=  I_1 + I_2 + I_3 }$ and we will verify that the following statement holds:
\begin{equation}\label{eq14}
\vert R_1 \vert \leq \frac{c(u_0,t)}{\vert x \vert^{n+\varepsilon}}, \quad \vert x \vert \to +\infty, \qquad \varepsilon>0.  
\end{equation}
To estimate term $I_1$ we need  Lemma \ref{Der-Ker-estimates}.
Since  $K_{\alpha,\beta}(t,\cdot) \in \mathcal{C}^{1}(\R)$, we write $K_{\alpha,\beta}(t, x-y) - K_{\alpha,\beta} (t,x) = - y  \ \partial_x K_{\alpha,\beta} (t, x- \theta y )$,  for some  $0< \theta < 1$.  Then, by this identity and using the first  point of Lemma \ref{Der-Ker-estimates} we get
\begin{eqnarray*}
 I_1 & \leq & 	\int_{\vert y \vert < \frac{\vert x \vert}{2}} \vert  (K_{\alpha,\beta} (t, x-y)- K_{\alpha,\beta} (t, x)) \vert \vert  u_0(y) \vert  dy \leq \int_{\vert y \vert < \frac{\vert x \vert}{2}} \vert y \vert \vert \partial_x K_{\alpha,\beta} (t, x- \theta y ) \vert \vert u_0(y) \vert d y \\
 & \lesssim  &  e^{c_{\eta_1}  t } \int_{\vert y \vert < \frac{\vert x \vert}{2}} \frac{\vert y \vert \vert u_0 (y) \vert}{\vert x - \theta y \vert^{n+1}} dy.  
 \end{eqnarray*}
We study now the expression $\ds{\frac{1}{\vert x -\theta y \vert^{n +1}}}$. As we have $0<\theta < 1$, and moreover, as we have $\vert y \vert < \frac{\vert x \vert }{2}$, then we can write $\vert x -\theta y \vert \geq \vert x \vert  - \theta \vert y \vert \geq \vert x \vert - \vert y \vert \geq \frac{\vert x \vert}{2}$; and thus  we get $\ds{\frac{1}{\vert x -\theta y \vert^{n+1}} \lesssim  \  \frac{1}{\vert x \vert^{n+1}}}$. With this inequality and recalling that the initial datum verifies
 $\ds{\vert u_0 (y) \vert \leq \frac{c}{1+\vert y \vert^{\kappa }}}$ (with $ \kappa > n$), we can write
  $$  e^{c_{\eta_1} \,  t } \int_{\vert y \vert < \frac{\vert x \vert}{2}} \frac{\vert y \vert \vert u_0 (y) \vert}{\vert x - \theta y \vert^{n +1}} dy \lesssim  \frac{ e^{c_{\eta_1} \,  t }}{\vert x \vert^{n + 1}} \int_{\vert y \vert < \frac{\vert x \vert}{2}} \frac{\vert y \vert }{1+ \vert y \vert^{\kappa}} dy  \lesssim   \frac{ e^{c_{\eta_1} \, t }}{\vert x \vert^{n+1}},$$
hence   we have
 \begin{equation}\label{eq15}
 I_1 \lesssim  \frac{ e^{c_{\eta_1} \, t }}{\vert x \vert^{n+1}}, \quad \vert x \vert \to +\infty.   
\end{equation}
For  the term  $I_2$, as $\ds{\vert u_0 (y) \vert \leq \frac{c_0}{\vert y \vert^{\kappa}}}$ (for $\vert y \vert$ large enough) and  as we have $\vert y \vert > \frac{\vert x \vert}{2}$, then  we write 
\begin{eqnarray*}
 I_2 &\leq & \int_{\vert y \vert > \frac{\vert x \vert}{2}} \vert K_{\alpha,\beta}(t,x-y) \vert \vert u_0 (y) \vert d y \lesssim   \int_{\vert y \vert > \frac{\vert x \vert}{2}} \frac{\vert K_{\alpha,\beta}(t, x-y) \vert}{\vert y \vert^{\kappa}} dy \lesssim \frac{1}{\vert x \vert^{\kappa}} \int_{\vert y \vert > \frac{\vert x \vert}{2}} \vert K_{\alpha,\beta}(t, x-y) \vert \\
 &\lesssim & \frac{1}{\vert x \vert^{\kappa}} \Vert K_{\alpha,\beta} (t,\cdot) \Vert_{L^1},
\end{eqnarray*}
but,  by Proposition \ref{Prop-Kernel-4} we have $\Vert K_{\alpha,\beta} (t,\cdot) \Vert_{L^1} \lesssim   \frac{e^{c_{\eta_1} \, t }}{ t^{\frac{1}{\alpha}}}$, and for $\kappa>n$ we get 
\begin{equation}\label{eq16}
I_2 \lesssim   \frac{e^{c_{\eta_1} \,t } }{t^{\frac{1}{\alpha}}} \frac{1}{\vert x \vert^{\kappa}}, \quad \vert x \vert \to +\infty. \end{equation}
Finally, in order to study  the term $I_3$,  recall first  that   by the estimate  (\ref{Kernel-pointwise-estimates})   for $\vert x \vert $ enough we have $\ds{\vert K_{\alpha,\beta} (t,x) \vert \lesssim  \frac{e^{c_{\eta_{1}} \, t }}{t^{\frac{1}{\alpha}}} \frac{1}{\vert x \vert^{n}}}$. Moreover, recall that the initial datum verifies $\ds{\vert u_0 (y) \vert \lesssim  \frac{1}{1+ \vert y \vert^{\kappa}}}$ (with $\kappa=n+ \varepsilon$). Then we write 
\begin{eqnarray}\label{eq17}\nonumber
I_3 &\lesssim &  \frac{e^{c_{\eta_1} \, t }}{ t^{\frac{1}{\alpha}}} \frac{1}{\vert x \vert^{n}} \int_{\vert y \vert > \frac{\vert x \vert}{2}} \vert u_0(y) \vert dy  \lesssim  \frac{e^{c_{\eta_1} \, t }}{ t^{\frac{1}{\alpha}}} \frac{1}{\vert x \vert^{n}} \int_{\vert y \vert > \frac{\vert x \vert}{2}} \frac{1}{1+\vert y \vert^{n+\varepsilon}} dy \\
& \lesssim &  \frac{e^{c_{\eta_1} \, t }}{t^{\frac{1}{\alpha}}} \frac{1}{\vert x \vert^{n+\varepsilon}} \int_{\vert y \vert > \frac{\vert x \vert}{2}} \frac{1}{1+\vert y \vert^{n}} dy \lesssim   \frac{e^{c_{\eta_1} \, t }}{ t^{\frac{1}{\alpha}}} \frac{1}{\vert x \vert^{n+\varepsilon}} \int_{\R} \frac{1}{1+\vert y \vert^{n}} dy \lesssim  \frac{e^{c_{\eta_1} \, t }}{ t^{\frac{1}{\alpha}}} \frac{1}{\vert x \vert^{n + \varepsilon}}. 
\end{eqnarray}
Thus, the desired estimate (\ref{eq14}) follows from (\ref{eq15}), (\ref{eq16}), and (\ref{eq17}); and we have the wished profile given in (\ref{eq18}). 

\medskip

Now, we focus on the nonlinear term on the right-hand side of the equation (\ref{Equation-Integral}).  The first and the second nonlinear terms: $\partial_x(u^2)$ and $\partial^{2}_{x}(u^2)$, can be estimated as follows.  For $t>0$ and $x \in \R$ fix we write  
\begin{eqnarray*}
\gamma_1 \int_{0}^{t} K_{\alpha,\beta} (t- \tau, \cdot) \ast \partial_x (u^2)(\cdot,\tau) d\tau  =   \gamma_1 \int_{0}^{t} \int_{\R} \partial_y K_{\alpha,\beta} (t- \tau, x-y)   \, u^2(\tau,y) dy \, d\tau=(a). 
\end{eqnarray*}
Then,  by the second point  of Lemma \ref{Der-Ker-estimates}, and since   by  (\ref{Decaying-Solution}) we have  $\ds{\vert u(\tau, y) \vert^2 \lesssim  \frac{1}{\tau^ {\frac{2}{\alpha}} \, (1+\vert y \vert^{2 n})}}$, then we get 
\begin{eqnarray*}
(a)&\lesssim &  \gamma_1  \int_{0}^{t} \frac{e^{c_{\eta_1} (t-\tau) }}{(t-\tau)^{\frac{2}{\alpha}} \, \tau^{\frac{2}{\alpha}}} \int_{\R} \frac{1}{1+ \vert x-y \vert^{n+ 1}} \frac{1}{1+\vert y \vert^{2n}}    d y \, d \tau \\
& \lesssim &  \gamma_1   e^{c_{\eta_1} t}  \int_{0}^{t} \frac{1}{(t-\tau)^{\frac{2}{\alpha}} \, \tau^{\frac{2}{\alpha}}}   d \tau  \left( \int_{\R} \frac{1}{1+ \vert x-y \vert^{n + 1}} \frac{1}{1+\vert y \vert^{2n }}   d y  \right) \\
&\lesssim & \gamma_1  e^{c_{\eta_1} t} \left(\int_{0}^{t} \frac{1}{(t-\tau)^{\frac{2}{\alpha}} \, \tau^{\frac{2}{\alpha}}} \right) \, \frac{1}{1+\vert x \vert^{n +1}}. 
\end{eqnarray*}
As $\alpha>2$, this integral computes down as $\ds{  \int_{0}^{t} \frac{ d \tau  }{(t-\tau)^{\frac{2}{\alpha}} \, \tau^{\frac{2}{\alpha}}}  \lesssim \frac{1}{t^{\frac{4}{\alpha} -1}}}$. 
\\ 

On the other hand, we write  
\begin{eqnarray*}
\gamma_2 \int_{0}^{t} K_{\alpha,\beta} (t- \tau, \cdot) \ast \partial_{x}^{2} (u^2)(\cdot,\tau) d\tau & =&  \gamma_2 \int_{0}^{t} \int_{\R} \partial_{y}^2  K_{\alpha,\beta} (t- \tau, x-y) \,   u^2(\tau,y) dy d\tau \\
&\lesssim &   \int_{0}^{t} \frac{e^{c_{\eta_1} (t-\tau) }}{(t-\tau)^{\frac{3}{\alpha}} \, \tau^{\frac{2}{\alpha}}} \int_{\R} \frac{1}{1+ \vert x-y \vert^{n+ 2}} \frac{1}{1+\vert y \vert^{2n}}    d y \, d \tau \\
& \lesssim &    e^{c_{\eta_1} t}  \int_{0}^{t} \frac{1}{(t-\tau)^{\frac{3}{\alpha}} \, \tau^{\frac{2}{\alpha}}}   d \tau  \left( \int_{\R} \frac{1}{1+ \vert x-y \vert^{n + 2}} \frac{1}{1+\vert y \vert^{2n}}   d y  \right) \\
&\lesssim &   e^{c_{\eta_1} t} \left(\int_{0}^{t} \frac{1}{(t-\tau)^{\frac{3}{\alpha}} \, \tau^{\frac{2}{\alpha}}} \right) \, \frac{1}{1+\vert x \vert^{n +2}}.
\end{eqnarray*}
As $\alpha>\frac{7}{2}$, then this integral computes down as $\ds{  \int_{0}^{t} \frac{ d \tau  }{(t-\tau)^{3 /\alpha} \, \tau^{2 / \alpha}}  \leq \frac{c}{t^{5 / \alpha -1}}}$.

Finally, we must study the third nonlinear term $(\partial_x u)^2$, which must be treated differently from the previous ones. 
\begin{Remarque} When studying this nonlinear term in the same fashion  as the previous ones we obtain the integral $\ds{\int_{0}^{t} \frac{d \tau }{(t-\tau)^{\frac{1}{\alpha}} \tau^{\frac{4}{\alpha}}}}$, which converges as long as $\alpha > 4$. But this constraint excludes the physically relevant value $\alpha=4$. 
\end{Remarque}
\begin{Remarque}
The more precise analysis on the term $(\partial_x u)^2$ which we shall perform will allow us to prove an interesting optimal criterion of the 
pointwise decaying  of  solutions. 
\end{Remarque} 
We shall prove the following identity:
\begin{equation}\label{gamma3eq18}
\int_{0}^{t} K_{\alpha,\beta}(t - \tau ,\cdot) \ast \gamma_3 (\partial_x u)^{2} (\tau,\cdot) \,d\tau =  \gamma_3\, K_{\alpha,\beta}(t,x) \int_{0}^{t} \| u(\tau,\cdot) \|^{2}_{\dot{H}^1} d \tau  +  R_{2}(t,x),  \quad \vert x \vert \to +\infty,   
\end{equation} with $\ds{\vert R_2(t,x) \vert = o(t)\left( \frac{1}{| x |^n} \right)}$. Indeed,   we follow the same ideas of the identity (\ref{Split-Linear}) to  write
\begin{equation*}
\gamma_3 \int_{0}^{t} K_{\alpha,\beta}(t - \tau ,x) \ast (\partial_x u)^{2} (\tau,x) \,d\tau = \gamma_3 \int_{0}^{t} K_{\alpha,\beta}(t - \tau ,x) \| u(\tau,\cdot)\|^{2}_{\dot{H}^1} \,d\tau + J_1(t) + J_2 (t) + J_3 (t), 
\end{equation*}
where 
\begin{equation*}\label{J1}
J_1(t)= - \gamma_3 \int_{0}^{t} \int_{\vert y \vert < \frac{\vert x \vert}{2}} (K_{\alpha,\beta} (t- \tau, x-y)- K_{\alpha,\beta} (t-\tau, x)) (\partial_y u)^{2}(y,\tau)  dy \  d\tau,
\end{equation*}
\begin{equation*}\label{J2}
J_2(t)= -  \gamma_3 \int_{0}^{t}  \int_{\vert y \vert > \frac{\vert x \vert}{2}} K_{\alpha,\beta}(t-\tau, x-y)  ((\partial_y u)^{2}) (\tau,y) dy \ d\tau,
\end{equation*}
and 
\begin{equation*}\label{J3}
J_3(t)=  \gamma_3 \int_{0}^{t} K_{\alpha,\beta}(t-\tau,x) \left( \int_{\vert y \vert > \frac{\vert x \vert}{2}}(\partial_y u)^{2}(\tau,y) dy\right) \, d\tau.   
\end{equation*}
Moreover, we write 
\begin{equation*}
\begin{split}
&\gamma_3 \int_{0}^{t} K_{\alpha,\beta}(t - \tau ,\cdot) \ast (\partial_x u)^{2} (\tau,x) \,d\tau\\
= &\, \gamma_3\, K_{\alpha,\beta}(t,x) \int_{0}^{t}\| u(\tau,\cdot)\|^{2}_{\dot{H}^1} dt  -  \gamma_3\, K_{\alpha,\beta}(t,x) \int_{0}^{t}\| u(\tau,\cdot)\|^{2}_{\dot{H}^1} dt+  \int_{0}^{t} K_{\alpha,\beta}(t - \tau ,x) \| u(\tau,\cdot)\|^{2}_{\dot{H}^1} \,d\tau  \\
&\,+ J_1(t) + J_2 (t) + J_3 (t), \\
=&\, J_1(t) + J_2 (t) + J_3 (t)  - \gamma_3\,  \int_{0}^{t} (K_{\alpha,\beta}(t,x)-K_{\alpha,\beta}(t-\tau,x)\| u(\tau,\cdot)\|^{2}_{\dot{H}^1} d \tau\\
=&\, J_1(t) + J_2 (t) + J_3 (t) + J_4(t). 
\end{split}
\end{equation*}

As before, we  we define $\ds{ R_2 =  J_1(t) + J_2(t) + J_3(t)+J_4(t)}$ and we will verify that the following statement holds:
\begin{equation}\label{eq14-bis}
\vert R_2 \vert \leq \frac{c(u,t)}{\vert x \vert^{n+\varepsilon}}, \quad \varepsilon>0, \quad \vert x \vert \to +\infty. 
\end{equation} 

By an analysis similar to the one done above, for the term  $J_1(t)$  we have 
 \begin{eqnarray*}
J_1 (t) &\lesssim&
\gamma_3 \int_{0}^{t} \int_{\vert y \vert < \frac{\vert x \vert}{2}} \vert  (K_{\alpha,\beta} (t-\tau, x-y)- K_{\alpha,\beta} (t-\tau, x)) \vert \vert  (\partial_y u)^{2} (t,y) \vert  dy \  d\tau \\
&\lesssim&  \gamma_3 \int_{0}^{t} \int_{\vert y \vert < \frac{\vert x \vert}{2}} \vert y \vert \vert \partial_x K_{\alpha,\beta} (t-\tau, x- \theta y ) \vert \vert (\partial_y u)^{2} (\tau,y) \vert d y  \  d\tau\\
 &\lesssim&   \gamma_3  \int_{0}^{t} \frac{e^{c_{\eta_1} \, (t-\tau)}}{(t-\tau)^{\frac{2}{\alpha}}} \int_{\vert y \vert < \frac{\vert x \vert}{2}} \frac{\vert y \vert \vert (\partial_y u)^{2} (\tau,y) \vert}{\vert x - \theta y \vert^{n+1}} dy \ d\tau \\ 
&\lesssim&  \frac{\gamma_3}{\vert x \vert^{n + 1}}  \int_{0}^{t} \frac{e^{c_{\eta_1} \, (t-\tau)}}{(t-\tau)^{\frac{2}{\alpha}}}
 \int_{\vert y \vert < \frac{\vert x \vert}{2}} \frac{\vert y \vert }{1+ \vert y \vert^{n+1}} \vert (\partial_y u)^{2} (\tau,y) \vert d y  \  d\tau \\
 &\lesssim&  \frac{\gamma_3}{\vert x \vert^{n + 1}}  \int_{0}^{t} \frac{e^{c_{\eta_1} \, (t-\tau)}}{(t-\tau)^{\frac{2}{\alpha}}} \|u(\tau,\cdot) \|_{\dot{H}^1}^{2} d\tau.
 \end{eqnarray*}
Hence we obtain 
\begin{equation}\label{eq15-bis}
 J_1(t) \leq  \frac{C(t,u) }{\vert x \vert^{n+1}}, \quad \vert x \vert \to +\infty.   
\end{equation}

For the term $J_2(t)$, recall the estimate   $\ds{\vert ((\partial_y u)) (\tau,y) \vert \leq \frac{c}{\tau^{\frac{2}{\alpha}}} \frac{1}{\vert y \vert^{n+1}}}$ (for $\vert y \vert$ large enough) and moreover, since $\vert y \vert > \frac{\vert x \vert}{2}$, then  we write 
\begin{equation*}
\begin{split}
J_2(t) \lesssim &\,  \gamma_3 \int_{0}^{t} \int_{\vert y \vert > \frac{\vert x \vert}{2}} \vert K_{\alpha,\beta}(t-\tau,x-y) \vert ((\partial_x u)^{2}) (\tau,y) \vert  \vert d y d \tau \\
\lesssim &\,  \gamma_3 \int_{0}^{t} \int_{\vert y \vert > \frac{\vert x \vert}{2}} \frac{\vert K_{\alpha,\beta}(t-\tau, x-y) \ \vert (\partial_x u) (\tau,y)}{\vert y \vert^{n+1}} dy \\
\lesssim &\,  \frac{\gamma_3}{\vert x \vert^{n+1}}  \int_{0}^{t}
 \frac{1}{\tau^{\frac{2}{\alpha}}} \Vert K_{\alpha,\beta} (t-\tau,\cdot) \Vert_{L^2} \ \|u(\tau)\|_{\dot{H}^1} \,d\tau \\
 \lesssim &\,  \frac{\gamma_3}{\vert x \vert^{n+1}}  \int_{0}^{t}
 \frac{1}{\tau^{\frac{2}{\alpha}}}  \frac{1}{(t-\tau)^{\frac{1}{\alpha}}} \ \|u(\tau)\|_{\dot{H}^1} \,d\tau.
\end{split}
\end{equation*}
We thus have 
\begin{equation}\label{eq16-bis}
 J_2(t) \lesssim \frac{C(t,u) }{\vert x \vert^{n+1}}, \quad \vert x \vert \to +\infty.   
\end{equation}
In order to study  the term $J_3(t)$,  recall first  that always  by  (\ref{Kernel-pointwise-estimates})  for $\vert x \vert $ enough enough we have the estimate $\ds{\vert K_{\alpha,\beta} (t,x) \vert \leq  C  \frac{e^{c_{\eta_1} \, t }}{t^{\frac{1}{\alpha}}} \frac{1}{\vert x \vert^{n}}}$. Then we write
\begin{eqnarray*}
J_3 (t) &\lesssim &  \frac{\gamma_3}{\vert x \vert^{n}} 
\int_{0}^{t}  \frac{e^{c_{n_1} \, (t- \tau) }}{ \tau^{\frac{2}{\alpha}}} \Big( \int_{\vert y \vert > \frac{\vert x \vert}{2}} 
 \frac{1}{1 +\vert y \vert^{n+1}}\vert  (\partial_x u) (\tau,y) \vert dy \Big) d\tau \\
 &\lesssim& \frac{\gamma_3}{\vert x \vert^{n+1}} 
\int_{0}^{t}  \frac{e^{c_{\eta_1} \, (t- \tau) }}{ \tau^{\frac{2}{\alpha}}} \Big( \int_{\vert y \vert > \frac{\vert x \vert}{2}} \frac{1}{1 +\vert y \vert^{n}} \vert ((\partial_x u)) (\tau,y) \vert dy \Big) d\tau \\
& \lesssim & \frac{\gamma_3}{\vert x \vert^{n+1}} 
\int_{0}^{t}  \frac{e^{c_{\eta_1} \, (t- \tau) }}{ \tau^{\frac{2}{\alpha}}} \|u(\tau,\cdot)\|_{\dot{H}^{1}} d\tau. 
\end{eqnarray*}
Hence, 
\begin{equation}\label{eq16-bis-bis} 
 J_3(t) \lesssim \frac{C(t,u) }{\vert x \vert^{n+1}}, \quad \vert x \vert \to +\infty.   
\end{equation} 
Finally, we must estimate the term $J_4(t)$. We use the mean value theorem (in the time variable) and for a time  $t-\tau \leq \tau_1 \leq t$ we write 
\begin{equation}\label{J4}
\begin{split}
J_4(t) \leq &\,  \gamma_3 \int_{0}^{t} | K_{\alpha,\beta}(t,x)-K_{\alpha,\beta}(t-\tau,x) | \, \| u(\tau,\cdot)\|^{2}_{\dot{H}^1} d \tau \\
\leq &\, \gamma_3 \int_{0}^{t} | \partial_t K_{\alpha,\beta}(\tau_1, x) |\, | \tau | \,  \| u(\tau,\cdot)\|^{2}_{\dot{H}^1} d \tau.
\end{split}
\end{equation}
At this point, we need to estimate the expression $\ds{| \partial_t K_{\alpha,\beta}(\tau_1, x) |}$:
\begin{Lemme} Let $\alpha>\beta \geq 1$ with $\alpha>2$.  Let $K_{\alpha,\beta}$ the kernel given in (\ref{kernel}) with $m(\xi)=\vert \xi \vert$ or $ m(\xi)=-\vert \xi \vert^2$. Moreover, let $n \geq 2$ be the parameter defined in (\ref{parameter-n}). For $t>0$ the following estimate hold:
\begin{equation}\label{Time-Der-Kernel}
| \partial_t K_{\alpha,\beta}(t,x) | \leq C \frac{e^{\eta_1\, t}}{| x |^{n+1}}, \quad |x|\to +\infty,  \end{equation}
for two constants $C,\eta_1>0$ depending on $\alpha$ and $\beta$. 
\end{Lemme}
\pv Recall that the kernel $K_{\alpha,\beta}(t,x)$ solves the equation (\ref{Equation-Kernel}) and for $t>0$ fixed we can write 
\begin{equation*}
| \partial_t K_{\alpha,\beta}(t,x) | \leq | D(\partial_x K_{\alpha,\beta}(t,x))| + |D^{\alpha}_{x} K_{\alpha,\beta}(t,x)|+ | D^{\beta}_{x} K_{\alpha,\beta}(t,x)|.  
\end{equation*} 
Each term on the right-hand side is  essentially a derivative of the Kernel $K_{\alpha,\beta}(t,x)$ in the spatial variable. Consequently, by following the same ideas in the proof of \cite[Lemma $4.2$]{CorJar0}, for $x\neq 0$ we get the following estimates
\begin{equation*}
| D(\partial_x K_{\alpha,\beta}(t,x))| \leq \begin{cases}\vspace{2mm} 
C\frac{e^{\eta_1\, t}}{|x|^{n+2}}, \quad m(\xi)=|\xi|,\\
C\frac{e^{\eta_1\, t}}{|x|^{n+3}}, \quad m(\xi)=-|\xi|^2, 
\end{cases}
\end{equation*}
\begin{equation*}
 |D^{\alpha}_{x} K_{\alpha,\beta}(t,x)| \leq C \frac{e^{\eta_1 \, t }}{|x|^{n+[\alpha]}},   
\end{equation*}
and 
\begin{equation*}
 |D^{\alpha}_{x} K_{\alpha,\beta}(t,x)| \leq C \frac{e^{\eta_1 \, t }}{|x|^{n+[\beta]}},   
\end{equation*}
where, as before  $[\alpha]$ and $[\beta]$ denote the integer part of the parameters $\alpha$ and $\beta$. Thereafter, recall that $\alpha >2$  and $\beta \geq 1$. Then,  each expression above is controlled by the term $\ds{C\frac{e^{\eta_1\, t}}{|x|^{n+1}}}$ when $|x|\to +\infty$; and we thus obtain the wished estimate (\ref{Time-Der-Kernel}). \finpv 

\medskip  

Once we have the estimate (\ref{Time-Der-Kernel}), we get back to the estimate (\ref{J4}) to finally obtain 
\begin{equation}\label{Estim-J4}
J_4(t) \lesssim \frac{\gamma_3}{|x|^{n+1}} \int_{0}^{t}  e^{\eta_1\, \tau_1}\, |\tau |\, \| u(\tau,\cdot)\|^{2}_{\dot{H}^1} d \tau \lesssim \frac{\gamma_3C(t,u)}{|x|^{n+1}}, \qquad |x|\to +\infty. 
\end{equation}

With the estimates (\ref{eq15-bis}), (\ref{eq16-bis}),  (\ref{eq16-bis-bis}), and (\ref{Estim-J4}) at our disposal, we obtain the wished identity (\ref{gamma3eq18}). This identity together with the identity  (\ref{Split-Linear}) yield  the asymptotic profile  (\ref{eq18}). Theorem  \ref{Th-Asymptotic-Profile} is proven. \finpv

\subsection*{Proof of Corollary \ref{Corollary-1}}
By the asymptotic profile (\ref{Asymptotic}) and by the identity (\ref{Kernel-pointwise}), for $t>0$ fixed and for $|x|$ large enough we write 
\begin{equation}\label{Estimate-below}
\begin{split}
|u(t,x)| = &\, \left|K_{\alpha,\beta}(t,x) \left[ \int_{\R}u_0(y)dy+\gamma_3 \int_{0}^{t}\| u(\tau,\cdot)\|^{2}_{\dot{H}^1} d \tau \right] + R(t,x) \right| \\
\geq &\, |K_{\alpha,\beta}(t,x)| \left| \int_{\R}u_0(y)dy+\gamma_3 \int_{0}^{t}\| u(\tau,\cdot)\|^{2}_{\dot{H}^1} d \tau  \right|-|R(t,x)|\\
=&\, \frac{| I(t) |}{|x|^n}  \left| \int_{\R}u_0(y)dy+\gamma_3 \int_{0}^{t}\| u(\tau,\cdot)\|^{2}_{\dot{H}^1} d \tau  \right|-|R(t,x)|.
\end{split}
\end{equation}
Recall that $| R(t,x)|\leq \frac{c_2(t,u)}{|x|^{n+\varepsilon}}$ with $0<\varepsilon\leq 1$ (hence we have $| R(t,x)|=o (1/ |x|^n)$) and for the quantity $\frac{|I(t)|}{2}\left| \int_{\R}u_0(y)dy+\gamma_3 \int_{0}^{t}\| u(\tau,\cdot)\|^{2}_{\dot{H}^1} d \tau  \right|>0$ there exists $M>0$ such that for $|x|>M$ we have 
\begin{equation*}
 | R(t,x)|\leq  \frac{|I(t)|}{2} \left| \int_{\R}u_0(y)dy+\gamma_3 \int_{0}^{t}\| u(\tau,\cdot)\|^{2}_{\dot{H}^1} d \tau  \right| \frac{1}{|x|^n}.  
\end{equation*}
We get back to the previous estimate to obtain 
\begin{equation*}
 \frac{|I(t)|}{2} \left| \int_{\R}u_0(y)dy+\gamma_3 \int_{0}^{t}\| u(\tau,\cdot)\|^{2}_{\dot{H}^1} d \tau  \right| \frac{1}{|x|^n} \leq |u(t,x)|, \quad |x|\to +\infty,    
\end{equation*}
hence we set 
\begin{equation}\label{c3}
c_3(u_0,\gamma_3, t ,u)=  \frac{|I(t)|}{2} \left| \int_{\R}u_0(y)dy+\gamma_3 \int_{0}^{t}\| u(\tau,\cdot)\|^{2}_{\dot{H}^1} d \tau  \right|. 
\end{equation}
Corollary \ref{Corollary-1} is proven. \finpv 

\subsection*{Proof of Corollary \ref{Corollary-2}}

The proof follows very similar ideas of the previous one. Indeed, in the case $\gamma_3 =0$ and $\int_{\R}u_0(y)dy \neq 0$ by the estimate (\ref{Estimate-below}) and  for $|x|$ large enough we have
\begin{equation*}
  \frac{|I(t)|}{2}\left| \int_{\R}u_0(y)dy  \right|\frac{1}{|x|^n} \leq |u(t,x)|,  
\end{equation*}
where we set the quantity 
\begin{equation}\label{c4}
 c_4(u_0,t)=   \frac{|I(t)|}{2}\left| \int_{\R}u_0(y)dy  \right|. 
\end{equation}

On the other hand, in case $\gamma_3=0$ and $\int_{\R}u_0(y)dy=0$ by the identity (\ref{Asymptotic}) we obtain the estimate (\ref{Decay1Imp}). Corollary \ref{Corollary-2}  is proven. \finpv 

\begin{appendices}
\section*{Appendix}\label{AppendixA}
A proof of identity (\ref{Iden}). First,  we recall that the bilinear for $B(\cdot , \cdot)$ is given in (\ref{def-B}). Moreover, for the sake of simplicity, we shall write $K_{\alpha,\beta}(t,\cdot)\ast u_0= \tilde{u_0}(t,\cdot)$ and $K_{\alpha,\beta}(t,\cdot)\ast v_0= \tilde{v_0}(t,\cdot)$. Then, we have 
\begin{equation*}
\begin{split}
g(t,\xi) =&\, \int_{0}^{t} e^{-f(\xi)(t-\tau)}\, \mathcal{F}\left( B\big( \tilde{u_0}, \tilde{v_0} \big) \right)(\tau, \xi)\,  d \tau\\
=&\,  \int_{0}^{t} e^{-f(\xi)(t-\tau)}\, \left( \gamma_1\, i \xi (\widehat{\tilde{u_0}}\ast \widehat{\tilde{v_0}}) -\gamma_2\, \xi^2 (\widehat{\tilde{u_0}}\ast \widehat{\tilde{v_0}}) + \gamma_3 (i \xi\widehat{\tilde{u_0}}\ast i \xi  \widehat{\tilde{v_0}})  \right)  (\tau, \xi) d \tau \\
=&\,  \int_{0}^{t} e^{-f(\xi)(t-\tau)}\,  \left((\gamma_1\,  i\xi  - \gamma_2\, \xi^2)\, \int_{\R} e^{f(\xi-\eta) \tau} \widehat{u_0}(\xi-\eta)\,e^{f(\eta) \tau} \widehat{v_0}(\eta)  d \eta \right. \\
&\, \left.  -\gamma_3 \int_{\R} (\xi-\eta) e^{f(\xi-\eta) \tau} \widehat{u_0}(\xi-\eta)\, \eta\, e^{f(\eta) \tau} \widehat{v_0}(\eta)  d \eta   \right) d \tau \\
=&\, \int_{0}^{t} e^{-f(\xi)(t-\tau)}\, \int_{\R} \left[ \gamma_1\, i \xi -\gamma_2\,\xi^2-\gamma_3(\xi-\eta)\eta\right]\, e^{-f(\xi-\eta)\tau}\, e^{-f(\eta)\tau} \widehat{u_0}(\xi-\eta)\, \widehat{v_0}(\eta) d \eta\, d \tau\\
=& \, \int_{\R} \left[\gamma_1 i \xi -\gamma_2\,\xi^2-\gamma_3(\xi-\eta)\eta\right] \widehat{u_0}(\xi-\eta)\, \widehat{v_0}(\eta) \left( \int_{0}^{t} e^{-f(\xi)(t-\tau)}\, e^{-f(\xi-\eta)\tau}\, e^{-f(\eta) \tau} d\tau \right)d \eta. 
\end{split}
\end{equation*}
where, the integral in the time variable computes down as 
\[  \int_{0}^{t} e^{-f(\xi)(t-\tau)}\, e^{-f(\xi-\eta)\tau}\, e^{-f(\eta) \tau} d\tau= \frac{e^{-f(\eta)t-f(\xi-\eta)t}-e^{-f(\xi)t}}{f(\xi)-f(\eta)-f(\xi-\eta)}.\]
\end{appendices} 	
 
	  

\begin{thebibliography}{100}

\bibitem{Bao-2} Bao-Feng Feng, \& T. Kawahara. \emph{Multi-hump stationary waves for a Korteweg-deVries equation with
nonlocal perturbations}. Physica D. 137: 237-246 (2000).
\bibitem{Bao} Bao-Feng Feng, \& T. Kawahara. \emph{Temporal evolutions and stationary waves for dissipative Benjamin-Ono
equation}. Phys. D 139 pp. 301-318 (2000).

\bibitem{Biagioni} H. A. Biagioni, J. L. Bona, R. Iorio and M. Scialom. \emph{On the Korteweg-de Vries-Kuramoto-Sivashinsky equation}, Adv. Diff. Eq. 1:1-20 (1996).

\bibitem{Benjamin} T. B. Benjamin. \emph{Internal waves of permanent form in fluids of great depth}. J. Fluid Mech. 29, 559-592 (1967)

\bibitem{Chris-Ve} C.I. Christov \& M.G. Velarde. \emph{Dissipative solitons}. Physica D: Nonlinear Phenomena. 86 (32): 323–347 (1995). 
\bibitem{Chris-Ve-2} C.I.  Christov, \& M.G. Velarde. \emph{On localized solutions of an equation governing Bénard–Marangoni convection}. Appl. Math. Model. 17, 311–320 (1993).
\bibitem{Coclite} G.M. Coclite \& L. di Ruvo. \emph{Well-posedness results for the Kuramoto–Velarde equation}. Bollettino dell’Unione Matematica Italiana  14:659–679 (2021). 
\bibitem{CorJar0} M. F. Cortez \& O. Jarr\'in. \emph{On decay properties and asymptotic behavior of solutions to a non-local
perturbed KdV equation}. Nonlinear Analysis 187:365-396 (2019). 	
\bibitem{CorJar1} M.F. Cortez \& O. Jarr\'in. \emph{Spatial behavior of solutions for a large class of non-local PDE's arising from stratified flows}.   Differential and  Integral Equations 34(9/10): 539-594 (2021).
\bibitem{Fonseca} G. Fonseca, R. Pastr\'an and G. Rodr\'iguez-Blanco. \emph{The IVP for a nonlocal perturbation of the Benjamin-Ono equation in classical and weighted Sobolev spaces}. Journal of Mathematical Analysis and Applications, Volume 476, Issue 2: 391-425 (2019).  
\bibitem{Gar}
P.L. Garcia-Ybarra, J.L. Castillo \& M.G. Velarde. \emph{Benard–Marangoni convection with a deformable interface and poorly conducting boundaries}. Phys. Fluids, 30: 2655–2661 (1987).
\bibitem{Gar1}
P.L. Garcia-Ybarra, J.L. Castillo \& M.G. Velarde. \emph{A nonlinear evolution equation for Benard–Marangoni convection with deformable boundary}. Phys. Lett. A, 122: 107–110 (1987).
\bibitem{Hoop} 
 : A.P. Hooper, R. Grimshaw. \emph{Nonlinear instability at the interface between two viscous fluids}, Phys. Fluids 28 (1985) 37–45.
\bibitem{Hym}
J.M. Hyman \& B. Nicolaenko. \emph{Coherence and chaos in Kuramoto–Velarde equation}. In: Grandall, M.G., Rabinovitz, P.H., Turner, R.E.L. (eds.) Directions in Partial Differential Equations, pp. 89–111 (1987).
\bibitem{Kenig-Ponce-Vega} C. E. Kenig, G. Ponce and L. Vega. \emph{Small solutions to nonlinear Schr\"odinger equation}.  Ann.
Inst. H. Poincaré Anal. Non Linaire 10: 255-288 (1993).

\bibitem{KdV} D. J. Korteweg \& G. de Vries. \emph{On the change of form of long waves advancing in a rectangular canal, and on a new type of long stationary waves}. Philos. Mag. 39: 422–443 (1895).


\bibitem{LinaresPonce} F. Linares \& G. Ponce. \emph{Introduction to Nonlinear Dispersive Equations}. Springer Science \& Business Media (2009).

\bibitem{OST0}  L.A. Ostrovsky, Y.A. Stepanyams \& L.S. Tsimring. \emph{Nonlinear stage of the shearing instability in a stratified liquid of finite depth}. Fluid Dyn. 17: 540-546 (1983).

\bibitem{OST} L.A. Ostrovsky, Y.A. Stepanyams \& L.S. Tsimring. \emph{Radiation instability in a stratified shear flow}. Int. J. Nonlinear Mech 19: 151-161 (1984).

\bibitem{OST1} L.A. Ostrovsky, S.A. Rybak \& L.Sh. Tsimring. \emph{Negative energy waves in hydrodynamics}. Sov. Phys. Usp. 29:1040-1052 (1986).


\bibitem{Pap1}
D.T. Papageorgiou, Y. S. Smyrlis. \emph{The route to chaos for the Kuramoto–Sivashinsky equation}. Theor. Comp. Fluid Dyn., 3, (1), 15–42 (1991).

\bibitem{Pastran} R. Pastr\'an and O. Ria\~{n}o. \emph{Well-posedness for Fractional Growth-Dissipative Benjamin-Ono Equations}.  arXiv:1902.06868 (2019). 


\bibitem{Pilod} D. Pilod. \emph{Sharp well-posedness results for the Kuramoto-Velarde equation}. Communications on Pure and Applied Analysis, 7(4): 867-881 (2008).  


\bibitem{Pokhozhaev} S. I.  Pokhozhaev. \emph{On blow-up of solutions of the Kuramoto-Sivashinsky equation}.Russian Academy of Sciences Sbornik Mathematics 199(9):1355 (2008). 


\bibitem{Qian}  S. Qian, Y. C. Lee and H. H. Chen. \emph{A study of nonlinear dynamical models of plasma turbulence}, Phys. Fluids B1, 1:  87–98 (1989).
\bibitem{Norm}
C. Normand, M.G. Velarde. \emph{Convection}. Sci. Am, 243: 92–108 (1980).
\bibitem{Oer}
H.Oertel. Jr \& Zierep.  \emph{Convective Transport and Instability Phenomena}. Braun, Karlsruhe (1982).
\bibitem{Wang} H. Wang and \& A. Esfahani. \emph{Well-posedness results for the Ostrovsky, Stepanyams and Tsimring equation at the critical regularity}. Nonlinear Analysis: Real World Applications, Volume 44: 347-364 (2018).
\end{thebibliography}
\end{document}